\documentclass[a4paper,12pt]{article}
\usepackage{amsfonts}
\usepackage[xdvi]{graphicx}
\usepackage[dvips]{color}
\usepackage{amsmath,amssymb}

\newcommand{\bm}[1]{\mbox{\boldmath $#1$}}
\newcommand{\R}{{\mathbb R}}
\newcommand{\C}{{\mathbb C}}

\newcommand{\Z}{{\mathbb Z}}

\makeatletter

\@addtoreset{equation}{section}

\newcounter{def}[section]
\renewcommand{\thedef}{\stepcounter{def}\thesection.\@arabic\c@def }

\begin{document}
\setlength{\baselineskip}{24pt}
\begin{center}
\textbf{\LARGE{The first, second and fourth Painlev\'{e} equations on weighted projective spaces}}
\end{center}

\setlength{\baselineskip}{14pt}

\begin{center}
Institute of Mathematics for Industry, Kyushu University, Fukuoka,
819-0395, Japan

\large{Hayato CHIBA} \footnote{E mail address : chiba@imi.kyushu-u.ac.jp}
\end{center}
\begin{center}
Apr 16, 2014
\end{center}

\begin{center}
\textbf{Abstract}
\end{center}

The first, second and fourth Painlev\'{e} equations are studied by means of dynamical systems theory and
three dimensional weighted projective spaces $\C P^3(p,q,r,s)$ with suitable weights $(p,q,r,s)$
determined by the Newton diagrams of the equations or the versal deformations of vector fields.
Singular normal forms of the equations,
a simple proof of the Painlev\'{e} property and symplectic atlases of the spaces of 
initial conditions are given with the aid of the orbifold structure of $\C P^3(p,q,r,s)$.
In particular, for the first Painlev\'{e} equation, a well known Painlev\'{e}'s transformation 
is geometrically derived, which proves to be the Darboux coordinates of a certain algebraic surface with a 
holomorphic symplectic form.
The affine Weyl group, Dynkin diagram and the Boutroux coordinates are also studied from a view point 
of the weighted projective space.


% \tableofcontents

%%%%%%%%%%%%%%%%%%%%%%%%%%%%%%%%%%%%%%%%%%%%%%%%%%%%%%%%%%%%%%%%%%%%%%%%%%%%%%%%%%%%%%%%%%%%%%%%%%
%%%%%%%%%%%%%%%%%%%%%%%%%%%%%%%%%%%%%%%%%%%%%%%%%%%%%%%%%%%%%%%%%%%%%%%%%%%%%%%%%%%%%%%%%%%%%%%%%%

\section{Introduction}

The first, second and fourth Painlev\'{e} equations in Hamiltonian forms are given by
\begin{eqnarray}
(\text{P}_\text{I}) \left\{ \begin{array}{l}
\displaystyle \frac{dx}{dz} = 6y^2 + z \\[0.2cm]
\displaystyle \frac{dy}{dz} = x,  \\
\end{array} \right.
\label{P1}
\end{eqnarray}
\begin{eqnarray}
(\text{P}_\text{II}) \left\{ \begin{array}{l}
\displaystyle \frac{dx}{dz} = 2y^3 + yz +\alpha  \\[0.2cm]
\displaystyle \frac{dy}{dz} = x,  \\
\end{array} \right.
\label{P2}
\end{eqnarray}
\begin{eqnarray}
(\text{P}_\text{IV}) \left\{ \begin{array}{l}
\displaystyle \frac{dx}{dz} = -x^2 + 2xy + 2xz - 2\theta _{\infty}    \\[0.2cm]
\displaystyle \frac{dy}{dz} = -y^2+2xy-2yz-2\kappa_0,  \\
\end{array} \right.
\label{P4}
\end{eqnarray}
with Hamiltonian functions
\begin{eqnarray*}
& & H_\text{I} = \frac{1}{2}x^2 - 2y^3 - zy, \\
& & H_\text{II} = \frac{1}{2}x^2 - \frac{1}{2}y^4 - \frac{1}{2}zy^2 - \alpha y,  \\
& &  H_\text{IV} = -xy^2 + x^2y - 2xyz - 2\kappa_0 x + 2\theta _\infty y, 
\end{eqnarray*}
where $\alpha , \theta _\infty$ and $\kappa_0 \in \C$ are parameters.
These equations are investigated by means of the weighted projective spaces $\C P^3(p,q,r,s)$
with natural numbers $p,q,r,s$ given by
\begin{eqnarray*}
(\text{P}_\text{I}) & & (p,q,r,s) = (3,2,4,5), \\
(\text{P}_\text{II}) & & (p,q,r,s) = (2,1,2,3), \\
(\text{P}_\text{IV}) & & (p,q,r,s) = (1,1,1,2).
\end{eqnarray*}
These numbers will be determined by the Newton diagrams of the equations or the versal deformations
of a certain class of dynamical systems.
The weighted projective space $\C P^3(p,q,r,s)$ is a three dimensional compact orbifold (toric variety)
with singularities, see Sec.2 for the definition.

($\text{P}_\text{I}$), ($\text{P}_\text{II}$) and ($\text{P}_\text{IV}$) are given as 
differential equations on $\C P^3(p,q,r,s)$, which is regarded as a compactification of 
the original phase space $\C^3_{(x,y,z)}$ of the Painlev\'{e} equations.
The Painlev\'{e} equations are invariant under the $\Z_s$ action of the form
\begin{equation}
(x,y,z) \mapsto (\omega ^px, \omega ^qy, \omega ^rz), \quad \omega := e^{2\pi i/s},
\label{1-4}
\end{equation}
with $p,q,r,s$ as above.
As a result, it turns out that ($\text{P}_\text{I}$), ($\text{P}_\text{II}$) and ($\text{P}_\text{IV}$)
are well defined as meromorphic differential equations on $\C P^3(p,q,r,s)$.

The space $\C P^3(p,q,r,s)$ is decomposed as
\begin{equation}
\C P^3(p,q,r,s) = \C^3/\Z_s \,\, \cup \,\, \C P^2(p,q,r), \quad (\text{disjoint}).
\label{1-5}
\end{equation}
This means that $\C P^3(p,q,r,s)$ is a compactification of $\C^3/\Z_s$ obtained by attaching a 2-dim weighted
projective space $\C P^2(p,q,r)$ at infinity.
The Painlev\'{e} equations ($\text{P}_\text{J}$),\, ($\text{J} = \text{I}, \text{II},\text{IV}$) divided by the 
$\Z_s$ action are given on $\C^3/\Z_s$, and the 2-dim space 
$\C P^2(p,q,r)$ describes the behavior of ($\text{P}_\text{J}$) near infinity 
(i.e. $x=\infty$ or $y=\infty$ or $z=\infty$).
On the ``infinity set" $\C P^2(p,q,r)$, there exist several singularities of the foliation
defined by solutions of the equation.
Some of them correspond to movable poles of ($\text{P}_\text{J}$),
and the others correspond to the irregular singular point $z=\infty$.
Local properties of these singularities of the foliation will be investigated by means of dynamical systems theory.
Our main results include  
\begin{itemize}

 \item the fact that the Painlev\'{e} equations are locally transformed into integrable systems near 
movable singularities,

 \item a simple proof of the fact that any solutions of ($\text{P}_\text{J}$) are meromorphic on $\C$,

 \item a simple construction of the symplectic atlas of Okamoto's space of initial conditions,
 
\item for ($\text{P}_\text{I}$), a geometric interpretation of Painlev\'{e}'s coordinates defined by
\begin{equation}
\left\{ \begin{array}{l}
x = uw^3 - 2w^{-3} - \frac{1}{2}z w - \frac{1}{2} w^2 \\
y = w^{-2},  \\
\end{array} \right.
\label{1-6}
\end{equation}
which was introduced in his original work \cite{Pai}
to prove the Painlev\'{e} property of ($\text{P}_\text{I}$),

\item a geometric interpretation of Boutroux's coordinates introduced in \cite{Bou} to 
investigate the irregular singular point of ($\text{P}_\text{I}$) and ($\text{P}_\text{II}$).
 \end{itemize}

In Sec.2, the Newton diagram of the Painlev\'{e} equation will be introduced 
to find a suitable weight of the weighted projective space $\C P^3(p,q,r,s)$.
Furthermore, it is shown that the Painlev\'{e} equations are obtained from certain problems
of dynamical systems theory.
Such a relationship between the Painlev\'{e} equations and dynamical systems proposes
normal forms of the Painlev\'{e} equations because for dynamical systems (germs of vector fields),
the normal form theory have been well developed.

In Sec.3, with the aid of the orbifold structure of $\C P^3(p,q,r,s)$ and 
the Poincar\'{e} linearization theorem,
it will be shown that ($\text{P}_\text{I}$), ($\text{P}_\text{II}$) and ($\text{P}_\text{IV}$) are 
locally transformed into integrable systems near each movable singularities.
For example, ($\text{P}_\text{I}$) and ($\text{P}_\text{II}$) can be transformed into the 
equations $y'' = 6y^2$ and $y'' = 2y^3$, respectively.
See Sec.3 for the result for ($\text{P}_\text{IV}$).
This fact was first obtained by \cite{CosCos} for ($\text{P}_\text{I}$), in which the transformed 
equation $y'' = 6y^2$ is called the singular normal form.
Our proof is based on the Poincar\'{e} linearization theorem and it is easily 
applied to other Painlev\'{e} equations, including ($\text{P}_\text{III}$), ($\text{P}_\text{V}$)
($\text{P}_\text{VI}$) and higher order Painlev\'{e} equations.
By using this result, a simple proof of the Painlev\'{e} property is proposed;
that is, a new proof of the fact that any solutions of ($\text{P}_\text{I}$), 
($\text{P}_\text{II}$) and ($\text{P}_\text{IV}$) are meromorphic on $\C$ will be given.

In Sec.4, the weighted blow-up will be introduced to construct the spaces of initial conditions.
For a polynomial system, a manifold $E(z)$ parameterized by $z\in \C$ is called the space of initial conditions
if any solutions give global holomorphic sections on the fiber bundle $\mathcal{P} = \{ (x,z) \, | \, x\in E(z) , z\in \C\}$ over $\C$.
It is remarkable that only one, two and three times blow-ups are sufficient to obtain 
the spaces of initial conditions for ($\text{P}_\text{I}$), ($\text{P}_\text{II}$) and ($\text{P}_\text{IV}$),
respectively, if we use suitable weights, while Okamoto performed blow-ups (without weights) eight times to
obtain the space of initial conditions \cite{Oka}.
Further, our method easily provides a symplectic atlas of the space of initial conditions.
Then, each Painlev\'{e} equation is characterized as a unique Hamiltonian system on the
space of initial conditions admitting a holomorphic symplectic form.
Symplectic atlases of the spaces of initial conditions were first obtained by 
Takano \textit{et al.} \cite{Tak1, Tak2, Tak3} only for ($\text{P}_\text{II}$) to ($\text{P}_\text{VI}$),
while left open for ($\text{P}_\text{I}$).
In the present paper, the orbifold structure plays an important role to obtain
a symplectic atlas for ($\text{P}_\text{I}$), see also Iwasaki and Okada \cite{Iwa}
for the orbifold setting of ($\text{P}_\text{I}$).  

By the weighted blow-up of $\C P^3(3,2,4,5)$ for ($\text{P}_\text{I}$),
we will recover the famous Painlev\'{e}'s coordinates (\ref{1-6})
in a purely geometric manner.
Painlev\'{e} found the coordinate transformation (\ref{1-6}) in an analytic way to prove the Painlev\'{e} property of 
($\text{P}_\text{I}$) (see \cite{Gro}).
From our approach based on the weighted projective space, Painlev\'{e}'s coordinates prove to be nothing 
but the Darboux coordinates of the nonsingular algebraic surface $M(z)$ defined by
\begin{eqnarray*}
V^2 =  UW^4 + 2z W^3 + 4W,
\end{eqnarray*}
which admits a holomorphic symplectic form, where $z \in \C$ is an independent variable of ($\text{P}_\text{I}$) 
and it is a parameter of the surface.
Our space of initial conditions is obtained by glueing $\C^2_{(x,y)}$ (the original space for dependent variables)
and the surface $M(z)$ by a symplectic mapping.
Then, ($\text{P}_\text{I}$) is a Hamiltonian system with respect to the symplectic form.
Since (\ref{1-6}) is a one-to-two transformation, an orbifold setting is essential to give a 
geometric meaning to Painlev\'{e}'s coordinates;
the orbifold $\C P^3 (3,2,4,5)$ provides a natural $\Z_2$-action which makes (\ref{1-6}) a one-to-one transformation.

In Sec.5, the characteristic indices for ($\text{P}_\text{I}$), ($\text{P}_\text{II}$) and ($\text{P}_\text{IV}$)
will be defined.
A few simple properties such as a relation with the Kovalevskaya exponents and the weights of $\C P^3 (p,q,r,s)$ will be given. 

In Sec.6, the Boutroux coordinates will be introduced.
It is shown that the weighted blow-ups of $\C P^3(p,q,r,s)$ constructed in Sec.4 also includes the space of 
initial conditions written in the Boutroux coordinates.
Further, we will show that autonomous Hamiltonian systems are embedded in the Boutroux coordinates.

In Sec.7, the extended affine Weyl group for ($\text{P}_\text{II}$) and ($\text{P}_\text{IV}$)
will be considered.
The action of the group on the original chart $\C^3_{(x,y,z)}$ is extended to a birational 
transformation on $\C P^3(p,q,r,s)$.
It is proved that on the ``infinity set", $\C P^2(p,q,r)$, the foliation defined by
an autonomous Hamiltonian system is invariant under the 
automorphism group $\mathrm{Aut} (X)$, where $X = A^{(1)}_1$ for ($\text{P}_\text{II}$)
and $X = A^{(1)}_2$ for ($\text{P}_\text{IV}$).

In Sec.8, a cellular decomposition of the weighted blow-ups of $\C P^3(p,q,r,s)$ will be given.
We will show that the weighted blow-ups of $\C P^3(p,q,r,s)$ is naturally decomposed into
the fiber space for ($\text{P}_\text{J}$) (a fiber bundle over $\C$ whose fiber is the space of 
initial conditions), a certain elliptic fibration over the moduli space of complex tori,
and the projective curve $\C P^1$.
We also show that the extended Dynkin diagrams of type $\tilde{E}_8, \tilde{E}_7$ and $\tilde{E}_6$ are hidden 
in the weighted blow-ups of $\C P^3(p,q,r,s)$.

An approach using toric varieties is also applicable to the third, fifth, sixth Painlev\'{e} equations
and higher order Painlev\'{e} equations, which will appear in a forthcoming paper.

%%%%%%%%%%%%%%%%%%%%%%%%%%%%%%%%%%%%%%%%%%%%%%%%%%%%%%%%%%%%%%%%%%%%%%%%%%%%%%%%%%%%%%%%%%%%%%%%%%
%%%%%%%%%%%%%%%%%%%%%%%%%%%%%%%%%%%%%%%%%%%%%%%%%%%%%%%%%%%%%%%%%%%%%%%%%%%%%%%%%%%%%%%%%%%%%%%%%%

\section{Weighted projective spaces}

In this section, a weighted projective space $\C P^3(p,q,r,s)$ is defined and the first, second and fourth Painlev\'{e} equations
are given as meromorphic equations on $\C P^3(p,q,r,s)$ for suitable integers $p,q,r,s$.
Such integers $p,q,r,s$ will be found via the Newton diagrams of the equations.
We also give a relationship between the Painlev\'{e} equations and the normal form theory of dynamical systems,
which proposes normal forms of the Painlev\'{e} equations.

%%%%%%%%%%%%%%%%%%%%%%%%%%%%%%%%%%%%%%%%%%%%%%%%%%%%%%%%%%%%%%%%%%%%%%%%%%%%%%%%%%%%%%%%%%%%%%%%%%

\subsection{Newton diagram}

Let us consider the system of polynomial differential equations 
\begin{eqnarray}
\frac{dx}{dz} = f_1(x,y,z), \quad \frac{dy}{dz} = f_2(x,y,z).
\label{2-1}
\end{eqnarray}
The exponent of a monomial $x^iy^jz^k$ included in $f_1$ is defined by $(i-1,j,k+1)$,
and by $(i,j-1,k+1)$ for one in $f_2$.
Each exponent specifies a point of the integer lattice in $\R^3$.
The Newton polyhedron of the system (\ref{2-1}) is the convex hull of the union of the positive quadrants
$\R_+^3$ with vertices at the exponents of the monomials which appear in the system.
The Newton diagram of the system is the union of the compact faces of its Newton polyhedron.
Suppose that the Newton diagram consists of only one compact face.
Then, there is a tuple of positive integers $(p_1, p_2, r, s)$ such that
the compact face lies on the plane $p_1 x + p_2 y + rz = s$ in $\R^3$.
In this case, the function $f_i\, (i=1,2)$ satisfies 
\begin{eqnarray*}
f_i(\lambda ^{p_1}x, \lambda ^{p_2}y , \lambda ^r z)
 = \lambda ^{s-r+p_i} f_i (x_1, \cdots ,x_m,z), 
\end{eqnarray*}
for any $\lambda \in \C$.

We also consider the perturbation of the system (\ref{2-1}) of the form
\begin{eqnarray}
\frac{dx}{dz} = f_1(x,y,z)+g_1(x,y,z), \quad \frac{dy}{dz} = f_2(x,y,z)+g_2(x,y,z).
\label{2-1b}
\end{eqnarray}
Suppose that $g_i(\lambda ^{p_1}x , \lambda ^{p_2}y , \lambda ^{r}z) \sim o(\lambda ^{s-r+p_i})$ for $i=1,2$ as $\lambda \to \infty$.
This implies that exponents of any monomials included in $g_i$ lie on the lower side of the plane $p_1 x + p_2 y + rz = s$.

The Newton polyhedron of the first Painlev\'{e} equation (\ref{P1}) is defined by three points
$(-1,2,1), (-1,0,2)$ and $(1,-1,1)$.
Hence, the Newton diagram consists of the unique face which lies on the plane $3x+2y + 4z = 5$.
One of the normal vector to the plane is given by $\mathbf{e}_0 = (-3/5, -2/5, -4/5)$.
Put $\mathbf{e}_1 = (1,0,0), \mathbf{e}_2 = (0,1,0)$ and $\mathbf{e}_3 = (0,0,1)$.
Then, the toric variety defined by the fan made up of the cones generated by all proper subsets
of $\{ \mathbf{e}_0, \mathbf{e}_1,\mathbf{e}_2,\mathbf{e}_3\}$ is the weighted projective space
$\C P^3(3,2,4,5)$ \cite{Cox}.

Next, let us consider the second Painlev\'{e} equation (\ref{P2}) with 
$f = (2y^3 + yz, x)$ and $g= (\alpha ,0)$.
The Newton polyhedron of $f = (f_1, f_2)$ is defined by three points $(-1,3,1), (-1,1,2), (1,-1,1)$,
and the Newton diagram is given by the unique face on the plane $2x + y+ 2z = 3$.
The associated toric variety is $\C P^3(2,1,2,3)$.

For the fourth Painlev\'{e} equation (\ref{P4}),
put $f = (-x^2 + 2xy + 2xz, -y^2+2xy-2yz)$ and $g= (- 2\theta _{\infty} ,-2\kappa_0)$.
The Newton diagram of $f = (f_1, f_2)$ is given by the unique face on the plane $x + y+ z = 2$ 
passing through the exponents $(1,0,1), (0,1,1)$ and $(0,0,2)$. 
The associated toric variety is $\C P^3(1,1,1,2)$.

In what follows, the weights $(p,q,r,s)$ denote $(3,2,4,5), (2,1,2,3)$ and $(1,1,1,2)$, respectively,
for ($\text{P}_\text{I}$), ($\text{P}_\text{II}$) and ($\text{P}_\text{IV}$).

The weighted degree of a monomial $x^iy^jz^k$ with respect to the weight $(p,q,r)$ is defined by
$\mathrm{deg} (x^iy^jz^k) = pi+qj+rk$.
The weighted degree of a polynomial $f = \sum a_{ijk} x^iy^jz^k$ is defined by
\begin{eqnarray*}
\mathrm{deg} (f) = \max_{i,j,k} \{ \mathrm{deg} (x^iy^jz^k) \, | \, a_{ijk} \neq  0\}.
\end{eqnarray*}
For ($\text{P}_\text{I}$), ($\text{P}_\text{II}$) and ($\text{P}_\text{IV}$) with the weights
$(p,q,r) = (3,2,4), (2,1,2)$ and $(1,1,1)$, respectively,
the weighted degrees of the Hamiltonian functions are
\begin{eqnarray*}
\mathrm{deg}(H_\text{I}) = 6, \quad \mathrm{deg}(H_\text{II}) = 4, \quad \mathrm{deg}(H_\text{IV}) = 3.
\end{eqnarray*}
They satisfy $\mathrm{deg}(H_\text{J}) = s + 1\, (\text{J} = \text{I},\text{II},\text{IV})$.
Further, it will be shown that they coincide with the Kovalevskaya exponents (Sec.2.3) 
and the characteristic index $\lambda _1$ (Sec.5).

%%%%%%%%%%%%%%%%%%%%%%%%%%%%%%%%%%%%%%%%%%%%%%%%%%%%%%%%%%%%%%%%%%%%%%%%%%%%%%%%%%%%%%%%%%%%%%%%%%

\subsection{Weighted projective space}

Let $\widetilde{U}$ be a complex manifold and $\Gamma $ a finite group acting analytically and effectively on 
$\widetilde{U}$.
In general, the quotient space $\widetilde{U}/\Gamma $ is not a smooth manifold if the action has fixed points.
Roughly speaking, a (complex) orbifold $M$ is defined by glueing a family of such spaces 
$\widetilde{U}_\alpha /\Gamma _\alpha $;
a Hausdorff space $M$ is called an orbifold if there exist an open covering $\{ U_\alpha  \}$ of $M$ and 
homeomorphisms $\varphi _\alpha : U_\alpha \simeq \widetilde{U}_\alpha /\Gamma _\alpha $.
See \cite{Thu} for more details.
In this article, only quotient spaces of the form $\C^n/ \Z_p$ will be used.

Consider the weighted $\C^*$-action on $\C^{4}$ defined by
\begin{equation}
(x,y,z,\varepsilon ) \mapsto (\lambda ^{p}x , \lambda ^{q} y, \lambda ^rz, \lambda ^s \varepsilon ), \quad
\lambda \in \C^* := \C\backslash \{ 0\},
\end{equation}
where the weights $p,q,r,s$ are positive integers.
We assume $1\leq p,q,r \leq s$ without loss of generality.
Further, we suppose that any three numbers among them have no common divisors.
The quotient space
\begin{eqnarray*}
\C P^3(p,q,r,s) := \C^{4}\backslash \{ 0\} /\C^*
\end{eqnarray*}
gives a three dimensional orbifold called the weighted projective space.
The inhomogeneous coordinates of $\C P^3(p,q,r,s)$, which give an orbifold structure of $\C P^3(p,q,r,s)$,
are defined as follows.

The space $\C P^3(p,q,r,s)$ is defined by the equivalence relation on $\C^4\backslash \{ 0\}$
\begin{eqnarray*}
(x,y,z, \varepsilon ) \sim (\lambda^p x, \lambda ^qy, \lambda ^rz, \lambda ^s \varepsilon ).
\end{eqnarray*}
(i) When $x\neq 0$,
\begin{eqnarray*}
(x,y,z,\varepsilon ) \sim (1,\,\,  x^{-q/p}y, \,\, x^{-r/p}z,\,\, x^{-s/p}\varepsilon )
=:(1, Y_1, Z_1, \varepsilon _1).
\end{eqnarray*}
Due to the choice of the branch of $x^{1/p}$, we also obtain
\begin{eqnarray*}
(Y_1, Z_1, \varepsilon _1) \sim
(e^{-2q\pi i /p}Y_1, e^{-2r\pi i /p}Z_1, e^{-2s\pi i /p}\varepsilon _1),
\end{eqnarray*}
by putting $x \mapsto e^{2\pi i }x$.
This implies that the subset of $\C P^3(p,q,r,s)$ such that $x\neq 0$ is homeomorphic to $\C^3 / \Z_p$,
where the $\Z_p$-action is defined as above.

(ii) When $y\neq 0$, 
\begin{eqnarray*}
(x,y,z,\varepsilon ) \sim (y^{-p/q}x,\,\, 1 ,\,\, y^{-r/q}z,\,\, y^{-s/q}\varepsilon )
=:(X_2, 1, Z_2, \varepsilon _2).
\end{eqnarray*}
Because of the choice of the branch of $y^{1/q}$, we obtain
\begin{eqnarray*}
(X_2, Z_2, \varepsilon _2) \sim (e^{-2p\pi i/q}X_2, e^{-2r\pi i/q}Z_2, e^{-2s\pi i/q} \varepsilon _2).
\end{eqnarray*}
Hence, the subset of $\C P^3(p,q,r,s)$ with $y\neq 0$ is homeomorphic to $\C^3 / \Z_q$.

(iii) When $z\neq 0$,
\begin{eqnarray*}
(x,y,z,\varepsilon ) \sim (z^{-p/r}x,\,\, z^{-q/r}y ,\,\, 1,\,\, z^{-s/r}\varepsilon )
=: (X_3, Y_3, 1, \varepsilon _3).
\end{eqnarray*}
Similarly, the subset $\{ z \neq 0\} \subset \C P^3(p,q,r,s)$ is homeomorphic to $\C^3 / \Z_r$.

(iv) When $\varepsilon \neq 0$,
\begin{eqnarray*}
(x,y,z,\varepsilon ) \sim (\varepsilon ^{-p/s}x,\,\, \varepsilon ^{-q/s}y ,\,\, \varepsilon ^{-r/s}z ,\,\, 1)
=:(X_4, Y_4, Z_4, 1).
\end{eqnarray*}
The subset $\{ \varepsilon  \neq 0\} \subset \C P^3(p,q,r,s)$ is homeomorphic to $\C^3 / \Z_s$.

This proves that the orbifold structure of $\C P^3(p,q,r,s)$ is given by
\begin{eqnarray*}
\C P^3(p,q,r,s) = \C^3/\Z_p \,\, \cup \,\, \C^3/\Z_q \,\, \cup \,\, \C^3/\Z_r \,\, \cup \,\, \C^3/\Z_s.
\end{eqnarray*}
The local charts $(Y_1, Z_1, \varepsilon _1)$, $(X_2, Z_2, \varepsilon _2)$,
$(X_3, Y_3, \varepsilon _3)$ and $(X_4, Y_4, Z_4)$ defined above are called inhomogeneous coordinates
as the usual projective space.
Note that they give coordinates on the lift $\C^3$, not on the quotient $\C^3 / \Z_i\,\, (i=p,q,r,s)$.
Therefore, any equations written in these inhomogeneous coordinates should be invariant under
the corresponding $\Z_i$ actions.

In what follows, we use the notation $(x,y,z)$ for the fourth local chart instead of $(X_4, Y_4, Z_4)$
because the Painlev\'{e} equation will be given on this chart.

The transformations between inhomogeneous coordinates are give by
\begin{eqnarray}
\left\{ \begin{array}{rrrr}
x = & \varepsilon _1^{-p/s} = & X_2\varepsilon _2^{-p/s} =& X_3\varepsilon _3^{-p/s} \\
y = & Y_1\varepsilon _1^{-q/s} = & \varepsilon _2^{-q/s}=& Y_3\varepsilon _3^{-q/s}\\
z = & Z_1\varepsilon_1 ^{-r/s} = & Z_2\varepsilon _2^{-r/s}=& \varepsilon _3^{-r/s}.
\end{array} \right.
\label{2-3}
\end{eqnarray}

We give the differential equation defined on the $(x, y, z)$-coordinates as
\begin{eqnarray}
\frac{dx}{dz} =f(x, y, z), \quad \frac{dy}{dz} = g(x, y, z),
\label{2-4}
\end{eqnarray}
where $f$ and $g$ are rational functions.
By the transformation (\ref{2-3}), the above equation is rewritten as equations on the other
inhomogeneous coordinates 
$(Y_1, Z_1, \varepsilon _1)$, $(X_2, Z_2, \varepsilon _2)$ and $(X_3, Y_3, \varepsilon _3)$.
It is easy to verify that the equations written in the other inhomogeneous coordinates are rational
if and only if Eq.(\ref{2-4}) is invariant under the $\Z_s$-action
\begin{eqnarray}
(x, y, z) \mapsto (\omega ^px, \omega ^qy, \omega ^r z), \quad \omega = e^{2\pi i/s}.
\label{2-5}
\end{eqnarray}
In this case, the equations written in 
$(Y_1, Z_1, \varepsilon _1)$, $(X_2, Z_2, \varepsilon _2)$ and $(X_3, Y_3, \varepsilon _3)$ are 
invariant under the $\Z_p, \Z_q$ and $\Z_r$-actions, respectively.
Hence, a tuple of these four equations gives a well-defined rational differential equation on $\C P^3(p,q,r,s)$.

When $x = \infty$ or $y = \infty$ or $z = \infty$, we have $\varepsilon _1=0$ or 
$\varepsilon _2 = 0$ or $\varepsilon _3 = 0$.
In this case, the transformation (\ref{2-3}) results in
\begin{eqnarray}
\left\{ \begin{array}{ll}
Y_1 = X_2^{-q/p}  &  = Y_3 X_3^{-q/p}, \\
Z_1 = Z_2 X_2^{-r/p} & = X_3^{-r/p}. \\
\end{array} \right.
\label{2-6}
\end{eqnarray}
The space obtained by glueing three copies of $\C^2$ by the above relations gives the 
2-dim weighted projective space $\C P^2(p,q,r)$.
Thus, we have obtained the decomposition 
\begin{equation}
\C P^3(p,q,r,s) = \C^3/\Z_s \,\, \cup \,\, \C P^2(p,q,r), \quad (\text{disjoint}).
\label{2-7}
\end{equation}
On the covering space $\C^3$ of $\C^3 /\Z_s$, the coordinates $(x, y, z)$ is assigned and 
Eq.(\ref{2-4}) is given.
The equation on $\C P^2(p,q,r)$ is obtained by putting $\varepsilon _1 = 0$ or $\varepsilon _2 =0$
or $\varepsilon _3 = 0$, which describes the behavior of Eq.(\ref{2-4}) near infinity;
\begin{equation}
\C P^2 (p,q,r) = \{ \varepsilon _1 = 0\} \cup \{ \varepsilon _2 = 0\} \cup \{ \varepsilon _3 = 0\}.
\label{2-8}
\end{equation}

Now we give the first Painlev\'{e} equation (\ref{P1}) on the fourth local chart of $\C P^3(3,2,4,5)$.
By (\ref{2-4}), ($\text{P}_\text{I}$) is transformed into the following equations 
\begin{eqnarray}
& & \frac{dY_1}{d\varepsilon _1} = \frac{3 - 12Y_1^3 - 2Y_1Z_1}{\varepsilon _1 (-30 Y_1^2 - 5Z_1)},
\quad \frac{dZ_1}{d\varepsilon _1} = \frac{3\varepsilon _1 - 24Y_1^2 Z_1 - 4 Z_1^2}{\varepsilon _1 (-30 Y_1^2 - 5Z_1)}, \label{P11} \\[0.2cm]
& & \frac{dX_2}{d\varepsilon _2} = \frac{-12 - 2Z_2 + 3X_2^2}{5X_2\varepsilon _2}, \quad
\frac{dZ_2}{d\varepsilon _2} = \frac{-2\varepsilon _2  + 4 X_2Z_2}{5X_2\varepsilon _2}, \label{P12} \\[0.2cm]
& & \frac{dX_3}{d\varepsilon _3} = \frac{24Y_3^2 + 4 - 3X_3\varepsilon _3}{-5\varepsilon _3^2}, \quad
\frac{dY_3}{d\varepsilon _3} = \frac{4X_3 - 2Y_3\varepsilon _3}{-5\varepsilon _3^2}, \label{P13}
\end{eqnarray}
on the other inhomogeneous coordinates.
Although the transformations (\ref{2-3}) have branches, the above equations are rational
due to the symmetry (\ref{2-5}) of ($\text{P}_\text{I}$).
Hence, they define a rational ODE on $\C P^3(3,2,4,5)$ in the sense of an orbifold.

Next, we give the second Painlev\'{e} equation (\ref{P2}) on the fourth local chart of $\C P^3(2,1,2,3)$.
By (\ref{2-3}), ($\text{P}_\text{II}$) is transformed into the following equations
\begin{eqnarray}
& & \frac{dY_1}{d\varepsilon _1} 
= \frac{-2 + Y_1(2Y_1^3 + Y_1Z_1 + \alpha \varepsilon _1)}{3\varepsilon _1 (2Y_1^3 + Y_1Z_1 + \alpha \varepsilon _1)},
 \frac{dZ_1}{d\varepsilon _1} 
= \frac{-2\varepsilon _1 + 2Z_1(2Y_1^3 + Y_1Z_1 + \alpha \varepsilon _1)}{3\varepsilon _1 (2Y_1^3 + Y_1Z_1 + \alpha \varepsilon _1)},
\quad \quad \quad \label{P21} \\[0.2cm]
& & \frac{dX_2}{d\varepsilon _2} = \frac{2X_2^2 - (2+Z_2 + \alpha \varepsilon _2)}{3X_2\varepsilon _2}, \quad
\frac{dZ_2}{d\varepsilon _2} = \frac{2X_2Z_2 - \varepsilon _2}{3X_2\varepsilon _2}, \label{P22} \\[0.2cm]
& & \frac{dX_3}{d\varepsilon _3} = \frac{4Y_3^3 + 2Y_3 + 2\alpha \varepsilon _3 - 2X_3\varepsilon _3}{-3\varepsilon _3^2}, \quad
\frac{dY_3}{d\varepsilon _3} = \frac{2X_3 - Y_3\varepsilon _3}{-3\varepsilon _3^2} \label{P23},
\end{eqnarray}
on the other local charts.
They define a rational ODE on $\C P^3(2,1,2,3)$ because of the symmetry (\ref{2-5}) of ($\text{P}_\text{II}$).
 
Similarly, we give the fourth Painlev\'{e} equation (\ref{P4}) on the fourth local chart of $\C P^3(1,1,1,2)$.
The equations written in the other inhomogeneous coordinates are given by
\begin{eqnarray}
\left\{ \begin{array}{ll}
\displaystyle \frac{dY_1}{d\varepsilon _1} 
= \frac{-Y_1^2+2Y_1-2Y_1Z_1-2\kappa_0\varepsilon_1+Y_1(1-2Y_1-2Z_1+2\theta_\infty\varepsilon_1)}{2\varepsilon_1 (1-2Y_1-2Z_1+2\theta_\infty\varepsilon_1)},\\[0.3cm]
\displaystyle \frac{dZ_1}{d\varepsilon _1} 
= \frac{\varepsilon _1+Z_1(1-2Y_1-2Z_1+2\theta _\infty\varepsilon _1)}{2\varepsilon _1 (1-2Y_1-2Z_1+2\theta _\infty\varepsilon _1)}, 
\end{array} \right.
\label{P41}
\end{eqnarray}
\begin{eqnarray}
\left\{ \begin{array}{ll}
\displaystyle  \frac{dX_2}{d\varepsilon _2} 
= \frac{ -X_2^2 + 2X_2 + 2X_2Z_2-2\theta _\infty \varepsilon _2+X_2(1-2X_2+2Z_2+2\kappa_0\varepsilon _2)}{2\varepsilon _2 (1-2X_2+2Z_2+2\kappa_0\varepsilon _2)},\\[0.3cm]
\displaystyle \frac{dZ_2}{d\varepsilon _2}
 = \frac{\varepsilon _2 + Z_2 (1-2X_2+2Z_2+2\kappa_0\varepsilon _2)}{2\varepsilon _2 (1-2X_2+2Z_2+2\kappa_0\varepsilon _2)}, 
\end{array} \right.
\label{P42}
\end{eqnarray}
\begin{eqnarray}
\left\{ \begin{array}{ll}
\displaystyle  \frac{dX_3}{d\varepsilon _3} 
= \frac{-X_3^2+2X_3Y_3+2X_3-2\theta _\infty \varepsilon _3-X_3\varepsilon _3}{-2\varepsilon _3^2},\\[0.3cm]
\displaystyle \frac{dY_3}{d\varepsilon _3} = \frac{-Y_3^2+2X_3Y_3-2Y_3-2\kappa_0 \varepsilon _3-Y_3\varepsilon _3}{-2\varepsilon _3^2}.
\end{array} \right.
\label{P43}
\end{eqnarray}
They define a rational ODE on $\C P^3(1,1,1,2)$.

%%%%%%%%%%%%%%%%%%%%%%%%%%%%%%%%%%%%%%%%%%%%%%%%%%%%%%%%%%%%%%%%%%%%%%%%%%%%%%%%%%%%%%%%%%%%%%%%%%

\subsection{Laurent series of solutions}

Before starting the analysis of the Painlev\'{e} equations by using the weighted projective spaces,
it is convenient to write down Laurent series of solutions.
Since any solutions of ($\text{P}_\text{I}$), ($\text{P}_\text{II}$) and ($\text{P}_\text{IV}$) 
are meromorphic, a general solution admits the Laurent series with respect to $T:=z-z_0$, where $z_0$ is a movable pole.

For the first Painlev\'{e} equation ($\text{P}_\text{I}$), the Laurent series of a general solution is given by
\begin{equation}
\left(
\begin{array}{@{\,}c@{\,}}
x\\
y
\end{array}
\right) = \left(
\begin{array}{@{\,}c@{\,}}
\! -2 \! \\
0
\end{array}
\right) T^{-3} + 
\left(
\begin{array}{@{\,}c@{\,}}
0 \\
1
\end{array}
\right) T^{-2} - 
\left(
\begin{array}{@{\,}c@{\,}}
\! z_0/5 \!\\
0
\end{array}
\right) T - 
\left(
\begin{array}{@{\,}c@{\,}}
1/2 \\
\! z_0/10 \!
\end{array}
\right) T^2 + 
\left(
\begin{array}{@{\,}c@{\,}}
A_6 \\
-1/6 \!
\end{array}
\right) T^3 + \cdots ,
\end{equation} 
where $A_6$ is an arbitrary constant.

For the second Painlev\'{e} equation ($\text{P}_\text{II}$), the Laurent series are expressed 
in two ways as
\begin{eqnarray*}
& (i) &  \!\! \left(
\begin{array}{@{\,}c@{\,}}
x\\
y
\end{array}
\right)\! =\! \left(
\begin{array}{@{\,}c@{\,}}
1 \\
0
\end{array}
\right) T^{-2} - 
\left(
\begin{array}{@{\,}c@{\,}}
0 \\
1
\end{array}
\right) T^{-1} +
\left(
\begin{array}{@{\,}c@{\,}}
\! z_0/6 \!\\
0
\end{array}
\right) + 
\left(
\begin{array}{@{\,}c@{\,}}
(1\! -\!\alpha )/2 \\
z_0/6 
\end{array}
\right) T + 
\left(
\begin{array}{@{\,}c@{\,}}
A_4 \\
B_3
\end{array}
\right) T^2 + \cdots , \\
& (ii) &\!\! \left(
\begin{array}{@{\,}c@{\,}}
x\\
y
\end{array}
\right)\! =\! -\left(
\begin{array}{@{\,}c@{\,}}
1 \\
0
\end{array}
\right) T^{-2} + 
\left(
\begin{array}{@{\,}c@{\,}}
0 \\
1
\end{array}
\right) T^{-1} -
\left(
\begin{array}{@{\,}c@{\,}}
\! z_0/6 \!\\
0
\end{array}
\right) - 
\left(
\begin{array}{@{\,}c@{\,}}
 \! (1 \! + \! \alpha )/2 \! \\
z_0/6 
\end{array}
\right) T - 
\left(
\begin{array}{@{\,}c@{\,}}
A_4 \\
B_3
\end{array}
\right) T^2 + \cdots,
\end{eqnarray*}
where $B_3 = (1 - \alpha )/4$ for the first line, $B_3 =  (1 + \alpha )/4$ for the second line 
and $A_4$ is an arbitrary constant.

For the fourth Painlev\'{e} equation ($\text{P}_\text{IV}$), there are three types of the Laurent series given by
\begin{eqnarray*}
& (i) &  \!\! \left(
\begin{array}{@{\,}c@{\,}}
x\\
y
\end{array}
\right)\! =\! \left(
\begin{array}{@{\,}c@{\,}}
1 \\
0
\end{array}
\right) T^{-1} + 
\left(
\begin{array}{@{\,}c@{\,}}
z_0 \\
0
\end{array}
\right) +
\left(
\begin{array}{@{\,}c@{\,}}
\! (2+z_0^2-2\theta _\infty + 4\kappa_0)/3 \!\\
2\kappa_0
\end{array}
\right) T + 
\left(
\begin{array}{@{\,}c@{\,}}
A_3 \\
B_3 
\end{array}
\right) T^2 + \cdots , \\
& (ii) &  \!\! \left(
\begin{array}{@{\,}c@{\,}}
x\\
y
\end{array}
\right)\! =\! \left(
\begin{array}{@{\,}c@{\,}}
\! -1 \! \\
\! -1 \!
\end{array}
\right) T^{-1}\! + \! 
\left(
\begin{array}{@{\,}c@{\,}}
z_0 \\
\! -z_0 \!
\end{array}
\right) \!+\! \frac{1}{3}
\left(
\begin{array}{@{\,}c@{\,}}
\! 6-z_0^2 + 2\theta _\infty - 4\kappa_0 \!\\
\! -6-z_0^2-4\theta _\infty + 2\kappa_0 \!
\end{array}
\right) T \!+ \!
\left(
\begin{array}{@{\,}c@{\,}}
A_3 \\
B_3 
\end{array}
\right) T^2 + \cdots , \\
& (iii) &  \!\! \left(
\begin{array}{@{\,}c@{\,}}
x\\
y
\end{array}
\right)\! =\! \left(
\begin{array}{@{\,}c@{\,}}
0 \\
1
\end{array}
\right) T^{-1} - 
\left(
\begin{array}{@{\,}c@{\,}}
0 \\
z_0
\end{array}
\right) -
\left(
\begin{array}{@{\,}c@{\,}}
2\theta _\infty\\
\! (2-z_0^2-4\theta _\infty +2\kappa_0)/3\!
\end{array}
\right) T + 
\left(
\begin{array}{@{\,}c@{\,}}
A_3 \\
B_3 
\end{array}
\right) T^2 + \cdots .
\end{eqnarray*}
$A_3$ is an arbitrary constant and $B_3$ is a certain constant depending on $A_3$.

Let us consider a general system (\ref{2-1b}) satisfying the assumptions given in Sec.2.1;
the Newton diagram consists of one compact face that lies on the plane $p_1 x + p_2 y + rz = s$,
and $g_i$ satisfies $g_i(\lambda ^{p_1}x , \lambda ^{p_2}y , \lambda ^{r}z) \sim o(\lambda ^{s-r+p_i})$.
In this case, the system has a formal series solution of the form
\begin{equation}
\left\{ \begin{array}{l}
\displaystyle x(z) = \sum^\infty_{n=0} A_n(z-z_0)^{-p_1 + n},  \\
\displaystyle y(z) = \sum^\infty_{n=0}B_n(z-z_0)^{-p_2 + n}.
\end{array} \right.
\end{equation}
The coefficients $A_n$ and $B_n$ are determined by substituting the series into the equation.
If the series solution represents a general solution, it includes an arbitrary parameter other than $z_0$.
The Kovalevskaya exponent $\kappa$ is defined to be the least integer $n$ such that the coefficient $(A_n, B_n)$
includes an arbitrary parameter.
For the Laurent series solution of ($\text{P}_\text{I}$), $\kappa = 6$.
For ($\text{P}_\text{II}$), $\kappa = 4$ for both series,
and for ($\text{P}_\text{IV}$), $\kappa = 3$ for all Laurent series solutions.
Note that the Kovalevskaya exponents of them coincide with the weighted degrees of Hamiltonian functions given in Sec.2.1.
In Sec.5, it is shown that 
the Kovalevskaya exponent coincides with an eigenvalue of a Jacobi matrix of a certain vector field,
and the exponent is invariant under the automorphism of $\C P^3(p,q,r,s)$.
See \cite{Abl, Bor, Chi4,  HuYa, Yos1, Yos2} for more general definition and properties of the Kovalevskaya exponent.

%%%%%%%%%%%%%%%%%%%%%%%%%%%%%%%%%%%%%%%%%%%%%%%%%%%%%%%%%%%%%%%%%%%%%%%%%%%%%%%%%%%%%%%%%%%%%%%%%%

\subsection{Relation with dynamical systems theory}

In this section, a relationship between the Painlev\'{e} equations
and the normal form theory of dynamical systems is shown.
The Painlev\'{e} equations will be obtained from certain singular perturbed problems of vector fields.

Let us consider a singular perturbation problem of the form
\begin{equation}
\left\{ \begin{array}{l}
\dot{\bm{x}} = \bm{f}(\bm{x}, \bm{z}, \varepsilon ),  \\
\dot{\bm{z}} = \varepsilon  \bm{g}(\bm{x}, \bm{z}, \varepsilon ),  \\
\end{array} \right.
\label{2-22}
\end{equation}
where $\bm{x}\in \R^{m}, \bm{z}\in \R^n$, and $(\bm{f}, \bm{g})$ is a smooth vector field on $\R^{m+n}$
parameterized by a small parameter $\varepsilon \in \R$.
The dot $(\,\dot{\,\,}\,)$ denotes the derivative with respect to time $t\in \R$.
Such a system is called a fast-slow system because it is characterized by two different time scales;
fast motion $\bm{x}$ and slow motion $\bm{z}$.
This structure yields nonlinear phenomena such as a relaxation oscillation,
which is observed in many physical, chemical and biological problems.
See Grasman~\cite{Gra}, Hoppensteadt and Izhikevich~\cite{Hop} and references therein for applications of fast-slow systems.
The unperturbed system is defined by putting $\varepsilon =0$ as
\begin{equation}
\dot{\bm{x}} = \bm{f}(\bm{x}, \bm{z}, 0 ), \quad \dot{\bm{z}} = 0.
\end{equation}
Since $\bm{z}$ is a constant for the unperturbed system,
it is regarded as a parameter of the fast system $\dot{\bm{x}} = \bm{f}(\bm{x}, \bm{z}, 0 )$.

It is known that when $\bm{f} \sim O(1)$ as $\varepsilon \to 0$ in some region of $\R^{m+n}$,
the dynamics of (\ref{2-22}) is approximately governed by the first system $\dot{\bm{x}} = \bm{f}(\bm{x}, \bm{z}, 0 )$,
and when $\bm{f} \sim 0$ while $D\bm{f} \sim O(1)$, 
the dynamics of (\ref{2-22}) is approximately governed by the slow system 
$\dot{\bm{z}} = \varepsilon \bm{g}(\varphi(\bm{z}), \bm{z}, 0 )$, where $D\bm{f}$ is the derivative of $\bm{f}$ with respect to $\bm{x}$
and $\varphi$ is a function satisfying $\bm{f}(\varphi(\bm{z}), \bm{z}, 0 ) = 0$.
However, if both of $\bm{f}$ and $D\bm{f}$ are nearly equal to zero, both of the fast and slow motion should be
taken into account and a nontrivial dynamics may occur.
The condition
\begin{eqnarray*}
\bm{f}(\bm{x}_0, \bm{z}_0, 0 )=0, \quad D\bm{f}(\bm{x}_0, \bm{z}_0, 0 )=0
\end{eqnarray*}
implies that the first system $\dot{\bm{x}} = \bm{f}(\bm{x}, \bm{z}, 0 )$ undergoes a bifurcation
at $\bm{x}=\bm{x}_0$ with a bifurcation parameter $\bm{z}=\bm{z}_0$.
A type of a bifurcation almost determines the local dynamics of (\ref{2-22}) around $(\bm{x}, \bm{z})=(\bm{x}_0, \bm{z}_0 )$

For the most generic case, in which the fast system undergoes a saddle-node bifurcation,
it is well known that a local behavior of (\ref{2-22}) is governed by the Airy equation $d^2u/dz^2 = zu$.
In particular, the asymptotic analysis of the Airy function plays an important role, see \cite{Kru1, Gil}.
Chiba \cite{Chi2} found that when the fast system  $\dot{\bm{x}} = \bm{f}(\bm{x}, \bm{z}, 0 )$ undergoes 
a Bogdanov-Takens bifurcation, then a local behavior of (\ref{2-22}) is determined 
by the asymptotic analysis of Boutroux's tritronqu\'{e}e solution of the first Painlev\'{e} equation.
This result was applied to prove the existence of a periodic orbit and chaos in a certain biological system \cite{Chi3}.
Here, we will demonstrate how the first, second and fourth Painlev\'{e} equations are obtained from fast-slow systems.
For a normal form and versal deformation of germs of vector fields, the readers can refer to \cite{Cho}.

Suppose that a one-dimensional dynamical system $\dot{x} = f(x)$ lies on a codimension 1 bifurcation at $x=0$.
This means that $f$ satisfies
\begin{equation}
f(0) = f'(0) = 0, \quad f''(0) \neq 0.
\end{equation}
The normal form is given by $f(x) = x^2$, and its versal deformation is $\dot{x} = x^2 + z$ or $\dot{x} = x^2 + zx$
with a deformation parameter $z\in \R$.
The former is called a saddle-node bifurcation and the latter is a transcritical bifurcation.
The fast-slow system having the saddle-node as an unperturbed fast system is written by
\begin{equation}
\dot{x} = x^2 + z + \varepsilon f(x,z,\varepsilon ), \quad \dot{z} = \varepsilon g(x,z,\varepsilon ).
\label{2-25}
\end{equation}
We also assume the generic condition $g(0,0,0) \neq 0$ so that we can write $g(x,z,\varepsilon ) = 1 + O(x,z,\varepsilon )$
without loss of generality.
In order to investigate the local behavior of the system near $(x,z) = (0,0)$ for a small $\varepsilon $,
we rewrite Eq.(\ref{2-25}) as a three dimensional system
\begin{equation}
\left\{ \begin{array}{l}
\dot{x} = x^2 + z + \varepsilon f(x,z,\varepsilon )  \\
\dot{z} = \varepsilon  + \varepsilon \cdot O(x,z,\varepsilon )  \\
\dot{\varepsilon } = 0,
\end{array} \right.
\label{2-26}
\end{equation}
by adding the trivial equation $\dot{\varepsilon } = 0$.
For this system, we perform the weighted blow-up at the origin defined by
\begin{equation}
\left(
\begin{array}{@{\,}c@{\,}}
x \\
z \\
\varepsilon 
\end{array}
\right) = \left(
\begin{array}{@{\,}l@{\,}}
r_1 \\
r_1^2 Z_1 \\
r_1^3 \varepsilon _1
\end{array}
\right) = \left(
\begin{array}{@{\,}l@{\,}}
r_2 X_2 \\
r_2^2 \\
r_2^3 \varepsilon _2
\end{array}
\right) = \left(
\begin{array}{@{\,}l@{\,}}
r_3 X_3 \\
r_3^2 Z_3 \\
r_3^3
\end{array}
\right).
\end{equation}
The weight (exponents of $r_i$'s) $(1,2,3)$ can be found by the Newton diagram of Eq.(\ref{2-26}) as in Sec.2.1.
The exceptional divisor of the blow-up is $\C P^2 (1,2,3)$ given by the set $\{ r_1 = 0\} \cup \{ r_2 = 0\} \cup \{ r_3 = 0\}$.
On the $(X_3, Z_3, \varepsilon _3)$ chart, Eq.(\ref{2-26}) is written as
\begin{equation}
\dot{X}_3 = r_3X_3^2 + r_3Z_3 + O(r_3^2), \quad \dot{Z}_3 = r_3 + O(r_3^2), \quad \dot{r}_3 = 0.
\end{equation}
This is equivalent to
\begin{eqnarray*}
\frac{dX_3}{dZ_3} = \frac{X_3^2 + Z_3 + O(r_3)}{1 + O(r_3)}.
\end{eqnarray*}
In particular, on the exceptional divisor $\C P^2 (1,2,3)$, it is reduced to the Riccati equation
\begin{equation}
\frac{dX_3}{dZ_3} = X_3^2 + Z_3,
\end{equation}
which is equivalent to the Airy equation $u'' = -Z_3u$ by $X_3 = -u'/u$.
Similarly, if we use the transcritical bifurcation as the fast system
and apply the weighted blow-up with the weight $(1,1,2)$, we obtain the Hermite equation
$u'' - Z_3u' - \alpha u = 0$ on the exceptional divisor.
See \cite{Chi} for the analysis of the Airy equation based on the geometry of $\C P^2 (1,2,3)$.

The Painlev\'{e} equations are obtained from codimension 2 bifurcations in a similar manner.
Suppose that a 2-dim system of $(x,y)$ undergoes a generic codimension 2 bifurcation called the Bogdanov-Takens bifurcation.
The normal form is given by $\dot{x} = y^2 + xy,\, \dot{y} = x$; that is, the linear part has two zero eigenvalues.
Its versal deformation is
\begin{equation}
\left\{ \begin{array}{l}
\dot{x} = y^2 + xy + z,  \\
\dot{y} = x,  \\
\end{array} \right.
\end{equation}
where $z$ is a deformation parameter.
The fast-slow system having it as an unperturbed fast system is written by
\begin{equation}
\left\{ \begin{array}{l}
\dot{x} = y^2 + xy + z + \varepsilon f_1(x,y,z,\varepsilon )  \\
\dot{y} = x+ \varepsilon f_2(x,y,z,\varepsilon ), \\
\dot{z} = \varepsilon  + \varepsilon \cdot O(x,y,z,\varepsilon ), \\
\dot{\varepsilon } = 0,
\end{array} \right.
\label{2-31}
\end{equation}
where the trivial equation $\dot{\varepsilon } = 0$ is added as before.
For this system, we introduce the weighted blow-up with the weight $(3,2,4,5)$ defined by
\begin{equation}
\left(
\begin{array}{@{\,}c@{\,}}
x \\
y \\
z \\
\varepsilon 
\end{array}
\right) = \left(
\begin{array}{@{\,}l@{\,}}
r_1^3 \\
r_1^2 Y_1 \\
r_1^4 Z_1\\
r_1^5 \varepsilon _1
\end{array}
\right) = \left(
\begin{array}{@{\,}l@{\,}}
r_2^3 X_2\\
r_2^2  \\
r_2^4 Z_2\\
r_2^5 \varepsilon _2
\end{array}
\right) = \left(
\begin{array}{@{\,}l@{\,}}
r_3^3 X_3  \\
r_3^2 Y_3 \\
r_3^4  \\
r_3^5 \varepsilon _3
\end{array}
\right) = \left(
\begin{array}{@{\,}l@{\,}}
r_4^3 X_4 \\
r_4^2 Y_4 \\
r_4^4 Z_4\\
r_4^5 
\end{array}
\right).
\end{equation}
The weight $(3,2,4,5)$ can be obtained via the Newton diagram of the system (\ref{2-31}).
The exceptional divisor of the blow-up is the weighted projective space $\C P^3(3,2,4,5)$
given as the set $\{ r_1 =0 \} \cup \{ r_2 =0 \} \cup\{ r_3 =0 \} \cup\{ r_4 =0 \}$.
Transforming the system (\ref{2-31}) to the $(X_4, Y_4, Z_4, r_4)$ chart, we obtain
\begin{eqnarray*}
\left\{ \begin{array}{l}
\dot{X}_4 = r_4Y_4^2+r_4Z_4+ O(r_4^2), \\
\dot{Y}_4 = r_4X_4 + O(r_4^2),  \\
\dot{Z}_4 = r_4 + O(r_4^2).
\end{array} \right.
\end{eqnarray*}
As $r_4 \to 0$, it is reduced to the 
first Painlev\'{e} equation $dX_4/dZ_4 = Y_4^2 + Z_4,\, dY_4/dZ_4 = X_4$.
This implies that the dynamics on the divisor $\C P^3(3,2,4,5)$ is governed by the compactified 
first Painlev\'{e} equation, and a local behavior of the system (\ref{2-31}) near $(x,y,z) = (0,0,0)$
can be investigated by a global analysis of the first Painlev\'{e} equation.

Next, we consider a 2-dim system that undergoes a Bogdanov-Takens bifurcation with $\Z_2$-symmetry $(x,y) \mapsto (-x,-y)$.
The versal deformation of the normal form is given by
\begin{equation}
\left\{ \begin{array}{l}
\dot{x} = y^3 - xy^3 + zy,  \\
\dot{y} = x,  \\
\end{array} \right.
\end{equation}
with a deformation parameter $z$.
The fast-slow system having it as an unperturbed fast system is written by
\begin{equation}
\left\{ \begin{array}{l}
\dot{x} = y^3 - xy^3 + zy + \varepsilon f_1(x,y,z,\varepsilon )  \\
\dot{y} = x+ \varepsilon f_2(x,y,z,\varepsilon ), \\
\dot{z} = \varepsilon  + \varepsilon \cdot O(x,y,z,\varepsilon ), \\
\dot{\varepsilon } = 0.
\end{array} \right.
\label{2-34}
\end{equation}
For this system, we introduce the weighted blow-up with the weight $(2,1,2,3)$,
which is found by the Newton diagram of the system (\ref{2-34}).
On the $(X_4, Y_4, Z_4, r_4)$ chart, it provides
\begin{eqnarray*}
\left\{ \begin{array}{l}
\dot{X}_4 = r_4Y_4^3 + r_4Z_4Y_4 - r_4^3 X_4Y_4^3 + r_4 f_1(r_4^2X_4, r_4Y_4, r_4^2Z_4, r_4^3), \\
\dot{Y}_4 = r_4X_4 + O(r_4^2),  \\
\dot{Z}_4 = r_4 + O(r_4^2).
\end{array} \right.
\end{eqnarray*}
Note that $r_4 f_1(r_4^2X_4, r_4Y_4, r_4^2Z_4, r_4^3) = \alpha r_4 + O(r_4^2)$,
where $\alpha := f_1(0,0,0,0)$.
As $r_4 \to 0$, this system is reduced to the second Painlev\'{e} equation $X_4' = Y_4^3 + Y_4Z_4 + \alpha ,\, Y_4' = X_4$
with a parameter $\alpha $.

Finally, we consider a 2-dim system that undergoes a Bogdanov-Takens bifurcation with $\Z_3$-symmetry.
Using the complex variable $\eta \in \C$, the normal form of such a bifurcation is given by
$\dot{\eta} = \eta |\eta |^2 + \overline{\eta}^2$.
Note that this system is invariant under the $\Z_3$ action $\eta \mapsto e^{2\pi i/ 3} \eta$.
The versal deformation of the normal form is 
\begin{equation}
\dot{\eta} = \eta |\eta |^2 + \overline{\eta}^2 +\eta z.
\end{equation}
where $z\in \C$ is a deformation parameter.
Putting $\eta = x + iy, z = z_1 + iz_2$ yields
\begin{eqnarray*}
\left\{ \begin{array}{l}
\dot{x} = x^2 - y^2 + z_1x - z_2y + O(\eta^3),  \\
\dot{y} = -2xy + z_1y + z_2 x + O(\eta^3).  \\
\end{array} \right.
\end{eqnarray*}
We assume $z_1 = 0$ so that the above system may become a Hamiltonian system,
and change the notation $z_2 \mapsto z$ to obtain
\begin{eqnarray*}
\left\{ \begin{array}{l}
\dot{x} = x^2 - y^2 - zy + O(\eta^3),  \\
\dot{y} = -2xy  + z x + O(\eta^3).  \\
\end{array} \right.
\end{eqnarray*}
The fast-slow system having it as an unperturbed fast system is written by
\begin{equation}
\left\{ \begin{array}{l}
\dot{x} = x^2 - y^2 - zy + O(\eta^3) + \varepsilon f_1(x,y,z,\varepsilon )  \\
\dot{y} = -2xy  + z x + O(\eta^3) + \varepsilon f_2(x,y,z,\varepsilon ), \\
\dot{z} = \varepsilon  + \varepsilon \cdot O(x,y,z,\varepsilon ), \\
\dot{\varepsilon } = 0.
\end{array} \right.
\label{2-36}
\end{equation}
For this system, we introduce the weighted blow-up with the weight $(1,1,1,2)$.
Moving to the $(X_4,Y_4,Z_4)$ chart and putting $r_4 = 0$ as before, it turns out that the system (\ref{2-36})
is reduced the system
\begin{equation}
\left\{ \begin{array}{l}
\displaystyle \frac{dX_4}{dZ_4} = X_4^2 - Y_4^2 - Z_4Y_4 + \alpha _1  \\
\displaystyle \frac{dY_4}{dZ_4} = -2X_4Y_4  + Z_4Y_4 + \alpha _2, \\
\end{array} \right.
\label{2-37}
\end{equation}
where $\alpha_i := f_i(0,0,0,0)$ is a constant ($i=1,2$).
This is equivalent to the fourth Painlev\'{e} equation (\ref{P4}) through an affine transformation of $(X_4, Y_4)$.

\begin{table}[h]
\begin{center}
\begin{tabular}{|c|c|c|}
\hline
Bifurcation type & Exceptional divisor & Equation on the divisor \\ \hline \hline
saddle-node & $\C P^2(1,2,3)$ & Airy \\ \hline
transcritical & $\C P^2(1,1,2)$ & Hermite \\ \hline
Bogdanov-Takens (BT) & $\C P^3(3,2,4,5)$ & ($\text{P}_\text{I}$) \\ \hline
BT with $\Z_2$ symmetry & $\C P^3(2,1,2,3)$ &($\text{P}_\text{II}$) \\ \hline
BT with $\Z_3$ symmetry & $\C P^3(1,1,1,2)$ & ($\text{P}_\text{IV}$) \\ \hline
\end{tabular}
\end{center}
\caption{Differential equations obtained by the weighted blow-up of the fast-slow systems.}
\end{table}

The results are summarized in Table 1.
It is remarkable that the weights $(p,q,r,s)$ derived in Sec.2.1 via the Newton diagrams are determined only by the 
versal deformations of the codimension 2 bifurcations.
Further, the Painlev\'{e} equations are obtained in a compactified manner on $\C P^3(p,q,r,s)$,
which is the exceptional divisor of the weighted blow-up.

%%%%%%%%%%%%%%%%%%%%%%%%%%%%%%%%%%%%%%%%%%%%%%%%%%%%%%%%%%%%%%%%%%%%%%%%%%%%%%%%%%%%%%%%%%%%%%%%%%
%%%%%%%%%%%%%%%%%%%%%%%%%%%%%%%%%%%%%%%%%%%%%%%%%%%%%%%%%%%%%%%%%%%%%%%%%%%%%%%%%%%%%%%%%%%%%%%%%%

\section{Singular normal forms and the Painlev\'{e} property}

Recall that a singularity $z = z_*$ of a solution of a differential equation is called 
movable if the position of $z_*$ depends on the choice of an initial condition.
In this section, we give local analysis near such movable singularities based on the normal form theory.
Our purpose is to show that near movable singularities, the Painlev\'{e} equations are locally transformed into 
integrable systems called the singular normal form.
Further, we will give a new proof of the Painlev\'{e} property; any solutions of ($\text{P}_\text{I}$),
($\text{P}_\text{II}$) and ($\text{P}_\text{IV}$) are meromorphic on $\C$ (for this purpose,
we do not use the Laurent series given in Sec.2.3).  

%%%%%%%%%%%%%%%%%%%%%%%%%%%%%%%%%%%%%%%%%%%%%%%%%%%%%%%%%%%%%%%%%%%%%%%%%%%%%%%%%%%%%%%%%%%%%%%%%%

\subsection{The first Painlev\'{e} equation}

($\text{P}_\text{I}$) is given on the weighted projective space $\C P^3(3,2,4,5)$
as a tuple of equations (\ref{P1}), (\ref{P11}), (\ref{P12}) and (\ref{P13}). 
Coordinate transformations between inhomogeneous coordinates are given by
\begin{eqnarray}
\left\{ \begin{array}{rrrr}
x = & \varepsilon _1^{-3/5} = & X_2\varepsilon _2^{-3/5} =& X_3\varepsilon _3^{-3/5} \\
y = & Y_1\varepsilon _1^{-2/5} = & \varepsilon _2^{-2/5}=& Y_3\varepsilon _3^{-2/5} \\
z = & Z_1\varepsilon_1 ^{-4/5} = & Z_2\varepsilon _2^{-4/5}=& \varepsilon _3^{-4/5} .
\end{array} \right.
\label{3-1}
\end{eqnarray} 
Due to the orbifold structure of $\C P^3(3,2,4,5)$, local charts $(Y_1, Z_1, \varepsilon _1)$,
$(X_2, Z_2,\varepsilon _2)$ and $(X_3, Y_3, \varepsilon _3)$ should be divided by the $\Z_3, \Z_2$ and $\Z_4$
actions, respectively, defined by
\begin{eqnarray}
& & (Y_1, Z_1, \varepsilon _1) \mapsto (\omega Y_1, \omega ^2 Z_1, \omega \varepsilon _1), 
\quad \omega := e^{2\pi i/3},  \label{3-2} \\
& & (X_2, Z_2,\varepsilon _2) \mapsto (-X_2, Z_2, -\varepsilon _2), \label{3-3} \\
& & (X_3, Y_3, \varepsilon _3) \mapsto (i X_3, -Y_3, -i \varepsilon _3). \label{3-4}
\end{eqnarray} 
Indeed, Eqs.(\ref{P11}),(\ref{P12}),(\ref{P13}) are invariant under these actions.
For our purposes, it is convenient to rewrite Eqs.(\ref{P11}), (\ref{P12}) and (\ref{P13}) 
as 3-dim vector fields (autonomous ODEs) given by
\begin{equation}
\left\{ \begin{array}{l}
\displaystyle \dot{Y}_1 = 3 - 12Y_1^3 - 2Y_1Z_1,  \\
\displaystyle \dot{Z}_1 = 3\varepsilon _1 - 24Y_1^2 Z_1 - 4 Z_1^2,  \\
\displaystyle \dot{\varepsilon }_1 = \varepsilon _1 (-30 Y_1^2 - 5Z_1),
\end{array} \right.
\label{3-5}
\end{equation}
\begin{equation}
\left\{ \begin{array}{l}
\displaystyle \dot{X}_2 = (-12 - 2Z_2 + 3X_2^2)/X_2,  \\
\displaystyle \dot{Z}_2 = (-2\varepsilon _2  + 4 X_2Z_2)/X_2,  \\
\displaystyle \dot{\varepsilon }_2 = 5\varepsilon _2,
\end{array} \right.
\label{3-6}
\end{equation}
\begin{equation}
\left\{ \begin{array}{l}
\displaystyle \dot{X}_3 =24Y_3^2 + 4 - 3X_3\varepsilon _3,  \\
\displaystyle \dot{Y}_3 = 4X_3 - 2Y_3\varepsilon _3,  \\
\displaystyle \dot{\varepsilon }_3 = -5\varepsilon _3^2,
\end{array} \right.
\label{3-7}
\end{equation}
where $(\dot{\,\,}) = d/dt$ and $t$ is an additional parameter.

Recall the decomposition (\ref{2-7}) with (\ref{2-8}).
According to Eqs.(\ref{3-5}),(\ref{3-6}),(\ref{3-7}), the set $\C P^2 (3,2,4)$ is an invariant manifold
of the vector fields; that is, $\varepsilon _i(t) \equiv 0$ when $\varepsilon _i (0) = 0$ at an initial time.
The dynamics on the invariant manifold describes the behavior of ($\text{P}_\text{I}$) near infinity.
In particular, fixed points of the vector fields play an important role.
Vector fields (\ref{3-5}),(\ref{3-6}),(\ref{3-7}) have exactly two fixed points on $\C P^2 (3,2,4)$;
\\[0.2cm]
\textbf{(i).} $(X_2, Z_2, \varepsilon _2) = (\pm 2,0,0)$.

Due to the $\Z_2$ action on the $(X_2,Z_2, \varepsilon _2)$-coordinates,
two points $(2,0,0)$ and $(-2,0,0)$ represent the same point on $\C P^3(3,2,4,5)$
and it is sufficient to consider one of them.
We will show that this fixed point corresponds to movable singularities of ($\text{P}_\text{I}$).
By applying the normal form theory of vector fields to this point, 
we will construct the singular normal form of ($\text{P}_\text{I}$).
In Sec.4, the space of initial conditions for ($\text{P}_\text{I}$) is constructed by applying
the weighted blow-up at this point.
\\[0.2cm]
\textbf{(ii).} $(X_3, Y_3, \varepsilon _3) = (0, (-6)^{-1/2} ,0)$

Again, two points should be identified due to the $\Z_4$ action on $(X_3, Y_3, \varepsilon _3)$.
This fixed point corresponds to the irregular singular point of ($\text{P}_\text{I}$) because
$\varepsilon _3 = 0$ provides $z = \infty$.

Note that fixed points obtained from the $(Y_1, Z_1, \varepsilon _1)$-coordinates are the same as one of the above.
For example, the fixed point $(Y_1, Z_1, \varepsilon _1) = ((1/4)^{1/3},0,0)$ is the same as 
$(X_2, Z_2, \varepsilon _2) = (\pm 2,0,0)$ due to (\ref{3-1}).

At first, let us show that ($\text{P}_\text{I}$) is locally transformed into a linear system around $(X_2,Z_2,\varepsilon _2) = (2,0,0)$.
By putting $\hat{X}_2 = X_2 - 2$, Eq.(\ref{3-6}) is rewritten as
\begin{equation}
\left\{ \begin{array}{l}
\displaystyle \dot{\hat{X}}_2 = 6\hat{X}_2 - Z_2 + \frac{-3\hat{X}_2^2 + \hat{X}_2Z_2}{2+\hat{X}_2},  \\
\displaystyle \dot{Z}_2 = 4Z_2 - \varepsilon _2 + \frac{\hat{X}_2\varepsilon _2}{2+\hat{X}_2},  \\
\displaystyle \dot{\varepsilon }_2 = 5\varepsilon _2.
\end{array} \right.
\label{3-8}
\end{equation}
The origin is a fixed point with the Jacobi matrix
\begin{eqnarray}
J = \left(
\begin{array}{@{\,}ccc@{\,}}
6 & - 1 & 0 \\
0 & 4 & - 1 \\
0 & 0 & 5
\end{array}
\right).
\label{3-9}
\end{eqnarray}
The eigenvalues $\lambda =6,4,5$ satisfy the conditions on the Poincar\'{e} linearization theorem 
(for the convenience of the reader, we give the statement of this theorem in the end of this subsection).
Hence, there exists a neighborhood $U$ of the origin and a local analytic transformation defined on $U$ of the form
\begin{eqnarray*}
\left(
\begin{array}{@{\,}c@{\,}}
\hat{X}_2 \\
Z_2
\end{array}
\right) \mapsto
\left(
\begin{array}{@{\,}c@{\,}}
\hat{u}_2 \\
v_2
\end{array}
\right) = \left(
\begin{array}{@{\,}c@{\,}}
\hat{X}_2 + \varphi _1(\hat{X}_2,Z_2,\varepsilon _2) \\
Z_2 + \hat{\varphi} _2 (\hat{X}_2, Z_2, \varepsilon _2)
\end{array}
\right),
\end{eqnarray*}
such that Eq.(\ref{3-8}) is linearized as
\begin{equation}
\left\{ \begin{array}{l}
\displaystyle \dot{\hat{u}}_2 = 6\hat{u}_2 - v_2,  \\
\displaystyle \dot{v}_2 = 4v_2 - \varepsilon _2,  \\
\displaystyle \dot{\varepsilon }_2 = 5\varepsilon _2,
\end{array} \right.
\label{3-10}
\end{equation}
where local analytic functions $\varphi_1$ and $\hat{\varphi}_2$ satisfy
$\varphi_1, \hat{\varphi}_2 \sim O(|| \bm{X} ||^2)$, $\bm{X}:= (\hat{X}_2, Z_2, \varepsilon _2)$.
Note that we need not change $\varepsilon _2$ because the equation of $\varepsilon _2$ is already linear.
Furthermore, we have $\hat{\varphi}_2 (\hat{X}_2, Z_2, 0) = 0$ because the equation of $Z_2$ is linear 
when $\varepsilon _2  =0$.
Thus, we can set $\hat{\varphi}_2 = \varepsilon _2 \varphi_2$ and the above transformation takes the form
\begin{eqnarray}
\left(
\begin{array}{@{\,}c@{\,}}
\hat{X}_2 \\
Z_2
\end{array}
\right) \mapsto
\left(
\begin{array}{@{\,}c@{\,}}
\hat{u}_2 \\
v_2
\end{array}
\right) = \left(
\begin{array}{@{\,}c@{\,}}
\hat{X}_2 + \varphi _1(\hat{X}_2,Z_2,\varepsilon _2) \\
Z_2 + \varepsilon _2 \varphi _2 (\hat{X}_2, Z_2, \varepsilon _2)
\end{array}
\right),
\label{3-11}
\end{eqnarray}
with $\varphi _1 \sim O(|| \bm{X} ||^2)$ and $\varphi _2 \sim O(X_2, Z_2, \varepsilon _2)$.

For the linear system Eq.(\ref{3-10}), let us change coordinates as
$\hat{u}_2 = u_2 - 2$, and $\tilde{x} = u_2\varepsilon _2^{-3/5}, \tilde{y} = \varepsilon _2^{-2/5}, \tilde{z} = v_2 \varepsilon _2^{-4/5}$;
that is, we move to the original chart for ($\text{P}_\text{I}$).
We can verify that $(\tilde{x}, \tilde{y}, \tilde{z})$ satisfies the equation $d^2 \tilde{y}/d\tilde{z}^2 = 6\tilde{y}^2$,
whose solution can be expressed by the Weierstrass's elliptic function.
The relation between $(x,y,z)$ and $(\tilde{x}, \tilde{y}, \tilde{z})$ is
\begin{eqnarray*}
\left(
\begin{array}{@{\,}c@{\,}}
\tilde{x} \\
\tilde{y} \\
\tilde{z} 
\end{array}
\right) = \left(
\begin{array}{@{\,}c@{\,}}
x + y^{3/2} \varphi_1 (xy^{-3/2}-2, zy^{-2}, y^{-5/2}) \\
y \\
z + y^{-1/2}\varphi_2 (xy^{-3/2}-2, zy^{-2}, y^{-5/2})
\end{array}
\right).
\end{eqnarray*}
In particular, $\tilde{y} = y$ is not changed.
Now we have obtained
\\[0.2cm]
\textbf{Proposition \thedef.}
There is a local analytic transformation $(x,z) \mapsto (\tilde{x},\tilde{z})$ defined near $(\hat{X}_2, Z_2, \varepsilon _2) = (0,0,0)$
such that ($\text{P}_\text{I}$) $y'' = 6y^2 + z$ is transformed into the integrable system $y'' = 6y^2$.

This fact was first obtained by \cite{CosCos} for ($\text{P}_\text{I}$).
Our proof using the weighted projective space and the Poincar\'{e} theorem is also applicable to
the second Painlev\'{e} to sixth Painlev\'{e} equations to prove that they are locally transformed to
solvable systems.
Since Eq.(\ref{3-10}) is linear, we can construct two integrals explicitly as
\begin{eqnarray*}
C_1 = \varepsilon _2^{-4/5}v_2 + \varepsilon _2^{1/5}, \quad 
C_2 = \frac{1}{2}\varepsilon _2^{-1/5} + \varepsilon _2^{-6/5}\hat{u}_2 - \frac{1}{2}\varepsilon _2^{-6/5}v_2.
\end{eqnarray*}
By applying the transformations (\ref{3-11}), $\hat{X}_2 = X_2 - 2$ and (\ref{3-1}), 
we obtain the local integrals of ($\text{P}_\text{I}$) of the form
\begin{eqnarray}
\left\{ \begin{array}{ll}
\displaystyle C_1 = z+ y^{-1/2}+y^{-1/2} \varphi _2(xy^{-3/2}-2, zy^{-2}, y^{-5/2}),  \\[0.2cm]
\displaystyle C_2 = \frac{1}{2}y^{1/2}-2y^3 + xy^{3/2} - \frac{1}{2}yz + y^3 \varphi _1(\cdots ) - \frac{1}{2}y^{1/2} \varphi _2(\cdots ), \\
\end{array} \right.
\label{3-12}
\end{eqnarray}
Arguments of $\varphi _1$ and $\varphi _2$ in the second line are the same as that of the first line.
Now we give a new proof of the well known theorem:
\\[0.2cm]
\textbf{Theorem.\thedef.} Any solutions of ($\text{P}_\text{I}$) are meromorphic on $\C$.
\\

A well known proof of this result is essentially based on Painlev\'{e}'s argument modified 
by Hukuhara (\cite{OkaTak}, \cite{Hin}, see also \cite{Gro}).
Here, we will prove the theorem by applying the implicit function theorem to the above integrals.
\\[0.2cm]
\textbf{Proof.} 
Fix a solution $(x(z), y(z))$ of ($\text{P}_\text{I}$) with an initial condition
$(x(z_0), y(z_0)) = (x_0, y_0)$.
The existence theorem of solutions shows that the solution is holomorphic near $z_0$.
Let $B(z_0, R)$ be the largest disk of radius $R$ centered at $z_0$ such that 
the solution is holomorphic inside the disk.
Let $z_* \, (\neq \infty)$ be a singularity on the boundary of the disk
(if $R=\infty$, there remains nothing to prove).
The next lemma implies that the fixed point $(X_2, Z_2, \varepsilon _2) = (\pm 2,0,0)$
corresponds to the singularity $z_*$.
\\[0.2cm]
\textbf{Lemma.\thedef.} $(X_2, Z_2, \varepsilon _2) \to (\pm 2,0,0)$ as $z \to z_*$ along a curve $\Gamma $
inside the disk $B(z_0, R)$.
\\[0.2cm]
\textbf{Proof.} Suppose that there exists a sequence $\{z_n\}_{n=1}^\infty$ converging to $z_*$ 
on the curve $\Gamma $ such that both of $x(z_n)$ and $y(z_n)$ are bounded as $n\to \infty$.
Taking a subsequence if necessary, we can assume that $(x, y)$ converges to some point $(x_*, y_*)$.
Because of the existence theorem of solutions, a solution of ($\text{P}_\text{I}$) satisfying the initial condition
$(x_*, y_*, z_*)$ is holomorphic around this point, which contradicts with the definition of $z_*$.
Hence, either $x$ or $y$ diverges as $z \to z_*$.

(i) Suppose that $y \to \infty$ as $z \to z_*$.
We move to the $(X_2, Z_2, \varepsilon _2)$-coordinates.
Eq.(\ref{3-1}) provides
\begin{equation}
X_2 = xy^{-3/2}, \quad Z_2 = zy^{-2}, \quad \varepsilon _2 = y^{-5/2}.
\end{equation}
This immediately yields $Z_2 \to 0,\, \varepsilon _2 \to 0$ as $z \to z_*$.
Let us show $X_2 \to \pm 2$.
($\text{P}_\text{I}$) is a Hamiltonian system with the Hamiltonian function
$H = x^2/2 - 2y^3 - zy$.
Thus, the equality $H = -\int\! y(z)dz $ holds along a solution.
In the $(X_2, Z_2, \varepsilon _2)$-coordinates, this is written as
\begin{eqnarray*}
X_2(Z_2)^2 = 4 + \frac{2}{3}Z_2 - \frac{2}{3}\varepsilon _2(Z_2)^{6/5} \int^{Z_2}_{\xi}\! \varepsilon _2(z)^{-6/5}dz,
\end{eqnarray*}
where $(X_2(Z_2), \varepsilon _2(Z_2))$ is a solution of the ODE solved as a function of $Z_2$,
and $\xi$ is a certain nonzero number determined by the initial condition.
Since $Z_2 \to 0,\, \varepsilon _2 \to 0$ as $z \to z_*$, we obtain $X^2_2 \to 4$.

(ii) Suppose that $x \to \infty$ as $z \to z_*$.
In this case, we use the $(Y_1, Z_1, \varepsilon _1)$-coordinates given by
\begin{equation}
Y_1 = yx^{-2/3}, \quad Z_1 = zx^{-4/3}, \quad \varepsilon _1 = x^{-5/3}.
\end{equation}
By the assumption, we have $Z_1 \to 0,\, \varepsilon _1 \to 0$ as $z \to z_*$.
Then, we can show that $Y_1^3 \to 1/4$ as $z \to z_*$ by the same way as above.
This means that $(Y_1, Z_1, \varepsilon _1)$ converges to the fixed point $((1/4)^{1/3},0,0)$
of the vector field (\ref{3-5}).
It is easy to verify that this fixed point is the same point as $(X_2, Z_2, \varepsilon _2) = (\pm 2,0,0)$
if written in the $(X_2, Z_2, \varepsilon _2)$-coordinates. $\Box$

The sign of $X_2 = xy^{-3/2}$ depends on the choice of the branch of $y^{1/2}$
and two points $(2,0,0)$ and $(-2,0,0)$ are the essentially the same.
In what follows, we assume that $(X_2, Z_2, \varepsilon _2) \to (2,0,0)$ as $z \to z_*$.
Due to the above lemma, when $z$ is sufficiently close to $z_*$,
the solution is included in the neighborhood $U$, on which local holomorphic functions 
$\varphi _1$ and $\varphi _2$ are well defined.
Then, the integrals (\ref{3-12}) are available.
To apply the implicit function theorem, put
\begin{eqnarray}
\left(
\begin{array}{@{\,}c@{\,}}
w \\
u
\end{array}
\right) = \left(
\begin{array}{@{\,}c@{\,}}
y^{-1/2} \\
\frac{1}{2}y^{1/2}-2y^3 + xy^{3/2} - \frac{1}{2}yz
\end{array}
\right),
\end{eqnarray}
Then, (\ref{3-12}) takes the form
\begin{eqnarray*}
\left\{ \begin{array}{l}
\displaystyle  C_1 = z + w + w \varphi _2(uw^6+\frac{1}{2}zw^4-\frac{1}{2}w^5, z w^4, w^5)\\
\displaystyle  C_2 = u + w^{-6}\varphi _1 (\cdots ) - \frac{1}{2}w^{-1}\varphi _2 (\cdots ). \\
\end{array} \right.
\end{eqnarray*}
Note that $C_1 = z+ w + O(w^5) $ and $C_2 = u + O(w^2)$ as $w \to 0$.
Since $w\to 0$ as $z \to z_*$, the constant $C_1 = z_*$ is just the position of the singularity.
If we set
\begin{eqnarray*}
& & f_1(w,u, z) = z + w + w\varphi _2 (\cdots ) - z_* \\
& & f_2(w,u, z) = u + w^{-6}\varphi _1 (\cdots )- \frac{1}{2}w^{-1}\varphi _2 (\cdots ) -C_2,
\end{eqnarray*}
then $f_i(0,C_2,z_*) = 0$.
The Jacobi matrix of $(f_1, f_2)$ with respect to $(w,u)$ at $(w,u,z) = (0, C_2, z_*)$ is the identity matrix.
Hence, the implicit function theorem proves that there exists a local holomorphic
function $g(z)$ such that $f_1 = f_2 = 0$ is solved as $w = g(z) \sim O(z - z_*)$.
Since $y = w^{-2}$, $z = z_*$ is a pole of second order of $y$.
This completes the proof of Thm.3.2. $\Box$
\\

The Poincar\'{e} linearization theorem used in this subsection is stated as follows:
Let $Ax + f(x)$ be a holomorphic vector field on $\C^n$ with a fixed point $x=0$,
where $A$ is an $n\times n$ constant matrix and $f(x) \sim O(|x|^2)$ is a nonlinearity.
Let $\lambda _1 , \cdots , \lambda _n$ be eigenvalues of $A$.
We consider the following two conditions:
\\[0.2cm]
\textbf{(Nonresonance)} There are no $j\in \{ 1, \cdots ,n\}$ and non-negative integers $m_1, \cdots ,m_n$
satisfying the resonant condition
\begin{equation}
m_1 \lambda _1 + \cdots  + m_n \lambda _n = \lambda _j, \quad (m_1 + \cdots + m_n \geq 2).
\end{equation}
\textbf{(Poincar\'{e} domain)} The convex hull of $\{ \lambda _1, \cdots , \lambda _n\}$ in $\C$ does
not include the origin.
\\[0.2cm]
Suppose that $A$ is diagonal and eigenvalues satisfy the above two conditions.
Then, there exists a local analytic transformation $y = x + \varphi (x),\,\, \varphi (x) \sim O(|x|^2)$ defined 
near the origin such that the equation $dx/dt = Ax + f(x)$ is transformed into the linear system $dy/dt = Ay$.
See \cite{Cho} for the proof.

%%%%%%%%%%%%%%%%%%%%%%%%%%%%%%%%%%%%%%%%%%%%%%%%%%%%%%%%%%%%%%%%%%%%%%%%%%%%%%%%%%%%%%%%%%%%%%%%%%

\subsection{The second Painlev\'{e} equation}

($\text{P}_\text{II}$) is given on the weighted projective space $\C P^3(2,1,2,3)$
as a tuple of equations (\ref{P2}), (\ref{P21}), (\ref{P22}) and (\ref{P23}). 
Coordinate transformations between inhomogeneous coordinates are given by
\begin{eqnarray}
\left\{ \begin{array}{rrrr}
x = & \varepsilon _1^{-2/3} = & X_2\varepsilon _2^{-2/3} =& X_3\varepsilon _3^{-2/3} \\
y = & Y_1\varepsilon _1^{-1/3} = & \varepsilon _2^{-1/3}=& Y_3\varepsilon _3^{-1/3}\\
z = & Z_1\varepsilon_1 ^{-2/3} = & Z_2\varepsilon _2^{-2/3}=& \varepsilon _3^{-2/3}.
\end{array} \right.
\label{3-17}
\end{eqnarray} 
Due to the orbifold structure of $\C P^3(2,1,2,3)$, local charts $(Y_1, Z_1, \varepsilon _1)$
and $(X_3, Y_3, \varepsilon _3)$ should be divided by the $\Z_2$ actions defined by
\begin{eqnarray}
& & (Y_1, Z_1, \varepsilon _1) \mapsto (-Y_1, Z_1, - \varepsilon _1), \\
& & (X_3, Y_3, \varepsilon _3) \mapsto (X_3, -Y_3, -\varepsilon _3). \label{3-19}
\end{eqnarray} 
For our purposes, it is convenient to rewrite Eqs.(\ref{P21}), (\ref{P22}) and (\ref{P23}) 
as 3-dim vector fields (autonomous ODEs) given by
\begin{equation}
\left\{ \begin{array}{l}
\displaystyle \dot{Y}_1 = -2 + Y_1(2Y_1^3 + Y_1Z_1 + \alpha \varepsilon _1),  \\
\displaystyle \dot{Z}_1 = -2\varepsilon _1 + 2Z_1(2Y_1^3 + Y_1Z_1 + \alpha \varepsilon _1),  \\
\displaystyle \dot{\varepsilon }_1 =3\varepsilon _1 (2Y_1^3 + Y_1Z_1 + \alpha \varepsilon _1),
\end{array} \right.
\label{3-20}
\end{equation}
\begin{equation}
\left\{ \begin{array}{l}
\displaystyle \dot{X}_2 = 2X_2 - (2+Z_2 + \alpha \varepsilon _2)/X_2,  \\
\displaystyle \dot{Z}_2 = 2Z_2 - \varepsilon _2/X_2,  \\
\displaystyle \dot{\varepsilon }_2 = 3\varepsilon _2,
\end{array} \right.
\label{3-21}
\end{equation}
\begin{equation}
\left\{ \begin{array}{l}
\displaystyle \dot{X}_3 =4Y_3^3 + 2Y_3 + 2\alpha \varepsilon _3 - 2X_3\varepsilon _3,  \\
\displaystyle \dot{Y}_3 = 2X_3 - Y_3\varepsilon _3,  \\
\displaystyle \dot{\varepsilon }_3 = -3\varepsilon _3^2.
\end{array} \right.
\label{3-22}
\end{equation}
The set $\C P^2 (2,1,2)$ expressed as $\{ \varepsilon _1 = 0\} \cup \{ \varepsilon _2 =0 \} \cup \{ \varepsilon _3 = 0\}$
is an invariant manifold of the vector fields.
The dynamics on the invariant manifold describes the behavior of ($\text{P}_\text{II}$) near infinity.
Vector fields (\ref{3-20}),(\ref{3-21}),(\ref{3-22}) have four fixed points on the ``infinity set" $\C P^2 (2,1,2)$;
\\[0.2cm]
\textbf{(I).} $(X_2, Z_2, \varepsilon _2) = (\pm 1,0,0)$.

We will show later that these fixed points correspond to movable singularities of ($\text{P}_\text{II}$).
\\[0.2cm]
\textbf{(II).} $(X_2, Z_2, \varepsilon _2) = (0,-2,0)$ and $(X_3, Y_3, \varepsilon _3) = (0, 0 ,0)$.

In this case, it is easy to see from (\ref{3-17}) that $z=\infty$.
Thus, these fixed points correspond to the irregular singular point of ($\text{P}_\text{II}$).

Note that other fixed points represent the same point as one of the above.
For example, the fixed point $(Y_1, Z_1, \varepsilon _1) = (1^{1/4},0,0)$ is the same as 
$(X_2, Z_2, \varepsilon _2) = (\pm 1,0,0)$ due to the transformation (\ref{3-17}).

At first, we show that ($\text{P}_\text{II}$) is locally transformed into a linear system.
Putting $\hat{X}_2 = X_2 \pm 1$ for Eq.(\ref{3-21}) yields
\begin{equation}
\left\{ \begin{array}{l}
\displaystyle \dot{\hat{X}}_2 = 4\hat{X}_2 \pm Z_2 \pm \alpha \varepsilon _2 
   - \frac{2\hat{X}_2^2\pm \hat{X}_2Z_2 \pm \alpha \hat{X}_2\varepsilon _2}{\hat{X}_2 \mp 1},  \\
\displaystyle \dot{Z}_2 = 2Z_2 \pm \varepsilon _2 \mp \frac{\hat{X}_2\varepsilon _2}{\hat{X}_2 \mp 1},  \\
\displaystyle \dot{\varepsilon }_2 = 3\varepsilon _2.
\end{array} \right.
\label{3-23}
\end{equation}
The origin is a fixed point with the Jacobi matrix
\begin{eqnarray}
J = \left(
\begin{array}{@{\,}ccc@{\,}}
4 & \pm 1 & \pm \alpha  \\
0 & 2 & \pm 1 \\
0 & 0 & 3
\end{array}
\right).
\label{3-24}
\end{eqnarray}
Now we apply the normal form theory \cite{Cho} to the fixed point.
The eigenvalues $\lambda _1 = 4, \lambda _2 = 2, \lambda _3 = 3$ satisfy the resonance relation 
$2\lambda _2 = \lambda _1$.
However, Eq.(\ref{3-23}) does not include the corresponding resonance term.
Hence, Poincar\'{e}'s theorem on normal forms proves that 
there exists a neighborhood $U$ of the origin and a local analytic transformation defined on $U$ of the form
\begin{eqnarray*}
\left(
\begin{array}{@{\,}c@{\,}}
\hat{X}_2 \\
Z_2
\end{array}
\right) \mapsto
\left(
\begin{array}{@{\,}c@{\,}}
\hat{u}_2 \\
v_2
\end{array}
\right) = \left(
\begin{array}{@{\,}c@{\,}}
\hat{X}_2 + \varphi _1(\hat{X}_2,Z_2,\varepsilon _2) \\
Z_2 + \hat{\varphi} _2 (\hat{X}_2, Z_2, \varepsilon _2)
\end{array}
\right),
\end{eqnarray*}
such that Eq.(\ref{3-23}) is linearized as
\begin{equation}
\left\{ \begin{array}{l}
\displaystyle \dot{\hat{u}}_2 = 4\hat{u}_2 \pm v_2 \pm \alpha \varepsilon _2,  \\
\displaystyle \dot{v}_2 = 2v_2 \pm \varepsilon _2,  \\
\displaystyle \dot{\varepsilon }_2 = 3\varepsilon _2,
\end{array} \right.
\label{3-25}
\end{equation}
where local analytic functions $\varphi_1$ and $\hat{\varphi}_2$ satisfy
$\varphi_1, \hat{\varphi}_2 \sim O(|| \bm{X} ||^2)$, $\bm{X}:= (\hat{X}_2, Z_2, \varepsilon _2)$.
Note that we need not change $\varepsilon _2$ because the equation of $\varepsilon _2$ is already linear.
Furthermore, we have $\hat{\varphi}_2 (\hat{X}_2, Z_2, 0) = 0$ because the equation of $Z_2$ is linear 
when $\varepsilon _2  =0$.
Thus, we can set $\hat{\varphi}_2 = \varepsilon _2 \varphi_2$ and the above transformation takes the form
\begin{eqnarray}
\left(
\begin{array}{@{\,}c@{\,}}
\hat{X}_2 \\
Z_2
\end{array}
\right) \mapsto
\left(
\begin{array}{@{\,}c@{\,}}
\hat{u}_2 \\
v_2
\end{array}
\right) = \left(
\begin{array}{@{\,}c@{\,}}
\hat{X}_2 + \varphi _1(\hat{X}_2,Z_2,\varepsilon _2) \\
Z_2 + \varepsilon _2 \varphi _2 (\hat{X}_2, Z_2, \varepsilon _2)
\end{array}
\right),
\label{3-26}
\end{eqnarray}
with $\varphi _1 \sim O(|| \bm{X} ||^2)$ and $\varphi _2 \sim O(X_2, Z_2, \varepsilon _2)$.
Although $\varphi _1, \varphi _2$ and the neighborhood $U$ depend on the choice of the sign of
$\hat{X}_2 = X_2 \pm 1$, we need not distinguish them in this subsection.

For the linear system (\ref{3-25}), let us move to 
the original chart for ($\text{P}_\text{II}$); that is,
change coordinates as $\hat{u}_2 = u_2 \pm 1$ 
and $\tilde{x} = u_2\varepsilon _2^{-2/3}, \tilde{y} = \varepsilon _2^{-1/3}, \tilde{z} = v_2 \varepsilon _2^{-2/3}$.
Then, we obtain the solvable system
\begin{eqnarray}
\left\{ \begin{array}{l}
\displaystyle \frac{d\tilde{x}}{d\tilde{z}} = \pm 2\tilde{x}\tilde{y} + 4\tilde{y}^3 + \tilde{z}\tilde{y} + \alpha,  \\[0.2cm]
\displaystyle \frac{d\tilde{y}}{d\tilde{z}} = \mp \tilde{y}^2.  \\
\end{array} \right.
\label{3-27}
\end{eqnarray}
The coordinate transformation is given by
\begin{eqnarray*}
\left(
\begin{array}{@{\,}c@{\,}}
\tilde{x} \\
\tilde{y} \\
\tilde{z}
\end{array}
\right) = \left(
\begin{array}{@{\,}c@{\,}}
x + y^2\varphi_1 (xy^{-2} \pm 1, zy^{-2}, y^{-3})  \\
y \\
z + y^{-1} \varphi_2 (xy^{-2} \pm 1, zy^{-2}, y^{-3})
\end{array}
\right).
\end{eqnarray*}
Hence, we have obtained the next proposition.
\\[0.2cm]
\textbf{Proposition.\thedef.}
There exists a local analytic transformation $(x, z) \mapsto (\tilde{x},\tilde{z})$ defined near
$(X_2, Z_2, \varepsilon _2) = (\pm 1, 0,0)$ such that ($\text{P}_\text{II}$) 
is transformed into the solvable system (\ref{3-27}).

For ($\text{P}_\text{II}$) $d^2y/dz^2 = 2y^3 + yz + \alpha $, if $y$ becomes sufficiently large
for a finite $z$, we expect that the equation is well approximated by $d^2y/dz^2 \sim 2y^3$ .
Indeed, the second equation of (\ref{3-27}) provides $d^2\tilde{y}/d\tilde{z}^2 = 2\tilde{y}^3$.
\\

Since Eq.(\ref{3-25}) is linear, we can construct two integrals explicitly as
\begin{eqnarray*}
C_1 = \varepsilon _2^{-2/3}v_2 \mp \varepsilon _2^{1/3}, \quad 
C_2 =\varepsilon _2^{-4/3}\hat{u}_2 \pm \frac{1}{2}\varepsilon _2^{-4/3}v_2 
      +  \frac{1}{2}\varepsilon _2^{-1/3} \pm \alpha \varepsilon _2^{-1/3}.
\end{eqnarray*}
By applying the transformations (\ref{3-26}), $\hat{X}_2 = X_2 \pm 1$ and (\ref{3-17}), 
we obtain the local integrals of ($\text{P}_\text{II}$) of the form
\begin{eqnarray}
\left\{ \begin{array}{ll}
\displaystyle C_1 = z \mp y^{-1} + y^{-1} \varphi _2(xy^{-2}\pm 1, zy^{-2}, y^{-3}),  \\[0.2cm]
\displaystyle C_2 = \left( \frac{1}{2} \pm \alpha \right) y + xy^{2} \pm \frac{1}{2}y^2 z \pm y^4
       + y^4 \varphi _1(\cdots ) \pm \frac{1}{2}y \varphi _2(\cdots ). \\
\end{array} \right.
\label{3-28}
\end{eqnarray}
Arguments of $\varphi _1$ and $\varphi _2$ in the second line are the same as that of the first line.
The next theorem is proved by the same way as Thm.3.2 for ($\text{P}_\text{I}$). 
\\[0.2cm]
\textbf{Theorem.\thedef.} Any solutions of ($\text{P}_\text{II}$) are meromorphic on $\C$.
\\[0.2cm]
\textbf{Proof.} 
Fix a solution $(x(z), y(z))$ of ($\text{P}_\text{II}$) and take a disk $B(z_0,R)$ as in the proof of Thm.3.2.
Let $z_* \, (\neq \infty)$ be a singularity on the boundary of the disk.
The next lemma implies that the fixed point $(X_2, Z_2, \varepsilon _2) = (\pm 1,0,0)$
corresponds to the singularity $z_*$.
\\[0.2cm]
\textbf{Lemma.\thedef.} $(X_2, Z_2, \varepsilon _2) \to (1,0,0)$ or $\to (-1,0,0)$ 
as $z \to z_*$ along a curve $\Gamma $ inside the disk $B(z_0, R)$.
\\[0.2cm]
\textbf{Proof.} 
The same argument as the proof of Lemma 3.3 proves that either $x$ or $y$ diverges as $z \to z_*$.

(i) Suppose that $y \to \infty$ as $z \to z_*$.
We move to the $(X_2, Z_2, \varepsilon _2)$-coordinates.
Eq.(\ref{3-17}) provides
\begin{equation}
X_2 = xy^{-2}, \quad Z_2 = zy^{-2}, \quad \varepsilon _2 = y^{-3}.
\end{equation}
This immediately yields $Z_2 \to 0,\, \varepsilon _2 \to 0$ as $z \to z_*$.
Let us show $X_2 \to \pm 1$.
($\text{P}_\text{II}$) is a Hamiltonian system with the Hamiltonian function
$H = x^2/2 - y^4/2 - zy^2/2 - \alpha y$.
Thus, the equality $H = -\int\! y(z)^2dz/2$ holds along a solution.
In the $(X_2, Z_2, \varepsilon _2)$-coordinates, this equality is written as
\begin{eqnarray*}
X_2(Z_2)^2 -1 = Z_2  + 2\alpha \varepsilon_2 (Z_2)- \varepsilon _2(Z_2)^{4/3} \int^{Z_2}_{\xi}\! \varepsilon _2(z)^{-2/3}dz,
\end{eqnarray*}
where $(X_2(Z_2), \varepsilon _2(Z_2))$ is a solution of the ODE solved as a function of $Z_2$,
and $\xi$ is a certain nonzero number determined by the initial condition.
Since $Z_2 \to 0$ and $\varepsilon _2 \to 0$ as $z \to z_*$, we obtain $X^2_2 \to 1$.

(ii) Suppose that $x \to \infty$ as $z \to z_*$.
In this case, we use the $(Y_1, Z_1, \varepsilon _1)$-coordinates given by
\begin{equation}
Y_1 = yx^{-1/2}, \quad Z_1 = zx^{-1}, \quad \varepsilon _1 = x^{-3/2}.
\end{equation}
By the assumption, we have $Z_1 \to 0,\, \varepsilon _1 \to 0$ as $z \to z_*$.
Then, we can show that $Y_1^4 \to 1$ as $z \to z_*$ by the same way as above.
This means that $(Y_1, Z_1, \varepsilon _1)$ converges to the fixed point $(1^{1/4},0,0)$
of the vector field (\ref{3-20}).
It is easy to verify that this fixed point is the same point as $(X_2, Z_2, \varepsilon _2) = (\pm 1,0,0)$
if written in the $(X_2, Z_2, \varepsilon _2)$-coordinates. $\Box$

Due to the above lemma, when $z$ is sufficiently close to $z_*$,
the solution is included in the neighborhood $U$, on which local holomorphic functions 
$\varphi _1$ and $\varphi _2$ are well defined.
Then, the integrals (\ref{3-28}) are available.
To apply the implicit function theorem, put
\begin{eqnarray}
\left(
\begin{array}{@{\,}c@{\,}}
w \\
u
\end{array}
\right) = \left(
\begin{array}{@{\,}c@{\,}}
\mp y^{-1} \\
\left( \frac{1}{2} \pm \alpha \right) y + xy^{2} \pm \frac{1}{2}zy^2 \pm y^4
\end{array}
\right),
\end{eqnarray}
Then, (\ref{3-28}) takes the form
\begin{eqnarray*}
\left\{ \begin{array}{l}
\displaystyle  C_1 = z + w \mp w \varphi _2(uw^4 \pm ( \frac{1}{2} \pm \alpha )w^3 \mp \frac{1}{2}zw^2, z w^2, \mp w^3),\\
\displaystyle  C_2 = u + w^{-4}\varphi _1 (\cdots ) - \frac{1}{2}w^{-1}\varphi _2 (\cdots ). \\
\end{array} \right.
\end{eqnarray*}
Note that $C_1 = z+ w + O(w^3) $ and $C_2 = u + O(1)$ as $w \to 0$ because
$\varphi_1 \sim O(|| \bm{X} ||^2)$ and $\varphi_2 \sim O(\bm{X})$.
Let $a$ be a constant such that 
\begin{eqnarray*}
 w^{-4}\varphi _1 (uw^4 \pm ( \frac{1}{2} \pm \alpha )w^3 \mp \frac{1}{2}zw^2, z w^2, \mp w^3 )
 = az^2 + O(w),
\end{eqnarray*}
so that $C_2 = u + az^2 + O(w)$ as $w \to 0$.
Since $w\to 0$ as $z \to z_*$, the constant $C_1 = z_*$ is just the position of the singularity.
If we set
\begin{eqnarray*}
& & f_1(w,u, z) = z + w \mp w\varphi _2 (\cdots ) - z_* \\
& & f_2(w,u, z) = u + w^{-4}\varphi _1 (\cdots )- \frac{1}{2}w^{-1}\varphi _2 (\cdots ) -C_2,
\end{eqnarray*}
then $f_i(0,C_2-az_*^2, z_*) = 0$.
The Jacobi matrix of $(f_1, f_2)$ with respect to $(w,u)$ at $(w,u,z) = (0, C_2-az_*^2, z_*)$ is 
\begin{eqnarray*}
\left(
\begin{array}{@{\,}cc@{\,}}
1 & 0 \\
* & 1
\end{array}
\right).
\end{eqnarray*}
Hence, the implicit function theorem proves that there exists a local holomorphic
function $g(z)$ such that $f_1 = f_2 = 0$ is solved as $w = g(z) \sim O(z - z_*)$.
Since $y = \mp w^{-1}$, $z = z_*$ is a pole of first order of $y$.
This completes the proof of Thm.3.5. $\Box$

%%%%%%%%%%%%%%%%%%%%%%%%%%%%%%%%%%%%%%%%%%%%%%%%%%%%%%%%%%%%%%%%%%%%%%%%%%%%%%%%%%%%%%%%%%%%%%%%%%

\subsection{The fourth Painlev\'{e} equation}

($\text{P}_\text{IV}$) is given on the weighted projective space $\C P^3(1,1,1,2)$
as a tuple of equations (\ref{P4}), (\ref{P41}), (\ref{P42}) and (\ref{P43}). 
Coordinate transformations between inhomogeneous coordinates are given by
\begin{eqnarray}
\left\{ \begin{array}{rrrr}
x = & \varepsilon _1^{-1/2} = & X_2\varepsilon _2^{-1/2} =& X_3\varepsilon _3^{-1/2} \\
y = & Y_1\varepsilon _1^{-1/2} = & \varepsilon _2^{-1/2}=& Y_3\varepsilon _3^{-1/2}\\
z = & Z_1\varepsilon_1 ^{-1/2} = & Z_2\varepsilon _2^{-1/2}=& \varepsilon _3^{-1/2}.
\end{array} \right.
\label{3-32}
\end{eqnarray} 
We rewrite Eqs.(\ref{P41}), (\ref{P42}) and (\ref{P43}) as three-dimensional polynomial vector fields as before, which results in
\begin{equation}
\left\{ \begin{array}{l}
\displaystyle \dot{Y}_1 = 3Y_1 - 2\kappa_0\varepsilon _1 - Y_1^2 - 2Y_1Z_1 - Y_1(2Y_1 + 2Z_1 - 2\theta _\infty\varepsilon _1),  \\
\displaystyle \dot{Z}_1 = Z_1 + \varepsilon _1 - Z_1(2Y_1 + 2Z_1 - 2\theta _\infty\varepsilon _1),  \\
\displaystyle \dot{\varepsilon }_1 =2\varepsilon _1  - 2\varepsilon _1 (2Y_1 + 2Z_1 - 2\theta _\infty\varepsilon _1),
\end{array} \right.
\label{3-33}
\end{equation}
\begin{equation}
\left\{ \begin{array}{l}
\displaystyle \dot{X}_2 = 3X_2 -2\theta _\infty \varepsilon _2 -X_2^2 + 2X_2Z_2 - X_2( 2X_2 - 2Z_2 - 2\kappa_0\varepsilon _2),  \\
\displaystyle \dot{Z}_2 = Z_2 + \varepsilon _2 - Z_2 (2X_2 - 2Z_2 - 2\kappa_0\varepsilon _2),  \\
\displaystyle \dot{\varepsilon }_2 = 2\varepsilon _2- 2\varepsilon _2 (2X_2 - 2Z_2  - 2\kappa_0\varepsilon _2) ,
\end{array} \right.
\label{3-34}
\end{equation}
\begin{equation}
\left\{ \begin{array}{l}
\displaystyle \dot{X}_3 = -X_3^2+2X_3Y_3+2X_3-2\theta _\infty \varepsilon _3-X_3\varepsilon _3,  \\
\displaystyle \dot{Y}_3 = -Y_3^2+2X_3Y_3-2Y_3-2\kappa_0 \varepsilon _3-Y_3\varepsilon _3,  \\
\displaystyle \dot{\varepsilon }_3 = -2\varepsilon _3^2 .
\end{array} \right.
\label{3-35}
\end{equation}
These vector fields have seven fixed points on the ``infinity set" $\C P^2 (1,1,1)$;
\\[0.2cm]
\textbf{(I).} $(Y_1, Z_1, \varepsilon _1) = (0,0,0), (1,0,0)$ and $(X_2, Z_2, \varepsilon _2) = (0,0,0)$.
\\
\textbf{(II).} $(X_3, Y_3, \varepsilon _3) = (0, 0 ,0), (0,-2,0), (2,0,0)$ and $(2/3, -2/3, 0)$.

As in the case of ($\text{P}_\text{II}$), three fixed points in type (I)
correspond to movable singularities of ($\text{P}_\text{IV}$), 
and four fixed points in type (II) correspond to the irregular singular point $z=\infty$ of ($\text{P}_\text{IV}$).
Note that other fixed points represent the same point as one of the above.
For example, the fixed point $(X_2, Z_2, \varepsilon _2) = (1,0,0)$ is the same as 
$(Y_1, Z_1, \varepsilon _1) = (1,0,0)$.

By the same way as for ($\text{P}_\text{I}$) and ($\text{P}_\text{II}$),
we can show that ($\text{P}_\text{IV}$) is locally transformed into a solvable system,
and that any solutions of ($\text{P}_\text{IV}$) are meromorphic.
Suppose that a solution $(x(z), y(z))$ has a singularity $z = z_*$.
As in the proof of Lemma.3.6, we can show that as $z \to z_*$, a solution converges to
one of the fixed points listed in type (I) above (the proof is the same as Lemma 3.6 and omitted).
Indeed, the Laurent series (i),(ii),(iii) given in Sec.2.3 correspond to the fixed points 
$(Y_1, Z_1, \varepsilon _1) = (0,0,0), (1,0,0)$ and $(X_2, Z_2, \varepsilon _2) = (0,0,0)$, respectively.
Local analysis around these fixed points using the normal form theory is done in the same way as before.

Let us consider the fixed point $(Y_1, Z_1, \varepsilon _1) = (0,0,0)$.
The Jacobi matrix of the vector field at the origin is given by
\begin{equation}
J = \left(
\begin{array}{@{\,}ccc@{\,}}
3 & 0 & -2\kappa_0 \\
0& 1& 1 \\
0 & 0& 2
\end{array}
\right).
\end{equation}
We can confirm that the vector field does not have resonance terms.
Hence, due to Poincar\'{e}'s theorem, 
there exists a neighborhood $U$ of the origin and a local analytic transformation defined on $U$ of the form
\begin{eqnarray}
\left(
\begin{array}{@{\,}c@{\,}}
Y_1 \\
Z_1
\end{array}
\right) \mapsto
\left(
\begin{array}{@{\,}c@{\,}}
u_1 \\
v_1
\end{array}
\right) = \left(
\begin{array}{@{\,}c@{\,}}
Y_1 + \varphi _1(Y_1,Z_1,\varepsilon _1) \\
Z_1 + \varepsilon _1 \varphi _2 (Y_1, Z_1, \varepsilon _1)
\end{array}
\right),
\label{3-38}
\end{eqnarray}
such that Eq.(\ref{3-33}) is linearized as
\begin{equation}
\left\{ \begin{array}{l}
\displaystyle \dot{u}_1 = 3u_1 -2\kappa_0 \varepsilon _1,  \\
\displaystyle \dot{v}_1 = v_1 + \varepsilon _1,  \\
\displaystyle \dot{\varepsilon }_1 = 2\varepsilon _1,
\end{array} \right.
\label{3-39}
\end{equation}
with $\varphi _1 \sim O(|| \bm{Y} ||^2)$ and $\varphi _2 \sim O(\bm{Y} )$, 
where $\bm{Y} = (Y_1, Z_1, \varepsilon _1)$.
For this linear system, we move to the original chart for ($\text{P}_\text{IV}$) by
$\tilde{x} = \varepsilon _1^{-1/2}, \tilde{y} = u_1 \varepsilon _1^{-1/2}, \tilde{z} = v_1 \varepsilon _2^{-1/2}$.
Then, we obtain the solvable system
\begin{eqnarray}
\left\{ \begin{array}{l}
\displaystyle \frac{d\tilde{x}}{d\tilde{z}} = -\tilde{x}^2,  \\[0.2cm]
\displaystyle \frac{d\tilde{y}}{d\tilde{z}} = 2\tilde{x}\tilde{y} - 2\kappa_0.  \\
\end{array} \right.
\label{3-40}
\end{eqnarray}
The coordinate transformation is given by
\begin{eqnarray*}
\left(
\begin{array}{@{\,}c@{\,}}
\tilde{x} \\
\tilde{y} \\
\tilde{z}
\end{array}
\right) = \left(
\begin{array}{@{\,}c@{\,}}
x \\
y + x \varphi_1(x^{-1}y, x^{-1}z, x^{-2}) \\
z + x^{-1} \varphi_2(x^{-1}y, x^{-1}z, x^{-2}) 
\end{array}
\right).
\end{eqnarray*}
Hence, we have obtained the next proposition.
\\[0.2cm]
\textbf{Proposition.\thedef.}
There exists a local analytic transformation $(y, z) \mapsto (\tilde{y},\tilde{z})$ defined near
$(Y_1, Z_1, \varepsilon _1) = (0, 0,0)$ such that ($\text{P}_\text{IV}$) 
is transformed into the integrable system (\ref{3-40}).
\\

Since Eq.(\ref{3-39}) or (\ref{3-40}) is solvable, 
we can construct two local analytic integrals of ($\text{P}_\text{IV}$).
By applying the implicit function theorem for them, it is proved that if
a solution satisfies $(Y_1, Z_1, \varepsilon _1) \to (0,0,0)$ as $z \to z_*$,
then $z_*$ is a pole of first order.
The same procedure can be done for the other fixed points $(Y_1, Z_1, \varepsilon _1) = (1,0,0)$ and 
$(X_2, Z_2, \varepsilon _2) = (0,0,0)$ to prove that
\\[0.2cm]
\textbf{Theorem.\thedef.} Any solutions of ($\text{P}_\text{IV}$) are meromorphic on $\C$.
\\

The detailed calculation is the same as the proof Thm.3.5 and left to the reader.

%%%%%%%%%%%%%%%%%%%%%%%%%%%%%%%%%%%%%%%%%%%%%%%%%%%%%%%%%%%%%%%%%%%%%%%%%%%%%%%%%%%%%%%%%%%%%%%%%%
\subsection{Characterization of ($\text{P}_\text{I}$)}

In order to apply the Poincar\'{e} linearization theorem to Eq.(\ref{3-8}),
eigenvalues of the Jacobi matrix (\ref{3-9}) have to satisfy certain conditions and 
the other components of the matrix are not important.
However, to prove the meromorphy of solutions, the $(2,3)$-component of the Jacobi matrix also
plays an important role.
If the $(2,3)$-component of the Jacobi matrix were zero, the function $f_1(w,u,z)$
defined in the proof of Thm.3.2 becomes $f_1 = z + w \varphi_2(\cdots ) - z_*$ (i.e. the term $w$ does not appear).
As a result, the implicit function theorem is not applicable and we can not prove Thm.3.2.
To see the geometric role of the $(2,3)$-component, let us consider the dynamical system
\begin{equation}
\left\{ \begin{array}{l}
\dot{x} = 6y^2 + z,  \\
\dot{y} = x,   \\
\dot{z} = \beta,
\end{array} \right.
\label{3-41}
\end{equation}
where $\beta \in \C$ is a constant.
When $\beta \neq 0$, this is reduced to ($\text{P}_\text{I}$) by a suitable scaling.
This system defines a family of integral curves on $\C^3$.
We regard $\C^3$ as a vector bundle;
$z$-space is a base and $(x,y)$-space is a fiber.
As long as $\beta \neq 0$, each integral curve is a local section of the bundle,
while if $\beta = 0$, integral curves are tangent to a fiber and we can not solve the system 
as a function of $z$.
Now we change the coordinates by (\ref{3-1}) and $\hat{X}_2 = X_2 - 2$.
Then, Eq.(\ref{3-41}) is brought into the system
\begin{equation}
\left\{ \begin{array}{l}
\displaystyle \dot{\hat{X}}_2 = 6\hat{X}_2 - Z_2 + \frac{-3\hat{X}_2^2 + \hat{X}_2Z_2}{2+\hat{X}_2},  \\
\displaystyle \dot{Z}_2 = 4Z_2 - \beta \varepsilon _2 + \frac{\beta \hat{X}_2\varepsilon _2}{2+\hat{X}_2},  \\
\displaystyle \dot{\varepsilon }_2 = 5\varepsilon _2.
\end{array} \right.
\end{equation}
Hence, integral curves give local sections if and only if the $(2,3)$-component of the Jacobi matrix
of the above system is not zero.

This suggests that the $(2,3)$-component is closely related with the Painlev\'{e} property.
On the $(x, y, z)$-coordinates of $\C P^3(3,2,4,5)$, give the ODE
\begin{equation}
\frac{dx}{dz} = f(x, y, z), \quad \frac{dy}{dz} = g(x, y, z),
\label{3-43}
\end{equation}
where $f$ and $g$ are holomorphic in $x,y$ and meromorphic in $z$.
We suppose that this equation defines a meromorphic ODE on $\C P^3(3,2,4,5)$.
This means that the equations expressed in the other inhomogeneous coordinates are also meromorphic.
We will show later that these equations are rational
(recall that a meromorphic function on a projective space is rational).
Thus, there are relatively prime polynomials $h_1, h_2, h_3$ such that 
the equation written in the $(X_2, Z_2, \varepsilon _2)$-coordinates is given by
$dX_2/d\varepsilon _2 = h_1(X_2, Z_2, \varepsilon _2)/h_3(X_2, Z_2, \varepsilon _2)$
and $dZ_2/d\varepsilon _2 = h_2(X_2, Z_2, \varepsilon _2)/h_3(X_2, Z_2, \varepsilon _2)$.
As before, we introduce a vector field
\begin{eqnarray*}
\left\{ \begin{array}{l}
\dot{X}_2 = h_1(X_2, Z_2, \varepsilon _2),  \\
\dot{Z}_2 = h_2(X_2, Z_2, \varepsilon _2),  \\
\dot{\varepsilon }_2 = h_3(X_2, Z_2, \varepsilon _2).
\end{array} \right.
\end{eqnarray*}
We call it the associated vector field with $dX_2/d\varepsilon _2 = h_1/h_3, dZ_2/d\varepsilon _2 = h_2/h_3$.
The next theorem shows that ($\text{P}_\text{I}$) is characterized by the (i) 
geometry of $\C P^3(3,2,4,5)$ and (ii) a local condition at a fixed point.
Note that there are infinitely many equations satisfying only the condition (i) below.
It is remarkable that the condition (ii) seems to be very weak, however, it completely determines an equation.
\\[0.2cm]
\textbf{Theorem.\thedef.} 
Consider the ODE (\ref{3-43}), where $f$ and $g$ are holomorphic in $x,y$ and meromorphic in $z$.
Suppose the following two conditions:
\\[0.2cm]
\textbf{(i)} Eq.(\ref{3-43}) defines a meromorphic ODE on $\C P^3(3,2,4,5)$.
\\
\textbf{(ii)} The associated polynomial vector field in the $(X_2, Z_2, \varepsilon _2)$-coordinates
has a fixed point of the form $(X_2, Z_2, \varepsilon _2) = (X_*, 0,0)$.
Eigenvalues and the $(2,3)$-component of the Jacobi matrix at this point are not zero.
\\[0.2cm]
Then, Eq.(\ref{3-43}) is of the form
\begin{eqnarray}
\frac{dx}{dz} = a y^2 + bz, \quad \frac{dy}{dz} = cx,
\label{3-44}
\end{eqnarray}
where $a\neq 0, c\neq 0$ and $b$ are constants.
When $b\neq 0$, this is equivalent to ($\text{P}_\text{I}$), and 
when $b=0$, this is equivalent to the integrable equation $y'' = 6y^2$.
\\[0.2cm]
\textbf{Proof.} 
At first, we show that $f$ and $g$ are polynomial in $x, y$ and rational in $z$.
In the $(X_2, Y_2, \varepsilon _2)$-coordinates, Eq.(\ref{3-43}) is written by
\begin{eqnarray*}
\left\{ \begin{array}{l}
\displaystyle \frac{dX_2}{d\varepsilon _2} = \frac{1}{5\varepsilon _2}\left( 3X_2 - 2\varepsilon _2^{1/5} 
                 \frac{f(X_2\varepsilon _2^{-3/5}, \varepsilon _2^{-2/5}, Z_2\varepsilon _2^{-4/5})}
                     {g(X_2\varepsilon _2^{-3/5}, \varepsilon _2^{-2/5}, Z_2\varepsilon _2^{-4/5})}\right), \\
\displaystyle \frac{dZ_2}{d\varepsilon _2} = \frac{1}{5\varepsilon _2}\left( 4Z_2 - 2\varepsilon _2^{2/5} 
                 \frac{1}
                     {g(X_2\varepsilon _2^{-3/5}, \varepsilon _2^{-2/5}, Z_2\varepsilon _2^{-4/5})}\right).
\end{array} \right.
\end{eqnarray*}
By the condition (i), the right hand sides are meromorphic in $\varepsilon _2$.
In particular, $f(X_2\varepsilon _2^{-3/5}, \varepsilon _2^{-2/5}, Z_2\varepsilon _2^{-4/5})$
and $g(X_2\varepsilon _2^{-3/5}, \varepsilon _2^{-2/5}, Z_2\varepsilon _2^{-4/5})$ are meromorphic 
as a function of $\varepsilon _2^{1/5}$.
On the other hand, they are obviously meromorphic in $\varepsilon _2^{-1/5}$.
Thus, they are rational in $\varepsilon _2^{1/5}$.
Therefore, the right hand sides above are rational in $\varepsilon _2$.
This implies that $\varepsilon _2^{-2/5}g$ and $\varepsilon _2^{-1/5}f$ are rational in $\varepsilon _2$
and we can set
\begin{eqnarray*}
& & \varepsilon _2^{-1/5} f(X_2\varepsilon _2^{-3/5}, \varepsilon _2^{-2/5}, Z_2\varepsilon _2^{-4/5})
 = \frac{\sum h_i^1(X_2, Z_2) \varepsilon _2 ^i}{\sum h_i^2(X_2, Z_2) \varepsilon _2 ^i}, \\
& & \varepsilon _2^{-2/5} g(X_2\varepsilon _2^{-3/5}, \varepsilon _2^{-2/5}, Z_2\varepsilon _2^{-4/5})
 = \frac{\sum h_i^3(X_2, Z_2) \varepsilon _2 ^i}{\sum h_i^4(X_2, Z_2) \varepsilon _2 ^i},
\end{eqnarray*}
where all $\sum$ are finite sums, and $h_i^1, h_i^2, h_i^3, h_i^4$ are meromorphic in $X_2, Z_2$.
In the original chart, $f$ and $g$ are written as
\begin{eqnarray}
& & f(x, y, z) = \frac{\sum h_i^1(xy^{-3/2}, zy^{-2})y^{-5i/2-1/2}}
                         {\sum h_i^2(xy^{-3/2}, zy^{-2})y^{-5i/2}}, \nonumber\\
& & g(x, y, z) = \frac{\sum h_i^3(xy^{-3/2}, zy^{-2})y^{-5i/2-1}}
                         {\sum h_i^4(xy^{-3/2}, zy^{-2})y^{-5i/2}}.
\label{3-45}
\end{eqnarray}
Next, we move to the $(X_3, Y_3, \varepsilon _3)$-coordinates.
Eq.(\ref{3-43}) is written as
\begin{eqnarray}
\left\{ \begin{array}{l}
\displaystyle \frac{dX_3}{d\varepsilon _3} = \frac{1}{5\varepsilon _2^2}\left( 3X_3\varepsilon _3 - 4\varepsilon _3^{4/5} 
                 f(X_3\varepsilon _3^{-3/5}, Y_3\varepsilon _3^{-2/5}, \varepsilon _3^{-4/5}) \right), \\
\displaystyle \frac{dY_3}{d\varepsilon _3} = \frac{1}{5\varepsilon _3^2}\left( 2Y_3\varepsilon _3 - 4\varepsilon _3^{3/5} 
                 g(X_3\varepsilon _3^{-3/5}, Y_3\varepsilon _3^{-2/5}, \varepsilon _3^{-4/5})\right).
\end{array} \right.
\label{3-46}
\end{eqnarray}
Substituting (\ref{3-45}) yields
\begin{eqnarray}
\left\{ \begin{array}{l}
\displaystyle \frac{dX_3}{d\varepsilon _3} = \frac{1}{5\varepsilon _2^2}\left( 3X_3\varepsilon _3 - 4 
                \frac{\sum h_i^1 (X_3Y_3^{-3/2}, Y_3^{-2}) Y_3^{-5i/2-1/2}\varepsilon _3^{i+1}}
                    {\sum h_i^2 (X_3Y_3^{-3/2}, Y_3^{-2}) Y_3^{-5i/2}\varepsilon _3^{i}} \right), \\
\displaystyle \frac{dY_3}{d\varepsilon _3} = \frac{1}{5\varepsilon _3^2}\left( 2Y_3\varepsilon _3 - 4 
                  \frac{\sum h_i^3 (X_3Y_3^{-3/2}, Y_3^{-2}) Y_3^{-5i/2-1}\varepsilon _3^{i+1}}
                    {\sum h_i^4 (X_3Y_3^{-3/2}, Y_3^{-2}) Y_3^{-5i/2}\varepsilon _3^{i}} \right).
\end{array} \right.
\label{3-47}
\end{eqnarray}
The right hand sides are obviously rational in $\varepsilon _3$.
By the same argument as before, they are also rational in $Y_3$.
Hence, the right hand sides of Eq.(\ref{3-46}) are rational in $Y_3, \varepsilon _3$, and we can set
\begin{eqnarray*}
& & \varepsilon _3^{4/5} f(X_3\varepsilon _3^{-3/5}, Y_3\varepsilon _3^{-2/5}, \varepsilon _3^{-4/5})
 = \frac{\sum h_{ij}^5(X_3)Y_3^i \varepsilon _3 ^j}{\sum h_{ij}^6(X_3)Y_3^i \varepsilon _3 ^j}, \\
& & \varepsilon _3^{3/5} g(X_3\varepsilon _3^{-3/5}, Y_3\varepsilon _3^{-2/5}, \varepsilon _3^{-4/5})
 = \frac{\sum h_{ij}^7(X_3)Y_3^i \varepsilon _3 ^j}{\sum h_{ij}^8(X_3)Y_3^i \varepsilon _3 ^j},
\end{eqnarray*}
where $\sum$ are finite sums and $h_{ij}^5, h_{ij}^6, h_{ij}^7, h_{ij}^8$ are meromorphic in $X_3$.

Finally, we move to the $(Y_1, Z_1, \varepsilon _1)$-coordinates.
Repeating the same procedure, we can verify that 
$\varepsilon _1^{-1/5} f(\varepsilon _1^{-3/5}, Y_1\varepsilon _1^{-2/5}, Z_1\varepsilon _1^{-4/5})$
and $\varepsilon _1^{-2/5} g(\varepsilon _1^{-3/5}, Y_1\varepsilon _1^{-2/5}, Z_1\varepsilon _1^{-4/5})$
are rational in $Y_1,Z_1, \varepsilon _1$.
This proves that $f(x, y, z)$ and $g(x, y, z)$ are rational functions.
By the assumption, they are polynomial in $x$ and $y$.

Now we can write $g$ and $f/g$ as quotients of polynomials as
\begin{equation}
g(X,Y,Z) = \frac{\sum a_{ijk} X^iY^jZ^k}{\sum b_k Z^k}, \quad
\frac{f(X,Y,Z)}{g(X,Y,Z)} = \frac{\sum p_{ijk} X^iY^jZ^k}{\sum q_{ijk} X^iY^jZ^k}.
\label{3-48}
\end{equation}
Our purpose is to determine coefficients $a_{ijk}, b_k, p_{ijk}$ and $q_{ijk}$ by the conditions (i) and (ii).
In the $(X_2, Y_2, Z_2)$-coordinates, Eq.(\ref{3-43}) with (\ref{3-48}) is given by
\begin{eqnarray}
\left\{ \begin{array}{l}
\displaystyle \frac{dX_2}{d\varepsilon _2} = \frac{1}{5\varepsilon _2}\left( 3X_2 - 2\varepsilon _2^{1/5} 
                 \frac{\sum p_{ijk} X_2^i Z_2^k \varepsilon _2^{-(3i+2j+4k)/5} }
                     {\sum q_{ijk} X_2^i Z_2^k \varepsilon _2^{-(3i+2j+4k)/5}}\right), \\
\displaystyle \frac{dZ_2}{d\varepsilon _2} = \frac{1}{5\varepsilon _2}\left( 4Z_2 - 2\varepsilon _2^{2/5} 
                 \frac{\sum b_k Z_2^k \varepsilon _2^{-4k/5}}
                     {\sum a_{ijk} X_2^i Z_2^k \varepsilon _2^{-(3i+2j+4k)/5}}\right).
\end{array} \right.
\label{3-49}
\end{eqnarray}
Due to the condition (i), the right hand sides are rational in $\varepsilon _2$.
This yields the conditions for coefficients as
\begin{equation}
\left\{ \begin{array}{l}
p_{ijk} \neq 0 \quad \text{only when} \quad 3i+2j+4k-1 = 5m + \delta, \quad (m=-1,0,\cdots , M),  \\
q_{ijk} \neq 0 \quad \text{only when} \quad 3i+2j+4k = 5m' + \delta, \quad (m'=0,1,\cdots , M'),  \\
b_{k} \neq 0   \quad \text{only when} \quad 4k-2 = 5n + \delta', \quad (n=-1,0,\cdots , N), \\
a_{ijk} \neq 0 \quad \text{only when} \quad 3i+2j+4k = 5n' + \delta', \quad (n'=0,1,\cdots , N'),
\end{array} \right.
\label{3-50}
\end{equation}
where $\delta , \delta ' \in \{ 0,1,2,3,4\}$.
More precisely, the first line means that $p_{ijk} \neq 0$ only when there are integers $m$ and $\delta $
such that $(i,j,k)$ satisfies $3i+2j+4k-1 = 5m + \delta$, where $\delta \in \{ 0,1,2,3,4\}$ is independent of $(i,j,k)$.
Since Eq.(\ref{3-49}) is rational, there is the largest integer $m$ satisfying $3i+2j+4k-1 = 5m + \delta$ and $p_{ijk} \neq 0$.
In (\ref{3-50}), the largest integer is denoted by $M$.
Integers $M', N$ and $N'$ play a similar role.
Then, Eq.(\ref{3-49}) is rewritten as
\begin{eqnarray}
\left\{ \begin{array}{l}
\displaystyle \frac{dX_2}{d\varepsilon _2} = \frac{1}{5\varepsilon _2}\left( 3X_2 - 2
                 \frac{\sum p_{ijk} X_2^i Z_2^k \varepsilon _2^{-m} }
                     {\sum q_{ijk} X_2^i Z_2^k \varepsilon _2^{-m'}}\right), \\
\displaystyle \frac{dZ_2}{d\varepsilon _2} = \frac{1}{5\varepsilon _2}\left( 4Z_2 - 2 
                 \frac{\sum b_k Z_2^k \varepsilon _2^{-n}}
                     {\sum a_{ijk} X_2^i Z_2^k \varepsilon _2^{-n'}}\right).
\end{array} \right.
\label{3-51}
\end{eqnarray}
In order to confirm the condition (ii), we shall rewrite it as a polynomial vector field.

\textbf{(I).} When $M>M'$ and $N>N'$, the associated polynomial vector field is of the form
\begin{equation}
\left\{ \begin{array}{l}
\dot{X}_2 
= 3X_2 (\sum a_{ijk} X_2^iZ_2^k \varepsilon _2^{N-n'}) (\sum q_{ijk} X_2^iZ_2^k \varepsilon _2^{M-m'}) \\
\qquad \qquad \qquad  -2 (\sum p_{ijk} X_2^iZ_2^k \varepsilon _2^{M-m}) (\sum a_{ijk} X_2^iZ_2^k \varepsilon _2^{N-n'}), \\[0.4cm]
\dot{Z}_2 
= 4Z_2 (\sum a_{ijk} X_2^iZ_2^k \varepsilon _2^{N-n'}) (\sum q_{ijk} X_2^iZ_2^k \varepsilon _2^{M-m'}) \\
\qquad \qquad \qquad  -2 (\sum b_{k} Z_2^k \varepsilon _2^{N-n}) (\sum q_{ijk} X_2^iZ_2^k \varepsilon _2^{M-m'}),  \\[0.4cm]
\dot{\varepsilon }_2 
= 5\varepsilon _2 (\sum a_{ijk} X_2^iZ_2^k \varepsilon _2^{N-n'}) (\sum q_{ijk} X_2^iZ_2^k \varepsilon _2^{M-m'}).
\end{array} \right.
\end{equation}
Because of the condition (ii), we seek a fixed point of the form $(X_*,0,0)$ with nonzero eigenvalues.
Since $N-n' > 0$ and $M-m' > 0$, the right hand side of the equation of $\varepsilon _2$ is of order $O(\varepsilon _2^3)$.
Hence, the Jacobi matrix at the point $(X_*,0,0)$ has a zero eigenvalue.
In a similar manner, we can verify that two cases $M>M', N\leq N'$ and $M\leq M', N>N'$ are excluded.
In these cases, the right hand side of the equation of $\varepsilon _2$ is of order $O(\varepsilon _2^2)$
and the Jacobi matrix has a zero eigenvalue.

\textbf{(II).} When $M\leq M'$ and $N\leq N'$, the associated polynomial vector field is of the form
\begin{equation}
\left\{ \begin{array}{l}
\dot{X}_2 
= 3X_2 (\sum a_{ijk} X_2^iZ_2^k \varepsilon _2^{N'-n'}) (\sum q_{ijk} X_2^iZ_2^k \varepsilon _2^{M'-m'}) \\
\qquad \qquad \qquad  -2 (\sum p_{ijk} X_2^iZ_2^k \varepsilon _2^{M'-m}) (\sum a_{ijk} X_2^iZ_2^k \varepsilon _2^{N'-n'}), \\[0.4cm]
\dot{Z}_2 
= 4Z_2 (\sum a_{ijk} X_2^iZ_2^k \varepsilon _2^{N'-n'}) (\sum q_{ijk} X_2^iZ_2^k \varepsilon _2^{M'-m'}) \\
\qquad \qquad \qquad  -2 (\sum b_{k} Z_2^k \varepsilon _2^{N'-n}) (\sum q_{ijk} X_2^iZ_2^k \varepsilon _2^{M'-m'}),  \\[0.4cm]
\dot{\varepsilon }_2 
= 5\varepsilon _2 (\sum a_{ijk} X_2^iZ_2^k \varepsilon _2^{N'-n'}) (\sum q_{ijk} X_2^iZ_2^k \varepsilon _2^{M'-m'}).
\end{array} \right.
\label{3-53}
\end{equation}
Suppose that this has a fixed point of the form $(X_*,0,0)$.
The $(2,3)$-component of the Jacobi matrix at this point is given by
\begin{eqnarray*}
-2 \cdot \frac{\partial }{\partial \varepsilon _2}\Bigl|_{(X_*, 0,0)}
(\sum b_{k} Z_2^k \varepsilon _2^{N'-n}) (\sum q_{ijk} X_2^iZ_2^k \varepsilon _2^{M'-m'}).
\end{eqnarray*}
We require that this quantity is not zero.

\textbf{(II-a).} Suppose that the polynomial $\sum b_{k} Z_2^k \varepsilon _2^{N'-n}$ includes a constant term.
This means that $b_k\neq 0$ when $k=0$ and $n= N'$.
Substituting it to the third condition of (\ref{3-50}) provides $-2 = 5N' + \delta'$.
This proves that $N' = -1$ and $\delta ' = 3$.
Then, the fourth condition of (\ref{3-50}) yields $3i+2j+4k = -2$.
Since there are no nonnegative integers $i,j,k$ satisfying this relation,
$a_{ijk} = 0$ for any $i,j,k$.
In this case, we obtain $\dot{\varepsilon }_2  = 0$ and the Jacobi matrix has a zero eigenvalue.

\textbf{(II-b).} When $\sum b_{k} Z_2^k \varepsilon _2^{N'-n}$ does not have a constant term,
it has to include a monomial $\varepsilon _2$. Otherwise, the $(2,3)$-component of the Jacobi matrix becomes zero.
This means that $b_k\neq 0$ when $k=0$ and $n= N'-1$.
The third condition of (\ref{3-50}) provides $-2 = 5(N'-1) + \delta'$,
which proves $N' = 0$ and $\delta ' = 3$.
Since $N\leq N'$, we have $N = -1$ or $N=0$. 
In this case, (\ref{3-50}) becomes
\begin{equation}
\left\{ \begin{array}{l}
b_{k} \neq 0   \quad \text{only when} \quad 4k-2 = -2\quad \text{or}\quad 3,\\
a_{ijk} \neq 0 \quad \text{only when} \quad 3i+2j+4k = 3.
\end{array} \right.
\end{equation}
Therefore, nonzero numbers among these coefficients are only $b_0$ and $a_{100}$.
This proves that $g(X,Y,Z)$ is given by $a_{100}X/b_0$.
In what follows, we put $a_{100}/b_0 = c$.

Since $g=cX$ and $f$ is polynomial in $X,Y$, $f/g$ can be written as
\begin{eqnarray*}
\frac{f(X,Y,Z)}{g(X,Y,Z)} = \frac{\sum p_{ijk}X^iY^jZ^k}{\sum q_{10k}XZ^k}.
\end{eqnarray*}
In this case, the equation of $\varepsilon _2$ in (\ref{3-22}) is given by
\begin{equation}
\dot{\varepsilon }_2 = 5a_{100}X_2^2\varepsilon _2 \cdot \left( \sum q_{10k}Z_2^k \varepsilon _2^{M'-m'} \right).
\end{equation}
The polynomial $\sum q_{10k}Z_2^k \varepsilon _2^{M'-m'}$ has to include a constant term so that 
the Jacobi matrix at $(X_*, 0,0)$ does not have a zero eigenvalue.
This means that $q_{10k} \neq 0$ when $k = 0$ and $m' = M'$.
Thus (\ref{3-50}) provides $3 = 5M' + \delta $.
This shows $M' = 0$ and $\delta =3$.
Then, (\ref{3-50}) becomes 
\begin{equation}
\left\{ \begin{array}{l}
p_{ijk} \neq 0   \quad \text{only when} \quad 3i+2j+4k-1 = 3,\\
q_{10k} \neq 0 \quad \text{only when} \quad 3+4k = 3.
\end{array} \right.
\end{equation}
Therefore, nonzero numbers among these coefficients are only $p_{020}, p_{001}$ and $q_{100}$.
Putting $a := cp_{020}/q_{100}$ and $b:= cp_{001}/q_{100}$, we obtain the ODE (\ref{3-44})
as a necessary condition for (i) and (ii).
It is straightforward to confirm that (\ref{3-44}) actually satisfies the conditions (i) and (ii)
when $a\neq 0, c\neq 0$. This completes the proof. $\Box$

%%%%%%%%%%%%%%%%%%%%%%%%%%%%%%%%%%%%%%%%%%%%%%%%%%%%%%%%%%%%%%%%%%%%%%%%%%%%%%%%%%%%%%%%%%%%%%%%%%
%%%%%%%%%%%%%%%%%%%%%%%%%%%%%%%%%%%%%%%%%%%%%%%%%%%%%%%%%%%%%%%%%%%%%%%%%%%%%%%%%%%%%%%%%%%%%%%%%%

\section{The space of initial conditions}

In this section, we construct the spaces of initial conditions for ($\text{P}_\text{I}$),
($\text{P}_\text{II}$) and ($\text{P}_\text{IV}$)
by the weighted blow-ups of $\C P^3(p,q,r,s)$.
For a polynomial system, a manifold $E(z)$ parameterized by $z\in \C$ is called the space of initial conditions
if any solutions give global holomorphic sections on the fiber bundle $\mathcal{P} = \{ (x,z) \, | \, x\in E(z) , z\in \C\}$ over $\C$.
For the Painlev\'{e} equations, it was first constructed by Okamoto \cite{Oka} by blow-ups of a Hirzebruch surface eight times and by
removing a certain divisor called vertical leaves.
Different approaches are proposed by Duistermaat and Joshi \cite{Dui} and Iwasaki and Okada \cite{Iwa}.
They also performed blow-ups many times.
See \cite{Sak, Sai} for algebro-geometric approach.

Here, we will obtain the spaces of initial conditions by weighted blow-ups only one time
for ($\text{P}_\text{I}$), two times for ($\text{P}_\text{II}$) and three times for ($\text{P}_\text{IV}$).
These numbers are the same as the numbers of types of Laurent series given in Sec.2.3.
We will easily find  a symplectic structure of the space of initial conditions.
For ($\text{P}_\text{I}$), we will recover Painlev\'{e}'s coordinates in a purely geometric manner.
We find a symplectic structure of the space of initial conditions and show that 
Painlev\'{e}'s coordinates are the Darboux coordinates of the symplectic structure.

%%%%%%%%%%%%%%%%%%%%%%%%%%%%%%%%%%%%%%%%%%%%%%%%%%%%%%%%%%%%%%%%%%%%%%%%%%%%%%%%%%%%%%%%%%%%%%%%%%

\subsection{The first Painlev\'{e} equation}

Recall that ($\text{P}_\text{I}$) written in the $(X_2, Z_2, \varepsilon _2)$-coordinates
has two fixed points $(\pm 2,0,0)$.
Putting $\hat{X}_2 = X_2\pm 2$, we obtain
\begin{equation}
\left\{ \begin{array}{l}
\displaystyle \dot{\hat{X}}_2 = 6\hat{X}_2 \pm Z_2 + \frac{\pm 3\hat{X}_2^2 + \hat{X}_2Z_2}{2 \mp \hat{X}_2}, \\
\displaystyle \dot{Z}_2 = 4Z_2 \pm \varepsilon _2 + \frac{\hat{X}_2\varepsilon _2}{2 \mp \hat{X}_2},  \\
\displaystyle \dot{\varepsilon }_2 = 5\varepsilon _2.
\end{array} \right.
\label{4-1}
\end{equation}
(In Sec.3, we used only $\hat{X}_2 = X_2- 2$).
The origin is a fixed point of the vector field and it is a singularity of the foliation
defined by integral curves.
We apply a blow-up to this point.
At first, we change the coordinates by the linear transformation
\begin{eqnarray*}
\left\{ \begin{array}{l}
\hat{X}_2 = u \mp \frac{1}{2}v - \frac{1}{2}w, \\
Z_2 = v,  \\
\varepsilon _2 = w.
\end{array} \right.
\end{eqnarray*}
This yields
\begin{eqnarray}
\left\{ \begin{array}{l}
\displaystyle \dot{u} = 6u + f_1(u,v,w),  \\
\displaystyle \dot{v} = 4v \pm w + f_2(u,v,w),  \\
\displaystyle \dot{w} = 5w,
\end{array} \right.
\label{4-2}
\end{eqnarray}
where $f_1$ and $f_2$ denote nonlinear terms.
Note that the linear part is \textit{not} diagonalized;
since we know that the $(2,3)$-component is important (Thm.3.9),
we remove only the $(1,2)$-component of the linear part of (\ref{4-1}).
Now we introduce the weighted blow-up with weights $6,4,5$, 
which are taken from eigenvalues of the Jacobi matrix, defined by the following transformations
\begin{eqnarray}
\left\{ \begin{array}{llll}
u & = u_1^6     & = v_2^6 u_2 & = w_3^6 u_3, \\
v & = u_1^4v_1  & = v_2^4     & = w_3^4 v_3, \\
w & = u_1^5w_1  & = v_2^5 w_2 & = w_3^5.  
\end{array} \right.
\label{4-3}
\end{eqnarray}
The exceptional divisor $\{ u_1 = 0\} \cup \{ v_2 = 0\} \cup \{ w_3= 0\}$
is a 2-dim weighted projective space $\C P^2(6,4,5)$ and the blow-up of $(u,v,w)$-space
is a (singular) line bundle over $\C P^2(6,4,5)$.
We mainly use the $(u_3, v_3, w_3)$-coordinates.
In the $(u_3, v_3, w_3)$-coordinates, Eq.(\ref{4-2}) is written as
\begin{equation}
\left\{ \begin{array}{l}
\displaystyle \frac{du}{dv} = \frac{1}{8}
\left( v^2w \pm 3vw^2 + 2w^3 \mp 8uvw^3 -10uw^4 + 12 u^2w^5 \right), \\[0.2cm]
\displaystyle \frac{dw}{dv} = \frac{1}{4}\left( \pm 4 \pm vw^4 + w^5 - 2uw^6 \right),  \\
\end{array} \right.
\label{4-4}
\end{equation}
where the subscripts for $u_3, v_3, w_3$ are omitted for simplicity.
The relation between the original chart $(x, y, z)$ and $(u_3, v_3, w_3)$ is
\begin{equation}
\left\{ \begin{array}{l}
\displaystyle x = u_3w_3^3 \mp 2w_3^{-3} \mp \frac{1}{2}v_3 w_3 - \frac{1}{2} w_3^2, \\
\displaystyle y = w_3^{-2},  \\
\displaystyle z = v_3,
\end{array} \right.
\label{4-5}
\end{equation}
or
\begin{equation}
\left\{ \begin{array}{l}
\displaystyle u_3 = xy^{3/2} \pm 2y^3 \pm \frac{1}{2}zy + \frac{1}{2}y^{1/2}, \\
\displaystyle w_3 = y^{-1/2},  \\
\displaystyle v_3 = z.
\end{array} \right.
\label{4-6}
\end{equation}
It is remarkable that the independent variable $z$ is not changed despite the fact that 
$z$ is changed in each step of transformations.
Now we have recovered Painlev\'{e}'s coordinates (\ref{4-5}) which was introduced in his paper
to prove the Painlev\'{e} property of ($\text{P}_\text{I}$).
He found this transformation by observing a Laurent series of a solution, 
see also Gromak, Laine and Shimomura \cite{Gro}.
At a first glance, $\C^2_{(x,y)} \cup \C^2_{(u_3, w_3)}$ glued by (\ref{4-5}) does not define a manifold
because (\ref{4-5}) is not one-to-one but one-to-two.
Nevertheless, we can show that (\ref{4-5}) defines a certain algebraic surface.
Recall that the $(X_2, Z_2,\varepsilon _2)$-space should be divided by the 
$\Z_2$ action $(X_2, Z_2,\varepsilon _2) \mapsto (-X_2, Z_2, -\varepsilon _2)$
due to the orbifold structure of $\C P^3(3,2,4,5)$.
This action induces a certain $\Z_2$ action on the $(u_3, v_3, w_3)$-space.
If we divide the $(u_3, v_3, w_3)$-space by the $\Z_2$ action, we can prove that 
(\ref{4-5}) becomes a one-to-one mapping.
Then, $\C^2_{(x,y)}$ and $\C^2_{(u_3, w_3)}/\Z_2$ are glued by (\ref{4-5}) to define 
an algebraic surface, which gives the space of initial conditions for ($\text{P}_\text{I}$).
The space $\C^2_{(u_3, w_3)}/ \Z_2$ is realized as a nonsingular algebraic surface defined by
\begin{equation}
M(z) : V^2 = UW^4 + 2z W^3 + 4W
\label{4-7}
\end{equation}
with the parameter $z$ (the proof is given below).
By using $(U,V,W)$, the relations (\ref{4-5}),(\ref{4-6}) are rewritten as
\begin{equation}
\left\{ \begin{array}{l}
\displaystyle x = VW^{-2} - \frac{1}{2}W,  \\
\displaystyle y = W^{-1}.  \\
\end{array} \right. ,\quad 
\left\{ \begin{array}{l}
\displaystyle V = xy^{-2} + \frac{1}{2}y^{-3}, \\
\displaystyle W = y^{-1}.
\end{array} \right.
\label{4-8}
\end{equation}
Hence, $\C^2_{(x,y)}$ and $M(z)$ glued by this relation defines a nonsingular algebraic surface 
denoted by $E(z)$.
The surface $M(z)$ admits a holomorphic symplectic form
\begin{equation}
-\frac{1}{W^4}dV \wedge dW = \frac{1}{4UW^3 + 6z W^2 + 4}dV \wedge dU.
\end{equation}
We can verify that
\begin{equation}
dx \wedge dy = -\frac{1}{W^4}dV \wedge dW.
\end{equation}
Thus, $E(z)$ also has a holomorphic symplectic form.
($\text{P}_\text{I}$) written by $(V,W)$ is 
\begin{equation}
\left\{ \begin{array}{ll}
\displaystyle \frac{dV}{dz} = 6 + zW^2 + \frac{1}{4W}(W^3 + 4V)(W^3 - 2V) 
 = W^4 \frac{\partial H}{\partial W}  \\[0.4cm]
\displaystyle \frac{dW}{dz} = \frac{1}{2}W^3 - V  = -W^4 \frac{\partial H}{\partial V},  \\
\end{array} \right.
\label{4-11}
\end{equation}
where $H$ is given by
\begin{eqnarray}
H &=& \frac{1}{2}x^2 - 2y^3 - zy \nonumber \\
&=& \frac{V^2}{2W^4} - \frac{V}{2W} + \frac{1}{8}W^2 - \frac{2}{W^3} - \frac{z}{W}.
\end{eqnarray}
Hence, Eq.(\ref{4-11}) is a Hamiltonian system with respect to the symplectic form
$-W^{-4}dV \wedge dW$ as well as the original ($\text{P}_\text{I}$) written by $(x, y)$.
Let $z_*$ be a pole of a solution of ($\text{P}_\text{I}$).
As $z \to z_*$, $x, y \to \infty$ and $V,W \to 0$.
By using (\ref{4-7}), it is easy to verify that Eq.(\ref{4-11}) is holomorphic even at $W = 0$.
This proves that $E(z) = \C^2_{(x,y)} \cup M(z)$ is the desired space of initial conditions.
Note that the system (\ref{4-4}) is already a Hamiltonian system with the Hamiltonian
\begin{eqnarray}
H = \frac{1}{8} \left( \pm 8u \pm 2uvw^4 + 2uw^5 - 2u^2 w^6 - \frac{1}{2}v^2 w^2 \mp vw^3 - \frac{1}{2}w^4\right) .
\end{eqnarray}
The transformation (\ref{4-5}) yields
\begin{equation}
dx \wedge dy = -2 du_3 \wedge dw_3.
\end{equation}
This means that $(u_3, w_3)$-coordinates are the Darboux coordinates for the form 
$-W^{-4}dV \wedge dW$ (if we remove the factor $-2$ by a suitable scaling).
These results are summarized as follows.
\\[0.2cm]
\textbf{Theorem.\thedef.} 

\noindent \textbf{(i)} The space $\C^2_{(u_3, w_3)}$ divided by the $\Z_2$ action induced from the orbifold 
structure of $\C P^3(3,2,4,5)$ gives the algebraic surface $M(z)$.
The space of initial conditions $E(z)$ for ($\text{P}_\text{I}$) is given by
$\C^2_{(x, x)} \cup M(z)$ glued by (\ref{4-8}).

\noindent \textbf{(ii)} $M(z)$ and $E(z)$ have holomorphic symplectic forms and ($\text{P}_\text{I}$) 
is a Hamiltonian system with respect to the form.
Painlev\'{e}'s coordinates defined by (\ref{4-5}) are the Darboux coordinates of the symplectic form on $M(z)$.

\noindent \textbf{(iii)} Consider an ODE (\ref{3-43}) defined on the $(x, y, z)$-coordinates,
where $f$ and $g$ are polynomials in $x$ and $y$.
If it is also expressed as a polynomial ODE in the Painlev\'{e}'s coordinates, then (\ref{3-43}) is ($\text{P}_\text{I}$). 
\\

Recently, a similar result is obtained by Iwasaki and Okada \cite{Iwa} by a different approach.
Symplectic atlases for the Painlev\'{e} equations are found by Takano et al. \cite{Tak1, Tak2, Tak3}
for the second Painlev\'{e} to sixth Painlev\'{e} equations, while left open for ($\text{P}_\text{I}$).
\\[0.2cm]
\textbf{Proof.}
Due to the orbifold structure of $\C P^3(3,2,4,5)$, the $(X_2, Z_2, \varepsilon _2)$-space
should be divided by the $\Z_2$ action 
$(X_2, Z_2, \varepsilon _2) \mapsto (-X_2, Z_2, -\varepsilon _2)$.
It is straightforward to show that in the $(u_3, v_3, w_3)$-coordinates, this action is written by
\begin{equation}
\left(
\begin{array}{@{\,}c@{\,}}
u_3 \\
v_3 \\
w_3
\end{array}
\right) \mapsto \left(
\begin{array}{@{\,}c@{\,}}
-u_3 + v_3 w_3^{-2} + 4w_3^{-6} \\
v_3 \\
-w_3
\end{array}
\right) .
\end{equation}
Since $v_3$ is fixed, we consider $\C^2_{(u_3, w_3)} /\Z_2$.
Polynomial invariants of this action are generated by
\begin{equation}
\left\{ \begin{array}{l}
\displaystyle U = u_3(u_3 w_3^6 - v_3w_3^4 - 4) + \frac{1}{4}w_3^2v_3^2,  \\
\displaystyle V = w_3^7 (u_3- \frac{1}{2}v_3 w_3^{-2} - 2w_3^{-6}) ,  \\
\displaystyle W = w_3^2.
\end{array} \right.
\end{equation}
They satisfy the equation (\ref{4-7}), which proves $\C^2_{(u_3, w_3)} /\Z_2 = M(z)$.
The rest of (i) and (ii) have already been shown.
Let us prove (iii).
In what follows, we omit the subscripts for $u_3, v_3, w_3$.
By (\ref{4-5}) with the upper sign, Eq.(\ref{3-43}) is written in Painlev\'{e}'s coordinates as
\begin{eqnarray}
\frac{d}{dz}\left(
\begin{array}{@{\,}c@{\,}}
u \\
w
\end{array}
\right) &=& \left(
\begin{array}{@{\,}c@{\,}}
\displaystyle w^{-3}\cdot f + \frac{3}{2}uw^2 \cdot g - \frac{1}{4} z \cdot g + 3w^{-4} \cdot g 
 - \frac{1}{2}w \cdot g + \frac{1}{2} w^{-2} \\
\displaystyle -\frac{1}{2} w^3 \cdot g
\end{array}
\right) \nonumber \\
&=:& T(f,g) + \left(
\begin{array}{@{\,}c@{\,}}
\frac{1}{2}w^{-2} \\
0
\end{array}
\right) =: \hat{T}(f,g).
\end{eqnarray}
Our purpose is to show that if $\hat{T}(f,g)$ is a polynomial, then Eq.(\ref{3-43}) is ($\text{P}_\text{I}$).
Since we know that $\hat{T}(f,g)$ is a polynomial for ($\text{P}_\text{I}$),
it is sufficient to show the uniqueness.
We define operators $T$ and $\hat{T}$ as above.
The operator $T$ is a linear mapping from the space of polynomials
$\C [x, y] \times \C [x, y]$ into the space of Laurent polynomials
$\C [u,w,w^{-1}] \times \C [u,w,w^{-1}] $ for each $z$.
Let 
\begin{equation}
\Pi : \C [u,w,w^{-1}] \times \C [u,w,w^{-1}] \to \C [u,w^{-1}] \times \C [u,w^{-1}]
\end{equation}
be the natural projection to the principle part.
If there are two pairs of polynomials $(f_1, g_1)$ and $(f_2, g_2)$ such that 
$\hat{T}(f_i, g_i),\, (i = 1,2)$ are polynomials, then
\begin{eqnarray*}
\hat{T}(f_1, g_1) - \hat{T}(f_2, g_2) = T(f_1 - f_2, g_1- g_2)
\end{eqnarray*}
is also a polynomial and $\Pi \circ T(f_1 - f_2, g_1- g_2) = 0$.
Hence, it is sufficient to prove $\mathrm{Ker}\, \Pi \circ T = \{ 0\}$.
For this purpose, we show that images of monomials of the form 
$(X^mY^n, 0), (0, X^mY^n),\,\, (m,n = 0,1,\cdots )$ are linearly independent.
They are calculated as
\begin{eqnarray*}
T (X^mY^n, 0) &=& \left(
\begin{array}{@{\,}c@{\,}}
w^{-2n-3} (uw^3 - 2w^{-3} - \frac{1}{2} z w - \frac{1}{2} w^2)^m \\
0
\end{array}
\right), \\
T(0, X^mY^n) &=& \left(
\begin{array}{@{\,}c@{\,}}
 * \\
-\frac{1}{2}w^{-2n+3} (uw^3 - 2w^{-3} - \frac{1}{2} z w - \frac{1}{2} w^2)^m 
\end{array}
\right).
\end{eqnarray*}
It is easy to verify that the principle parts of them are linearly independent. $\Box$
\\[0.2cm]
\textbf{Remark.\thedef.}
We have constructed the space of initial conditions by the weighted blow-up at the fixed point 
$(X_2, Z_2, \varepsilon _2) = (2,0,0)$.
We can also construct it by using the fixed point $(Y_1, Z_1, \varepsilon _1) = ((1/4)^{1/3} , 0,0)$, 
which represents the same point as $(X_2, Z_2, \varepsilon _2) = (2,0,0)$.
By the same procedure as before (an affine transformation and the weighted blow-up
in the $(Y_1, Z_1, \varepsilon _1)$-coordinates), we obtain the one-to-three transformation
\begin{equation}
\left\{ \begin{array}{l}
\displaystyle x = w_3^{-3},  \\
\displaystyle y = u_3 w_3^4 - \frac{2^{-1/3}}{3}v_3 w_3^2 + 2^{-2/3} w_3^{-2},  \\
\displaystyle z = v_3.
\end{array} \right.
\label{4-19}
\end{equation}
Due to the orbifold structure, the $(Y_1, Z_1, \varepsilon _1)$-space should be divided by the 
$\Z_3$ action (\ref{3-2}).
This induces the $\Z_3$ action in the $(u_3, v_3, w_3)$-space and we can show that $\C^2_{(x, y)}$
and $\C^2_{(u_3, w_3)}/\Z_3$ glued by (\ref{4-19}) gives the same algebraic surface $E(z)$ as before.

%%%%%%%%%%%%%%%%%%%%%%%%%%%%%%%%%%%%%%%%%%%%%%%%%%%%%%%%%%%%%%%%%%%%%%%%%%%%%%%%%%%%%%%%%%%%%%%%%%

\subsection{The second Painlev\'{e} equation}

Recall that ($\text{P}_\text{II}$) written in the $(X_2, Z_2, \varepsilon _2)$-coordinates
has two fixed points $(\pm 1,0,0)$.
Putting $\hat{X}_2 = X_2\pm 1$, we have obtained Eq.(\ref{3-23}).
The origin is a fixed point of (\ref{3-23}) and it is a singularity of the foliation
defined by integral curves.
We apply a blow-up to this point.
At first, we change the coordinates by the linear transformation
\begin{eqnarray}
\left\{ \begin{array}{l}
\hat{X}_2 = u \mp \frac{1}{2}v - \left( \frac{1}{2} \pm \alpha \right)w, \\
Z_2 = v,  \\
\varepsilon _2 = w.
\end{array} \right.
\label{4-20}
\end{eqnarray}
Then, we obtain
\begin{eqnarray}
\left\{ \begin{array}{l}
\displaystyle \dot{u} = 4u + f_1(u,v,w),  \\
\displaystyle \dot{v} = 2v \pm w + f_2(u,v,w),  \\
\displaystyle \dot{w} = 3w,
\end{array} \right.
\label{4-21}
\end{eqnarray}
where $f_1$ and $f_2$ denote nonlinear terms.
Now we introduce the weighted blow-up with weights $4,2,3$, 
which are taken from eigenvalues of the Jacobi matrix, defined by the following transformations
\begin{eqnarray}
\left\{ \begin{array}{llll}
u & = u_1^4     & = v_2^4 u_2 & = w_3^4 u_3, \\
v & = u_1^2v_1  & = v_2^2     & = w_3^2 v_3, \\
w & = u_1^3w_1  & = v_2^3 w_2 & = w_3^3.  
\end{array} \right.
\label{4-22}
\end{eqnarray}
The exceptional divisor $\{ u_1 = 0\} \cup \{ v_2 = 0\} \cup \{ w_3= 0\}$
is a 2-dim weighted projective space $\C P^2(4,2,3)$ and the blow-up of $(u,v,w)$-space
is a (singular) line bundle over $\C P^2(4,2,3)$.
Note that we performed the blow-ups at two points $(X_2, Z_2, \varepsilon _2) = (1,0,0)$ and $(-1,0,0)$.
If we want to distinguish the sign of $\hat{X}_2 = X_2 \pm 1$, the notation
$(u^{\pm}_i,v^{\pm}_i,w^{\pm}_i)$ for $i=1,2,3$ will be used.
In the $(u_3, v_3, w_3)$-coordinates, Eq.(\ref{4-21}) is written as
\begin{equation}
\left\{ \begin{array}{l}
\displaystyle \frac{du}{dv} = \frac{1}{4}
\left( 4uw-1\mp 2\alpha \right) \left( \mp v+2uw^2-(1\pm 2\alpha ) w\right)
 = -\frac{\partial \widetilde{H}}{\partial w}, \\[0.2cm]
\displaystyle \frac{dw}{dv} = \frac{1}{2}\left( -2uw^4+(1\pm 2\alpha ) w^3\pm vw^2 \pm 2 \right)
 = \frac{\partial \widetilde{H}}{\partial u},  \\
\end{array} \right.
\label{4-23}
\end{equation}
where the subscripts for $u_3, v_3, w_3$ are omitted for simplicity.
This is a Hamiltonian system with the Hamiltonian
\begin{equation}
\widetilde{H} = \frac{1}{2}(-u^2w^4+(1\pm 2\alpha )uw^3\pm uvw^2 \pm 2u) 
   \mp \frac{1}{4}(1\pm 2\alpha ) (vw\pm \frac{1}{2}(1\pm 2\alpha )w^2),
\end{equation}
as well as the original ($\text{P}_\text{II}$).
The relation between the original chart $(x, y, z)$ and $(u_3, v_3, w_3)$ is
\begin{equation}
\left\{ \begin{array}{l}
\displaystyle x = u_3w_3^2 \mp w_3^{-2} \mp \frac{1}{2}v_3 - (\frac{1}{2} \pm \alpha ) w_3, \\
\displaystyle y = w_3^{-1},  \\
\displaystyle z = v_3,
\end{array} \right.
\label{4-25}
\end{equation}
or
\begin{equation}
\left\{ \begin{array}{l}
\displaystyle u_3 = xy^{2} \pm y^4 \pm \frac{1}{2}zy^2 + (\frac{1}{2}\pm \alpha )y, \\
\displaystyle w_3 = y^{-1},  \\
\displaystyle v_3 = z.
\end{array} \right.
\label{4-26}
\end{equation}
It is remarkable that the system (\ref{4-23}) is polynomial,
and the independent variable $z$ is not changed despite the fact that 
$z$ is changed in each step of transformations.
Now we distinguish the choice of the sign of $\hat{X}_2 = X_2 \pm 1$.
For the upper sign $\hat{X}_2 = X_2 + 1$ and the lower one $\hat{X}_2 = X_2 - 1$,
we use the notation $(u_3^+, v_3^+, w_3^+)$ and $(u_3^-, v_3^-, w_3^-)$, respectively.
Eq.(\ref{4-26}) should be
\begin{eqnarray}
\left\{ \begin{array}{l}
\displaystyle u^+_3 = xy^{2} + y^4 + \frac{1}{2}zy^2 + (\frac{1}{2} + \alpha )y, \\
\displaystyle w^+_3 = y^{-1},  \\
\displaystyle v^+_3 = z.
\end{array} \right.
\left\{ \begin{array}{l}
\displaystyle u^-_3 = xy^{2} - y^4 - \frac{1}{2}zy^2 + (\frac{1}{2} - \alpha )y, \\
\displaystyle w^-_3 = y^{-1},  \\
\displaystyle v^-_3 = z.
\end{array} \right.
\label{4-27}
\end{eqnarray}
Thus, $\C^2_{(x,y)}, \C^2_{(u^+_3, w^+_3)}$ and $\C^2_{(u^-_3, w^-_3)}$ are glued by (\ref{4-27}) to define 
an algebraic surface $E(z)$, which gives the space of initial conditions for ($\text{P}_\text{II}$).
The transformation (\ref{4-27}) yields
\begin{equation}
dy \wedge dx = du^{\pm}_3 \wedge dw^{\pm}_3,
\end{equation}
where $z$ is regarded as a parameter, and
\begin{equation}
dy \wedge dx - dH \wedge dz = du^{\pm}_3 \wedge dw^{\pm}_3 - d\widetilde{H} \wedge dz,
\end{equation}
where $z$ is regarded as a coordinate.
These results are summarized as follows.
\\[0.2cm]
\textbf{Theorem.\thedef.} 

\noindent \textbf{(i)} 
The space of initial conditions $E(z)$ for ($\text{P}_\text{II}$) is given by
$\C^2_{(x, x)} \cup \C^2_{(u^+_3, w^+_3)}\cup \C^2_{(u^-_3, w^-_3)}$ glued by (\ref{4-27}).

\noindent \textbf{(ii)} The transformation (\ref{4-27}) is symplectic,
and ($\text{P}_\text{II}$) written in $(u_3^+, w_3^+)$ and $(u_3^-, w_3^-)$ are also polynomial Hamiltonian systems.

\noindent \textbf{(iii)} Consider a Hamiltonian system (\ref{3-43}) defined on the $(x, y, z)$-coordinates,
where $f$ and $g$ are polynomials in $x$ and $y$.
If it is also expressed as a polynomial Hamiltonian system in the $(u^{\pm}_3, w^{\pm}_3)$-coordinates, 
then (\ref{3-43}) is ($\text{P}_\text{II}$). 
\\

Part (iii) is proved in the same way as Thm.4.1 and omitted.
The same result is also obtained by Takano et al. \cite{Tak1, Tak2, Tak3}
using a slightly different coordinates.

%%%%%%%%%%%%%%%%%%%%%%%%%%%%%%%%%%%%%%%%%%%%%%%%%%%%%%%%%%%%%%%%%%%%%%%%%%%%%%%%%%%%%%%%%%%%%%%%%%

\subsection{The fourth Painlev\'{e} equation}

($\text{P}_\text{IV}$) have been written as three dimensional vector fields (\ref{3-33}), (\ref{3-34}) and (\ref{3-35}).
They have three fixed points $(Y_1, Z_1, \varepsilon _1) = (0,0,0), (1,0,0)$ and $(X_2, Z_2, \varepsilon _2) = (0,0,0)$,
which correspond to movable singularities.
Let us construct the space of initial conditions for ($\text{P}_\text{IV}$) 
by weighted blow-ups at the three fixed points.

\textbf{(i)} $(Y_1, Z_1, \varepsilon _1) = (0,0,0)$.

For Eq.(\ref{3-33}), we change the coordinates by the linear transformation
\begin{eqnarray}
\left\{ \begin{array}{l}
Y_1 = U_1 + 2\kappa_0 W_1, \\
Z_1 = V_1,  \\
\varepsilon _1 = W_1,
\end{array} \right.
\end{eqnarray}
which results in
\begin{eqnarray}
\left\{ \begin{array}{l}
\displaystyle \dot{U}_1 = 3U_1 + f_1(U_1, V_1, W_1),  \\
\displaystyle \dot{V}_1 = V_1 + W_1 + f_2(U_1, V_1, W_1),  \\
\displaystyle \dot{W}_1 = 2W_1,
\end{array} \right.
\label{4-32}
\end{eqnarray}
where $f_1$ and $f_2$ denote nonlinear terms.
We introduce the weighted blow-up with weights $3,1,2$ by
\begin{eqnarray}
\left\{ \begin{array}{llll}
U_1 & = u_1^3     & = v_2^3 u_2 & = w_3^3 u_3, \\
V_1 & = u_1v_1  & = v_2     & = w_3 v_3, \\
W_1 & = u_1^2w_1  & = v_2^2 w_2 & = w_3^2.  
\end{array} \right.
\label{4-33}
\end{eqnarray}
The exceptional divisor $\{ u_1 = 0\} \cup \{ v_2 = 0\} \cup \{ w_3= 0\}$
is a 2-dim weighted projective space $\C P^2(3,1,2)$.
In the $(u_3, v_3, w_3)$-coordinates, Eq.(\ref{4-32}) is written as
\begin{equation}
\left\{ \begin{array}{l}
\displaystyle \frac{du}{dv} = -(4\theta _\infty- 8 \kappa_0 ) uw + 2uv + 3u^2w^2 + 4\kappa_0(\kappa_0-\theta _\infty)
 = -\frac{\partial \widetilde{H}_1}{\partial w}, \\[0.2cm]
\displaystyle \frac{dw}{dv} = 1-2uw^3-2vw + (2\theta _\infty-4\kappa_0) w^2
 = \frac{\partial \widetilde{H}_1}{\partial u},  \\
\end{array} \right.
\label{4-34}
\end{equation}
where the subscripts for $u_3, v_3, w_3$ are omitted for simplicity.
This is a Hamiltonian system with the Hamiltonian
\begin{equation}
\widetilde{H}_1 = u-u^2w^3-2uvw+(2\theta _\infty-4\kappa_0) uw^2-4\kappa_0(\kappa_0-\theta _\infty)w
\end{equation}
as well as the original ($\text{P}_\text{IV}$).
The relation between the original chart $(x, y, z)$ and $(u_3, v_3, w_3)$ is
\begin{equation}
\left\{ \begin{array}{l}
\displaystyle x =w_3^{-1}, \\
\displaystyle y = u_3w_3^2 + 2\kappa_0 w_3,  \\
\displaystyle z = v_3.
\end{array} \right.
\label{4-36}
\end{equation}
We can verify the equalities
\begin{eqnarray*}
& & dx \wedge dy = du_3 \wedge w_3, \\
& & dx \wedge dy + dH\wedge dz = du_3 \wedge dw_3 + d \widetilde{H}_1 \wedge dz,
\end{eqnarray*}
where $z$ is regarded as a parameter in the former relation, and as a coordinate in the latter one.

\textbf{(ii)} $(Y_1, Z_1, \varepsilon _1) = (1,0,0)$.

In this case, putting $\hat{Y}_1 = Y_1 -1$ for Eq.(\ref{3-33}) yields
\begin{eqnarray}
\left\{ \begin{array}{l}
\displaystyle \dot{\hat{Y}}_1 = 3\hat{Y}_1 + 4Z_1 + (2\kappa_0 - 2\theta _\infty) \varepsilon _1
   + f_1(\hat{Y}_1, Z_1, \varepsilon _1),  \\
\displaystyle \dot{Z}_1 = Z_1 - \varepsilon _1 + f_2(\hat{Y}_1, Z_1, \varepsilon _1),  \\
\displaystyle \dot{\varepsilon }_1 = 2\varepsilon _1,
\end{array} \right.
\end{eqnarray}
where $f_1$ and $f_2$ denote nonlinear terms.
After a certain linear transformation $(\hat{Y}_1,Z_1, \varepsilon _1) \mapsto (U_2, V_2, W_2)$,
which removes the $(1,2)$ and $(1,3)$-components of the linear part as before, 
we introduce the weighted blow-up with weights $3,1,2$ by
\begin{eqnarray}
\left\{ \begin{array}{llll}
U_2 & = u_4^3     & = v_5^3 u_5 & = w_6^3 u_6, \\
V_2 & = u_4v_4  & = v_5     & = w_6 v_6, \\
W_2 & = u_4^2w_4  & = v_5^2 w_5 & = w_6^2.  
\end{array} \right.
\end{eqnarray}
In the $(u_6, v_6, w_6)$-coordinates, we obtain
\begin{equation}
\left\{ \begin{array}{l}
\displaystyle \frac{du}{dv} = 3u^2w^2-2uv-4(2\kappa_0-\theta _\infty -2) uw + 4(1-2\kappa_0+\theta _\infty +\kappa_0^2 - \kappa_0 \theta _\infty)
 = -\frac{\partial \widetilde{H}_2}{\partial w}, \\[0.2cm]
\displaystyle \frac{dw}{dv} = -1-2uw^3 + 2vw + 2(2\kappa_0-\theta _\infty -2)w^2
 = \frac{\partial \widetilde{H}_2}{\partial u},  \\[0.2cm]
\widetilde{H}_2 := -u-u^2w^3 + 2uvw + 2(2\kappa_0-\theta _\infty -2) uw^2 - 4(1-2\kappa_0+\theta _\infty +\kappa_0^2 - \kappa_0 \theta _\infty) w.
\end{array} \right.
\end{equation}
where the subscripts for $u_6, v_6, w_6$ are omitted for simplicity.
The relation with the original chart $(x, y, z)$ is
\begin{equation}
\left\{ \begin{array}{l}
\displaystyle x =w_6^{-1}, \\
\displaystyle y =w_6^{-1}+u_6w_6^2-2v_6-2(\kappa_0-\theta _\infty -1)w_6 ,  \\
\displaystyle z = v_6,
\end{array} \right.
\label{4-40}
\end{equation}
which is symplectic as the case (i).

\textbf{(iii)} $(X_2, Z_2, \varepsilon _2) = (0,0,0)$.

After a certain linear transformation $(X_2, Z_2, \varepsilon _2) \mapsto (U_3, V_3, W_3)$ as before, 
we introduce the weighted blow-up with weights $3,1,2$ by
\begin{eqnarray}
\left\{ \begin{array}{llll}
U_3 & = u_7^3     & = v_8^3 u_8 & = w_9^3 u_9, \\
V_3 & = u_7v_7  & = v_8     & = w_9 v_9, \\
W_3 & = u_7^2w_7  & = v_8^2 w_8 & = w_9^2.  
\end{array} \right.
\end{eqnarray}
In the $(u_9, v_9, w_9)$-coordinates, we obtain
\begin{equation}
\left\{ \begin{array}{l}
\displaystyle \frac{du}{dv} =3u^2w^2 -2uv -(4\kappa_0- 8\theta _\infty) uw + 4\theta _\infty (\theta _\infty - \kappa_0)
 = -\frac{\partial \widetilde{H}_3}{\partial w}, \\[0.2cm]
\displaystyle \frac{dw}{dv} = 1-2uw^3 + 2vw + (2\kappa_0 - 4\theta _\infty) w^2
 = \frac{\partial \widetilde{H}_3}{\partial u},  \\[0.2cm]
\widetilde{H}_3 := u-u^2w^3 + 2uvw + (2\kappa_0 - 4\theta _\infty) uw^2 - 4\theta _\infty (\theta _\infty-\kappa_0) w,
\end{array} \right.
\end{equation}
where the subscripts for $u_9, v_9, w_9$ are omitted.
The relation with the original chart $(x, y, z)$ is
\begin{equation}
\left\{ \begin{array}{l}
\displaystyle x =u_9w_9^2 + 2\theta _\infty w_9, \\
\displaystyle y =w_9^{-1},  \\
\displaystyle z = v_9,
\end{array} \right.
\label{4-44}
\end{equation}
which is symplectic as the case (i).

Thus, $\C^2_{(x,y)}, \C^2_{(u_3, w_3)}, \C^2_{(u_6, w_6)}$ and $\C^2_{(u_9, w_9)}$ are glued by
above symplectic transformations to define 
an algebraic surface $E(z)$, which gives the space of initial conditions for ($\text{P}_\text{IV}$).
These results are summarized as follows.
\\[0.2cm]
\textbf{Theorem.\thedef.} 

\noindent \textbf{(i)} 
The space of initial conditions $E(z)$ for ($\text{P}_\text{IV}$) is given by
$\C^2_{(x,y)} \cup \C^2_{(u_3, w_3)} \cup \C^2_{(u_6, w_6)} \cup \C^2_{(u_9, w_9)}$ glued by 
(\ref{4-36}), (\ref{4-40}) and (\ref{4-44}).

\noindent \textbf{(ii)} These transformations are symplectic.

\noindent \textbf{(iii)} Consider a Hamiltonian system (\ref{3-43}) defined on the $(x, y, z)$-coordinates,
where $f$ and $g$ are polynomials in $x$ and $y$.
If it is also expressed as polynomial Hamiltonian systems in the 
$(u_3, w_3), (u_6, w_6)$ and $(u_9, w_9)$-coordinates, 
then (\ref{3-43}) is ($\text{P}_\text{IV}$). 
\\

Part (iii) is proved in the same way as Thm.4.1 and omitted.

%%%%%%%%%%%%%%%%%%%%%%%%%%%%%%%%%%%%%%%%%%%%%%%%%%%%%%%%%%%%%%%%%%%%%%%%%%%%%%%%%%%%%%%%%%%%%%%%%%
%%%%%%%%%%%%%%%%%%%%%%%%%%%%%%%%%%%%%%%%%%%%%%%%%%%%%%%%%%%%%%%%%%%%%%%%%%%%%%%%%%%%%%%%%%%%%%%%%%

\section{Characteristic index}

As usual, the weight $(p,q,r,s)$ denotes $(3,2,4,5)$, $(2,1,2,3)$ and $(1,1,1,2)$ for ($\text{P}_\text{I}$),
($\text{P}_\text{II}$) and ($\text{P}_\text{IV}$), respectively.
Recall the decomposition
\begin{eqnarray*}
\C P^3(p,q,r,s) = \C^3 / \Z_s \cup \C P^2(p,q,r).
\end{eqnarray*} 
The Painlev\'{e} equation is defined on the covering space of $ \C^3 / \Z_s$, and $\C P^2(p,q,r)$
is attached at infinity.
We have seen that there are fixed points of the vector field on $\C P^2(p,q,r)$ 
which correspond to movable singularities.
In Sec.3 and 4, the eigenvalues $(\lambda _1, \lambda _2, \lambda _3)$ of the Jacobi matrices at the fixed
points play an important role, where $(\lambda _1, \lambda _2, \lambda_3) = (6,4,5), (4,2,3)$ 
and $(3,1,2)$ for ($\text{P}_\text{I}$), ($\text{P}_\text{II}$) and ($\text{P}_\text{IV}$), respectively.
In Sec.3, they are used to apply the Poincar\'{e} linearization theorem.
In Sec.4, they determine the weight of the weighted blow-up.
We call the eigenvalues $(\lambda _1, \lambda _2, \lambda _3)$ the characteristic index 
for ($\text{P}_\text{J}$), see Table.2.
Obviously, they are invariant under the automorphism on $\C P^3(p,q,r,s)$.
 
\begin{table}[h]
\begin{center}
\begin{tabular}{|c||c|c|c|c|}
\hline
 & $\C P^3(p,q,r,s)$ & $\lambda _1, \lambda _2, \lambda _3$ & $(\lambda _1 + \lambda _2)/\lambda _3$ & $\kappa$ \\ \hline \hline
($\text{P}_\text{I}$)  & $\C P^3(3,2,4,5)$ & $6,4,5$  & 2 & 6   \\ \hline
($\text{P}_\text{II}$)  & $\C P^3(2,1,2,3)$ & $4,2,3$ & 2 & 4   \\ \hline
($\text{P}_\text{IV}$)  & $\C P^3(1,1,1,2)$ & $3,1,2$ & 2 & 3  \\ \hline
\end{tabular}
\end{center}
\caption{The characteristic index $(\lambda _1, \lambda _2, \lambda _3)$ and the Kovalevskaya exponent $\kappa$.}
\end{table}

We observe the following properties:
\\[0.2cm]
\textbf{(i)} $r = \lambda _2$ and $s = \lambda _3$. \\
\textbf{(ii)} $\lambda _1 = \kappa = \mathrm{deg} (H_\text{J}) = s+1$, where $\kappa$ is the Kovalevskaya exponent 
and $\mathrm{deg} (H_\text{J})$ is a weighted degree of the Hamiltonian given in Sec.2.1.  \\
\textbf{(iii)} $(\lambda _1 + \lambda _2)/\lambda _3$ is an integer. \\
\textbf{(iv)} $p+q = s$. \\

In the forthcoming paper \cite{Chi4}, the properties (i) and (ii) are proved for a general $m$-dimensional
system satisfying certain conditions on the Newton diagram;
two numbers in the characteristic index coincide with $r$ and $s$ determined by the Newton diagram,
and the others coincide with the Kovalevskaya exponents.

Part (iii) and (iv) are related to the following proposition. 
\\[0.2cm]
\textbf{Proposition.\thedef.} For ($\text{P}_\text{I}$), ($\text{P}_\text{II}$) and ($\text{P}_\text{IV}$),
the symplectic form $dx \wedge dy + dH_J \wedge dz$ is a rational form on $\C P^3(p,q,r,s)$.

This can be proved by a straightforward calculation.
For example, on the $(Y_1, Z_1, \varepsilon _1)$ chart, we have
\begin{eqnarray*}
dx \wedge dy = d (\varepsilon _1^{-p/s}) \wedge d(Y_1\varepsilon _1^{-q/s}).
\end{eqnarray*}
In order for it to be rational in $\varepsilon _1$, $p+q$ should be a multiple of $s$,
and this is true for ($\text{P}_\text{I}$), ($\text{P}_\text{II}$) and ($\text{P}_\text{IV}$).
Similarly, we can verify that $dH \wedge dz$ is rational, which is a consequence of Part (iii) above.
%%%%%%%%%%%%%%%%%%%%%
% \begin{eqnarray*}
% dH _I \wedge dz &=& d (\frac{1}{2}\varepsilon _1^{-2p/s} - 2Y_1^3 \varepsilon _1^{-3q/s} - Z_1Y_1 \varepsilon _1^{-(q+r)/s})
%    \wedge d(Z_1 \varepsilon _1^{-r/s}) \\
% &=& d (\frac{1}{2}\varepsilon _1^{-6/5} - 2Y_1^3 \varepsilon _1^{-6/5} - Z_1Y_1 \varepsilon _1^{-6/5})
%     \wedge d(Z_1 \varepsilon _1^{-4/5})
% \end{eqnarray*} 
% is rational in $\varepsilon _1$.
% The integer $6$ in the exponent of $\varepsilon _1^{-6/5}$ is  the weighted degree of $H_I$ by the definition.
% Thus, we have
% \begin{eqnarray*}
% \frac{\lambda _1 + \lambda _2}{\lambda _3} = \frac{\mathrm{deg} (H_\text{J}) + r}{s}.
% \end{eqnarray*}
% This implies that the property (iii) above is a necessary condition for $dH_I \wedge dz$ to be rational.
%%%%%%%%%%%%%%%%%%%%%

The above properties (i) to (iv) are true for higher order first and second Painlev\'{e} equations.
For example, the characteristic indices and the Kovalevskaya exponents for the fourth order and sixth order
first Painlev\'{e} equations, $(\text{P}_\text{I})_4$ and $(\text{P}_\text{I})_6$, 
from the first Painlev\'{e} hierarchy $(\text{P}_\text{I})_{2n}$,
and the fourth order second Painlev\'{e} equation $(\text{P}_\text{II})_4$ from the second Painlev\'{e} hierarchy 
$(\text{P}_\text{II})_{2n}$ are shown in Table 3.
In Table 3, the equality $p_j + q_j = s$ holds and $\nu := (\lambda _1+ \cdots + \lambda _{2n})/\lambda _{2n+1}$ is an integer. 
Furthermore, $r=\lambda _{2n}, s=\lambda _{2n+1}$ and the Kovalevskaya exponents coincide with $\lambda _1, \cdots ,\lambda _{2n-1}$.
See \cite{Chi4} for the detail.
\begin{table}[h]
\begin{center}
\begin{tabular}{|c||c|c|c|c|}
\hline
$(\text{P}_\text{J})_{2n}$ & $\C P^{2n+1}(p_1, q_1, \cdots , p_n, q_n, r,s)$ & $\lambda _1, \cdots , \lambda _{2n+1}$ & $\nu$ & $\kappa$ \\ \hline \hline
$(\text{P}_\text{I})_4$  & $\C P^5(5,2,3,4,6,7)$ & $8,5,2,6,7$  & 3 & 8,5,2   \\ \hline
$(\text{P}_\text{I})_6$  & $\C P^7(7,2,5,4,3,6,8,9)$ & $10,7,5,4,2,8,9$ & 4 & 10,7,5,4,2   \\ \hline
$(\text{P}_\text{II})_4$  & $\C P^5(2,3,4,1,4,5)$ & $6,3,2,4,5$ & 3 & 6,3,2  \\ \hline
\end{tabular}
\end{center}
\caption{The characteristic index $(\lambda _1, \cdots , \lambda _{2n+1})$, the Kovalevskaya exponent $\kappa$
and $\nu := (\lambda _1+ \cdots + \lambda _{2n})/\lambda _{2n+1}$.}
\end{table}

%%%%%%%%%%%%%%%%%%%%%%%%%%%%%%%%%%%%%%%%%%%%%%%%%%%%%%%%%%%%%%%%%%%%%%%%%%%%%%%%%%%%%%%%%%%%%%%%%%
%%%%%%%%%%%%%%%%%%%%%%%%%%%%%%%%%%%%%%%%%%%%%%%%%%%%%%%%%%%%%%%%%%%%%%%%%%%%%%%%%%%%%%%%%%%%%%%%%%

\section{The Boutroux coordinates}

In his celebrated paper, Boutroux \cite{Bou} introduced the coordinate transformations
$y = u t^{2/5},\, z = t^{4/5}$ for ($\text{P}_\text{I}$), and $y = u t^{1/3},\, z = t^{2/3}$ for ($\text{P}_\text{II}$)
to investigate the asymptotic behavior of solutions around the essential singularity $z = \infty$.
His coordinate transformations are essentially the same as the third local chart of $\C P^3(p,q,r,s)$,
see Eq.(\ref{3-1}) and (\ref{3-17}), and put $\varepsilon _3 = 1/t$.
Hence, we call the third local chart $(X_3, Y_3, \varepsilon _3)$ of $\C P^3(p,q,r,s)$ the 
Boutroux coordinates even for ($\text{P}_\text{IV}$).
On the Boutroux coordinates, ($\text{P}_\text{I}$), ($\text{P}_\text{II}$) and ($\text{P}_\text{IV}$) 
expressed as the autonomous vector fields are given by Eq.(\ref{3-7}), (\ref{3-22}) and (\ref{3-35}), respectively.
Recall that the set $\C P^2(p,q,r) = \{ \varepsilon _1 = 0\} \cup \{ \varepsilon _2 = 0\}\cup \{ \varepsilon _3 = 0\}$
is attached at infinity of the original chart $(x,y,z)$.
In particular, the set $\{ \varepsilon _3 \, ( = z^{-s/r}) = 0 \}$ describes the behavior around the irregular singular point $z = 0$.
Putting $\varepsilon _3 = 0$ in Eq.(\ref{3-7}), (\ref{3-22}) and (\ref{3-35}), we obtain
\begin{equation}
\left\{ \begin{array}{l}
\displaystyle \dot{X}_3 =24Y_3^2 + 4 ,  \\
\displaystyle \dot{Y}_3 = 4X_3 ,  \\
\end{array} \right.
\label{Bou1}
\end{equation}
\begin{equation}
\left\{ \begin{array}{l}
\displaystyle \dot{X}_3 =4Y_3^3 + 2Y_3,  \\
\displaystyle \dot{Y}_3 = 2X_3,  \\
\end{array} \right.
\label{Bou2}
\end{equation}
and
\begin{equation}
\left\{ \begin{array}{l}
\displaystyle \dot{X}_3 = -X_3^2+2X_3Y_3+2X_3,  \\
\displaystyle \dot{Y}_3 = -Y_3^2+2X_3Y_3-2Y_3,  \\
\end{array} \right.
\label{Bou4}
\end{equation}
respectively.
It is remarkable that they are autonomous Hamiltonian systems with the Hamiltonian functions $\mathcal{H}_J$ shown in Table 4.

\begin{table}[h]
\begin{center}
\begin{tabular}{|c||c|c|c|}
\hline
 & space & $\mathcal{H}_J$ & symmetry \\ \hline \hline
($\text{P}_\text{I}$)   & $\C P^3(3,2,4,5)$ & $2X_3^2 - 8Y_3^3 - 4Y_3$  & $\Z_4$    \\ \hline
($\text{P}_\text{II}$)  & $\C P^3(2,1,2,3)$ & $X_3^2 - Y_3^4 - Y_3^2$  & $\Z_2\times \Z_2$  \\ \hline
($\text{P}_\text{IV}$)  & $\C P^3(1,1,1,2)$ & $X_3^2Y_3 - X_3Y_3^2 - 2X_3Y_3$  & $\mathfrak{S}_3$  \\ \hline
\end{tabular}
\end{center}
\caption{Hamiltonian functions defined on the set $\{ \varepsilon _3 = 0\}$. See Section 7 for the symmetry.}
\end{table}

Therefore, each leaf of the foliation on the set $\{ \varepsilon _3 = 0\} \subset \C P^2(p,q,r)$
is determined by the level set $\{ \mathcal{H}_J = c\},\, c\in \C $ of the Hamiltonian, 
which is an elliptic curve for a generic value of $c\in \C$.
In particular, $\mathcal{H}_I = c$ is the Weierstrass form and $\mathcal{H}_{II} = c$ is the Jacobi form.

In what follows, we will see that the space of initial conditions written in the Boutroux coordinates has
been already constructed by the weighted blow-ups introduced in Sec.4.

%%%%%%%%%%%%%%%%%%%%%%%%%%%%%%%%%%%%%%%%%%%%%%%%%%%%%%%%%%%%%%%%%%%%%%%%%%%%%%%%%%%%%%%%%%%%%%%%%%

\subsection{The first Painlev\'{e} equation}

We have proved that $(x, y)$-space and $(u_3, w_3)$-space are glued to give the 
space of initial conditions.
In this subsection, it is shown that $(X_3, Y_3)$-space and $(u_2, v_2)$-space give the 
space of initial conditions for ($\text{P}_\text{I}$) written in the Boutroux coordinates.
On the $(X_3, Y_3, \varepsilon _3)$-coordinates, ($\text{P}_\text{I}$) is given as
\begin{equation}
\frac{dX_3}{d\varepsilon _3} = \frac{1}{5\varepsilon _3^2} \left( -24Y_3^2 - 4 + 3X_3\varepsilon _3 \right), \quad
\frac{dY_3}{d\varepsilon _3} = \frac{1}{5\varepsilon _3^2} \left( -4X_3 + 2Y_3\varepsilon _3\right) .
\label{5-1}
\end{equation}
It has an irregular singular point $\varepsilon _3 = 0$.
Since $z = \varepsilon _3^{-4/5}$, this equation also has the Painlev\'{e} property with 
a possible branch point $\varepsilon _3 = 0$.

Recall that the $(u_2, v_2, w_2)$-coordinates are defined by (\ref{4-3}).
In these coordinates, ($\text{P}_\text{I}$) is written as
\begin{equation}
\left\{ \begin{array}{l}
\displaystyle \frac{du_2}{dw_2} = \frac{1}{5w_2^2} \left( -\frac{1}{2}v_2 - 6u_2^2v_2^5 + 6u_2w_2
             \mp \frac{3}{2}v_2^2w_2 + 5u_2 v_2^4w_2 \pm 4u_2v_2^3 - v_2^3w_2^2 \right),   \\
\displaystyle \frac{dv_2}{dw_2} = \frac{1}{5w_2^2}\left( \mp 4 \mp v_2^4 + 2u_2v_2^6 - v_2 w_2 - v_2^5 w_2 \right) .
\end{array} \right.
\label{5-2}
\end{equation}
This is a polynomial ODE with an irregular singular point $w_2 = 0$.
The relation between two charts are
\begin{equation}
\left\{ \begin{array}{l}
\displaystyle X_3 = v_2^{-3} (v_2^6 u_2 \mp \frac{1}{2}v_2^4 - \frac{1}{2}v_2^5 w_2 \mp 2),  \\
\displaystyle Y_3 = v_2^{-2},  \\
\displaystyle \varepsilon _3 = w_2,
\end{array} \right.
\label{5-3}
\end{equation}
or 
\begin{equation}
\left\{ \begin{array}{l}
\displaystyle u_2 = X_3Y_3^{3/2} + \frac{1}{2}Y_3^{1/2} \varepsilon _3 \pm \frac{1}{2}Y_3 \pm 2Y_3^3,  \\
\displaystyle v_2 = Y_3^{-1/2},  \\
\displaystyle w_2 = \varepsilon _3.
\end{array} \right.
\label{5-4}
\end{equation}
We should divide the $(u_2, v_2, w_2)$-space by the $\Z_2$ action as before.
The $\Z_2$ action induced from the orbifold structure is given by
\begin{equation}
\left(
\begin{array}{@{\,}c@{\,}}
u_2 \\
v_2 \\
w_2
\end{array}
\right) \mapsto \left(
\begin{array}{@{\,}c@{\,}}
-u_2 + v_2^{-2} + 4v_2^{-6} \\
-v_2 \\
w_2
\end{array}
\right) .
\end{equation}
We define invariants of this action to be
\begin{equation}
\left\{ \begin{array}{l}
\displaystyle U = u_2 (u_2 v_2^6 - v_2^4 - 4) + \frac{1}{4}v_2^2 \\
\displaystyle V = v_2^7( u_2 - \frac{1}{2} v_2^{-2} - 2v_2^{-6}) \\
\displaystyle W = v_2^2.
\end{array} \right.
\end{equation}
This defines a nonsingular algebraic surface $M = \C ^2_{(u_2, v_2)} / \Z_2$
\begin{equation}
M : V^2 = UW^4 + 2W^3 + 4W.
\end{equation}
Note that it is independent of a parameter.
Eqs.(\ref{5-2}) to (\ref{5-4}) are rewritten as
\begin{equation}
\left\{ \begin{array}{l}
\displaystyle X_3 = W^{-2}V - \frac{1}{2}W\varepsilon _3  \\
\displaystyle Y_3 = W^{-1},
\end{array} \right. \quad \left\{ \begin{array}{l}
\displaystyle V = X_3Y_3^{-2}+\frac{1}{2}Y_3^{-3}\varepsilon _3  \\
\displaystyle W = Y_3^{-1},
\end{array} \right.
\end{equation}
and 
\begin{equation}
\left\{ \begin{array}{l}
\displaystyle \frac{dV}{d\varepsilon _3} = \frac{1}{5\varepsilon _3^2} \left( 
8 \frac{V^2}{W} - \varepsilon _3 V - 2\varepsilon _3 VW^2 - 4W^2 - \varepsilon _2^2 W^5 - 24 \right), \\
\displaystyle \frac{dW}{d\varepsilon _3} = \frac{1}{5\varepsilon _3^2} \left( 
-2W\varepsilon _3 + 4V - 2W^3 \varepsilon _3 \right),
\end{array} \right.
\label{5-9}
\end{equation}
respectively.
Hence, $\C^2_{(X_3, Y_3)} \cup M$ gives the space of initial conditions for the Boutroux coordinates.
Note that Eqs.(\ref{5-1}) and (\ref{5-9}) are not Hamiltonian systems,
though they are reduced to Hamiltonian systems as $\varepsilon _3 \to 0$.

%%%%%%%%%%%%%%%%%%%%%%%%%%%%%%%%%%%%%%%%%%%%%%%%%%%%%%%%%%%%%%%%%%%%%%%%%%%%%%%%%%%%%%%%%%%%%%%%%%

\subsection{The second Painlev\'{e} equation}

On the $(X_3, Y_3, \varepsilon _3)$ chart, ($\text{P}_\text{II}$) is written as
\begin{equation}
\frac{dX_3}{d\varepsilon _3} = \frac{4Y_3^3 + 2Y_3 + 2\alpha \varepsilon _3 - 2X_3\varepsilon _3}{-3\varepsilon _3^2}, \quad
\frac{dY_3}{d\varepsilon _3} = \frac{2X_3 - Y_3\varepsilon _3}{-3\varepsilon _3^2}.
\label{5-10}
\end{equation}
The space of initial conditions for this system is also obtained by the weighted blow up as follows:
Recall that $(u_2^{\pm},v_2^{\pm},w_2^{\pm})$ is defined by (\ref{4-22}).
The relation between $(X_3,Y_3, \varepsilon _3)$ and $(u_2^{\pm},v_2^{\pm},w_2^{\pm})$ is given by
\begin{equation}
\left\{ \begin{array}{l}
\displaystyle u_2^{\pm} = X_3Y_3^2 \pm Y_3^4 \pm \frac{1}{2} Y_3^2 + (\frac{1}{2} \pm \alpha ) \varepsilon _3 Y_3  \\
v_2^{\pm} = Y_3^{-1} \\
w_2^{\pm} = \varepsilon _3.
\end{array} \right.
\label{5-11}
\end{equation}
It is remarkable that the independent variable $\varepsilon _3$ is not changed.
The equation written in $(u_2^{\pm},v_2^{\pm},\varepsilon _3)$ is
\begin{equation}
\left\{ \begin{array}{l}
\displaystyle \frac{du}{d\varepsilon }= \frac{-1}{3\varepsilon ^2}
\left(  4u^2v^3 \mp 2uv \pm (\frac{1}{2} \pm \alpha ) \varepsilon -6 (\frac{1}{2} \pm \alpha ) \varepsilon uv^2
-4 \varepsilon u + \frac{\varepsilon ^2 v}{2} \pm 2 \alpha  \varepsilon ^2v + 2\alpha ^2 \varepsilon ^2v \right) \\[0.4cm]
\displaystyle \frac{dv}{d\varepsilon }= \frac{-1}{3\varepsilon ^2}
\left( \pm 2 \pm v^2- 2uv^4 + \varepsilon v + \varepsilon v^3  \pm 2\alpha \varepsilon v^3  \right),  \\
\end{array} \right.
\label{5-12}
\end{equation}
where the subscript and the superscript for $u_2^{\pm},v_2^{\pm}, \varepsilon _3$ are omitted.
Since this equation is polynomial in $u$ and $v$, the space of initial conditions is 
obtained by glueing $\C^2_{(X_3, Y_3)}$, $\C^2_{(u_2^+, v_2^+)}$ and $\C^2_{(u_2^-, v_2^-)}$ by the relation (\ref{5-11}).
%%%%%%%%%%%%%%%%%%%%%%%%%%%%%%%%%%%%%%%%%%%%%%%%%%%%%%%%%%%%%%%%%%%%%%%%%%%%%%%%%%%%%%%%%%%%%%%%%%

\subsection{The fourth Painlev\'{e} equation}

We call the third local chart $(X_3, Y_3, \varepsilon _3)$ for ($\text{P}_\text{IV}$)
the Boutroux coordinates as in the first and second Painlev\'{e} equations.
In this chart, ($\text{P}_\text{IV}$) is written as
\begin{equation}
\left\{ \begin{array}{ll}
\displaystyle \frac{dX_3}{d\varepsilon _3} = \frac{-X_3^2+2X_3Y_3+2X_3-2\theta _\infty \varepsilon _3-X_3\varepsilon _3}{-2\varepsilon _3^2} \\[0.4cm]
\displaystyle \frac{dY_3}{d\varepsilon _3} = \frac{-Y_3^2+2X_3Y_3-2Y_3-2\kappa_0 \varepsilon _3-Y_3\varepsilon _3}{-2\varepsilon _3^2}. \\
\end{array} \right.
\end{equation}
In Sec.4.3, $(u_2, v_2, w_2), (u_5, v_5, w_5)$ and $(u_8, v_8, w_8)$ coordinates are defined through the weighted blow-ups.
The relations between them and $(X_3,Y_3, \varepsilon _3)$ are given by
\begin{equation}
\left\{ \begin{array}{ll}
X_3 = v_2^{-1} & = v_5^{-1} \\
Y_3 = u_2v_2^2 + 2\kappa_0 w_2 v_2 & = u_5v_5^2 + 2(1 + \theta _\infty - \kappa_0) w_5v_5 + v_5^{-1} -2  \\
\varepsilon _3 = w_2 & = w_5 ,
\end{array} \right.
\end{equation}
and
\begin{equation}
\left\{ \begin{array}{ll}
X_3 = u_8v_8^2 + 2\theta _\infty w_8 v_8  \\
Y_3 = v_8^{-1}  \\
\varepsilon _3 = w_8 
\end{array} \right.
\end{equation}
In particular, the independent variable $\varepsilon _3$ is not changed.
The equations written in $(u_2,v_2, w_2), (u_5, v_5, w_5)$ and $(u_8,v_8,w_8)$ are polynomial in $(u_2,v_2), (u_5, v_5)$ and $(u_8,v_8)$,
respectively.
Thus, the space of initial conditions is obtained by glueing 
$\C^2_{(X_3, Y_3)}, \C^2_{(u_2, v_2)}, \C^2_{(u_5, v_5)}$ and $\C^2_{(u_8, v_8)}$ by the above relation.

%%%%%%%%%%%%%%%%%%%%%%%%%%%%%%%%%%%%%%%%%%%%%%%%%%%%%%%%%%%%%%%%%%%%%%%%%%%%%%%%%%%%%%%%%%%%%%%%%%
%%%%%%%%%%%%%%%%%%%%%%%%%%%%%%%%%%%%%%%%%%%%%%%%%%%%%%%%%%%%%%%%%%%%%%%%%%%%%%%%%%%%%%%%%%%%%%%%%%

\section{The extended affine Weyl group}

It is known that there exists a group of rational transformations acting on $\C^3$ which 
changes the Painlev\'{e} equation to another Painlev\'{e} equation of the same type with
different parameters.
The transformation group is isomorphic to the extended affine Weyl group.
See \cite{Tsu} for the complete list of the actions of the groups for the second to the sixth Painlev\'{e} equations
written in Hamiltonian forms.
In this section, we study the actions of the extended affine Weyl groups for ($\text{P}_\text{II}$) and ($\text{P}_\text{IV}$). 

For a classical root system $R$, the affine Weyl group and the extended affine Weyl group are denoted by
$W(R^{(1)})$ and $\widetilde{W}(R^{(1)})$, respectively.
Let $G = \mathrm{Aut} (R^{(1)})$ be the Dynkin automorphism group of the extended Dynkin diagram.
We have $\widetilde{W}(R^{(1)}) \cong G \ltimes W(R^{(1)})$.

For the second Painlev\'{e} equation,
\begin{eqnarray*}
& & \widetilde{W}(A_1^{(1)}) \cong G \ltimes W(A_1^{(1)}) = \langle  s_1, \pi  \rangle, \\
& & G = \mathrm{Aut} (A_1^{(1)}) = \langle  \pi  \rangle \cong \Z_2.
\end{eqnarray*}
The action of the group is given in Table 5.
\begin{table}[h]
\begin{center}
\begin{tabular}{|c||c|c|c|}
\hline
      & $\alpha$     & $x$ & $y$ \\ \hline \hline
$s_1$ & $-\alpha +1$ & $\displaystyle x+\frac{(2\alpha -1)y}{y^2-x+z/2} + \frac{(\alpha -1/2)^2}{(y^2-x+z/2)^2}$  
               & $\displaystyle y+ \frac{\alpha -1/2}{y^2-x+z/2}$  \\ \hline
$\pi$ & $-\alpha $   & $-x$  & $-y$   \\ \hline
\end{tabular}
\end{center}
\caption{The action of the extended affine Weyl group for ($\text{P}_\text{II}$) .}
\end{table}

For the fourth Painlev\'{e} equation,
\begin{eqnarray*}
& & \widetilde{W}(A_2^{(1)}) \cong G \ltimes W(A_2^{(1)}) = \langle s_0, s_1, s_2, \sigma_1, \sigma_2 \rangle, \\
& & G = \mathrm{Aut} (A_2^{(1)}) = \langle  \sigma_1, \sigma_2 \rangle \cong \mathfrak{S}_3.
\end{eqnarray*}
The action of the group is given in Table 6.
\begin{table}[h]
\begin{center}
\begin{tabular}{|c||c|c|c|c|c|}
\hline
      & $\kappa_0$ & $\theta _\infty$ & $x$ & $y$ & $z$ \\ \hline \hline
$s_0$ & $1+\theta _\infty$  & $\kappa_0-1$ & $\displaystyle x+\frac{2(1-\kappa_0+\theta _\infty)}{x-y-2z}$
           & $\displaystyle y+\frac{2(1-\kappa_0+\theta _\infty)}{x-y-2z}$ & $z$  \\ \hline
$s_1$ & $-\kappa_0$ & $\theta _\infty-\kappa_0$ & $x-2\kappa_0/y$ & $y$ & $z$   \\ \hline
$s_2$ & $\kappa_0 - \theta _\infty$ & $-\theta _\infty$ & $x$ & $y-2\theta _\infty/x$ & $z$ \\ \hline
$\pi$ & $-\theta _\infty$ & $\kappa_0-\theta _\infty -1 $ & $-x+y+2z$ & $-x$ & $z$   \\ \hline
$\sigma_1$ & $-\theta _\infty$ & $-\kappa_0$ & $-iy$ & $-ix$ & $iz$   \\ \hline
$\sigma_2$ & $\kappa_0$ & $\kappa_0 - \theta _\infty - 1$ & $ix$ & $i(x-y-2z)$ & $iz$   \\ \hline
\end{tabular}
\end{center}
\caption{The action of the extended affine Weyl group for ($\text{P}_\text{IV}$) .}
\end{table}

In the next theorem, $R^{(1)}$ and $(p,q,r,s)$ denote $A^{(1)}_1$ and $(2,1,2,3)$ for ($\text{P}_\text{II}$),
and  $A^{(1)}_2$ and $(1,1,1,2)$ for ($\text{P}_\text{IV}$), respectively.
\\[0.2cm]
\textbf{Theorem.\thedef.}

(i) The transformation group $\widetilde{W}(R^{(1)})$ given above is extended to a rational transformation group on $\C P^3(p,q,r,s)$.

(ii) For each element $s$ in $W(R^{(1)})$, the action of $s$ on the infinity set $\C P^2(p,q,r)$ is trivial: $s|_{\C P^2(p,q,r)} = \mathrm{id}$.
Hence, the transformation group $\widetilde{W}(R^{(1)}) \cong \mathrm{Aut} (R^{(1)}) \ltimes W(R^{(1)})$ is reduced to 
$\mathrm{Aut} (R^{(1)})$ on $\C P^2(p,q,r)$.

(iii) The foliation on $\C P^2(p,q,r)$ is $\mathrm{Aut} (R^{(1)})$-invariant.
\\[0.2cm]
\textbf{Proof.} We give a proof for ($\text{P}_\text{II}$).
A proof for ($\text{P}_\text{IV}$) is done in the same way. 

(i) The original chart $(x,y,z)$ for ($\text{P}_\text{II}$) has to be divided by the $\Z_3$-action 
$(x,y,z)\mapsto (\omega ^2 x, \omega y, \omega ^2z),\, \omega :=e^{2\pi i/3}$ because of the orbifold structure, see Eq.(\ref{2-5}).
It is easy to verify that the actions of the generators $s_1$ and $\pi$ of $\widetilde{W}(A_1^{(1)})$ shown in Table 5 commute with 
the $\Z_3$-action, so that they induce actions on the quotient space $\C^3 / \Z_3$.
Further, these actions are extended to the whole space $\C P^3(2,1,2,3)$.
For example, the action of $s_1$ on the third local chart is expressed by
\begin{equation}
(X_3, Y_3,\varepsilon _3) \mapsto (X_3 + \frac{(2\alpha -1)Y_3\varepsilon _3}{Y_3^2-X_3+1/2} +\frac{(\alpha -1/2)^2 \varepsilon _3^2}{(Y_3^2-X_3+1/2)^2}, 
Y_3 + \frac{(\alpha -1/2)\varepsilon _3}{Y_3^2-X_3+1/2}, \varepsilon _3 ).
\label{7-1}
\end{equation}
This is rational and commutes with the action (\ref{3-19}) defining the orbifold structure.
It is straightforward to calculate the action on the other charts $(Y_1, Z_1, \varepsilon _1)$ and $(X_2, Z_2, \varepsilon _2)$,
which proves that the action $s_1$ is well defined on $\C P^3(2,1,2,3)$.
The same is also true for the action of $\pi$.

(ii) The set $\C P^2 (2,1,2)$ is given by $\{ \varepsilon _1 = 0\} \cup \{ \varepsilon _2 = 0\} \cup\{ \varepsilon _3 = 0\}$.
The action (\ref{7-1}) is reduced to the trivial action as $\varepsilon _3 \to 0$.
Similarly, the action of $s_1$ on the other charts $(Y_1, Z_1, \varepsilon _1)$ and $(X_2, Z_2, \varepsilon _2)$ becomes trivial
as $\varepsilon _1 \to 0$ and $\varepsilon _2 \to 0$.
On the other hand, the action of $\pi$ on $\C P^2 (2,1,2)$ is not trivial.
For example, the action of $\pi$ on $(X_3, Y_3,\varepsilon _3)$ is $(X_3, Y_3,\varepsilon _3)\mapsto (-X_3, -Y_3,\varepsilon _3)$.

(iii) The action of $\widetilde{W}(A_1^{(1)})$ transforms ($\text{P}_\text{II}$) into ($\text{P}_\text{II}$) with a different parameter.
However, the foliation on $\C P^2 (2,1,2)$ is independent of the parameter $\alpha $, see Eq.(\ref{Bou2}).
Thus, the foliation on  $\C P^2 (2,1,2)$ is $\mathrm{Aut} (A_1^{(1)})$-invariant. $\Box$
\\

In Table 4 in Sec.6, the symmetry groups of the foliations on $\C P^2 (p,q,r)$ generated by ($\text{P}_\text{I}$), ($\text{P}_\text{II}$) 
and ($\text{P}_\text{IV}$) are shown.
The foliation generated by $\mathcal{H}_I$ is invariant under the $\Z_4$-action (\ref{3-4}) which arises from the orbifold structure.
The foliation generated by $\mathcal{H}_{II}$ is invariant under the $\Z_2$-action (\ref{3-19}), and $\mathrm{Aut} (A_1^{(1)}) \cong \Z_2$
given by $(X_3, Y_3)\mapsto (-X_3, -Y_3)$.
Similarly, the foliation generated by $\mathcal{H}_{IV}$ is invariant under the action of $\mathrm{Aut} (A_2^{(1)}) \cong \mathfrak{S}_3$,
while there are no symmetry induced from the orbifold structure.

The foliation $\{ \mathcal{H}_{IV} = c\}$ for real $c \in \R$ is represented in Fig.1.
Note that $\mathrm{Aut} (A_2^{(1)})$ is isomorphic to the dihedral group $D_3$ of a triangle.
In Fig.1, the zero level set $\mathcal{H}_{IV} = 0$ consists of three lines, which creates a triangle.
$\mathrm{Aut} (A_2^{(1)})$ acts on the triangle as the  dihedral group $D_3$.

\begin{figure}
\begin{center}
\includegraphics[]{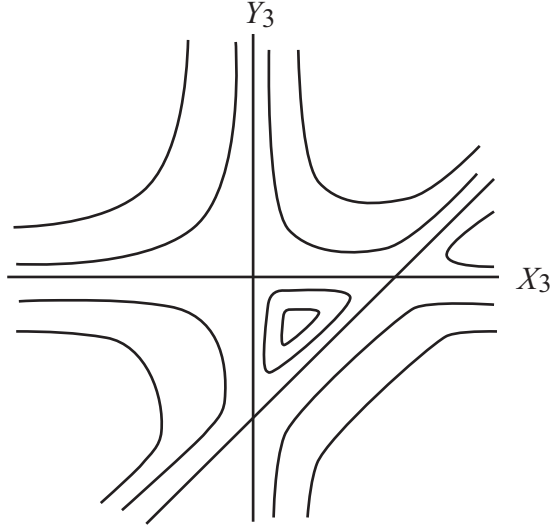}
\caption[]{The foliation $\{ \mathcal{H}_{IV} = c\}$ for $c \in \R$.}
\label{fig1}
\end{center}
\end{figure}

%%%%%%%%%%%%%%%%%%%%%%%%%%%%%%%%%%%%%%%%%%%%%%%%%%%%%%%%%%%%%%%%%%%%%%%%%%%%%%%%%%%%%%%%%%%%%%%%%%
%%%%%%%%%%%%%%%%%%%%%%%%%%%%%%%%%%%%%%%%%%%%%%%%%%%%%%%%%%%%%%%%%%%%%%%%%%%%%%%%%%%%%%%%%%%%%%%%%%

\section{Cellular decomposition and Dynkin diagrams}

In this section, the weighted blow-up of $\C P^3(3,2,4,5)$ with weights $6,4,5$ 
defined in Sec.4.1 is called the total space for ($\text{P}_\text{I}$) and denoted by $\mathcal{M}_I$.
We calculate the cellular decomposition of it.
It will be shown that $\mathcal{M}_I$ is decomposed into the fiber space $\mathcal{P} = \{ (x,z) \, | \, x\in E(z), z\in \C \}$ for 
($\text{P}_\text{I}$), an elliptic fibration over the moduli space of complex tori
defined by the Weierstrass equation and a projective line.
We also shows that the extended Dynkin diagram of type 
$\tilde{E}_8$ is hidden in the space $\mathcal{M}_I$.

%%%%%%%%%%%%%%%%%%%%%%%%%%%%%%%%%%%%%%%%%%%%%%%%%%%%%%%%%%%%%%%%%%%%%%%%%%%%%%%%%%%%%%%%%%%%%%%%%%

\subsection{The elliptic fibration}

Let us calculate a cellular decomposition of $\mathcal{M}_I$.
$\C P^3(3,2,4,5)$ is decomposed as (\ref{2-7}).
Furthermore, $\C P^2(3,2,4)$ is decomposed as $\C P^2(3,2,4) = \C^2/\Z_4 \cup \C P^1(3,2)$.
Since $\C P^1(3,2)$ is isomorphic to the Riemann sphere, we obtain
\begin{eqnarray*}
\C P^3(3,2,4,5) &=& \C^3 / \Z_5 \,\, \cup \,\, \C^2 / \Z_4 \,\, \cup \,\, \C \,\, \cup \,\, \{ p\},
\end{eqnarray*}
where $\{ p\}$ denotes the point $(X_2, Z_2, \varepsilon _2) = (2,0,0)$, at which the weighted blow-up
was performed.
In local coordinates, $\C^3 / \Z_5$ is the $(x, y, z)$-space divided by the $\Z_5$ action, and
$\C^2 / \Z_4$ is the set $\{ (X_3, Y_3, \varepsilon _3) \, | \, \varepsilon _3 = 0\}$ divided by the $\Z_4$ action.
On the set $\{ (X_3, Y_3, \varepsilon _3) \, | \, \varepsilon _3 = 0\}$, the foliation is defined by the Hamiltonian system (\ref{Bou1})
with the Hamiltonian function $H = 2X^2 - 8Y^3 - 4Y$.
This equation is actually invariant under the $\Z_4$ action given by (\ref{3-4}).

Next, due to the definition of the weighted blow-up, we have
\begin{eqnarray*}
& & \mathcal{M}_I 
= \C^3 / \Z_5 \,\, \cup \C^2 / \Z_4 \,\, \cup \,\, \C \,\, \cup \,\, \C P^2(6,4,5).
\end{eqnarray*}
Since $\C P^2(6,4,5) = \C^2 / \Z_5 \cup \C P^1(6,4)$, we obtain
\begin{eqnarray}
 = & &  \,\, \C^3 / \Z_5  \quad \cup \quad\quad \C^2 / \Z_4 \quad\quad\, \cup \quad \C \nonumber \\
&\cup & \,\, \C^2 / \Z_5 \quad \cup \,\, \C P^1(6,4)\backslash \{ q\} \,\, \cup \,\,\,\, \{ q\},
\label{7-2}
\end{eqnarray}
where $\{ q\}$ denotes the point $(u_1, v_1, w_1) = (0,0,0)$.
In local coordinates, $\C^2 / \Z_5$ is given as $\{ (u_3, v_3, w_3) \, | \, w_3 = 0\}$ divided by $\Z_5$.
This implies that the first column $\C^3 / \Z_5 \cup \C^2 / \Z_5$ is just the fiber space 
$\mathcal{P}$ for ($\text{P}_\text{I}$) divided by the $\Z_5$ action, the fiber space over $\C$ whose fiber
is the space of initial conditions.
The last column $\C \cup \{ q\}$ is the Riemann sphere.

Let us investigate the second column $\C^2 / \Z_4 \cup \C P^1(6,4)\backslash \{ q\}$.
On the space $\C^2 / \Z_4$, the equation (\ref{Bou1}) divided by the $\Z_4$ action is defined,
and $\C P^1(6,4)$ is attached at ``infinity".
Note that each integral curve $X^2 = 4Y^3 + 2Y - g_3$ of (\ref{Bou1}) defines
an elliptic curve, where $g_3$ is an integral constant; compare with the Weierstrass normal form
$X^2 = 4Y^3 - g_2 Y - g_3$.

The Weierstrass normal form defines a complex torus when $g_2^3 - 27 g_3^2 \neq 0$.
Two complex tori defined by $(g_2, g_3)$ and $(g_2', g_3')$ are isomorphic to one another 
if there is $\lambda \neq 0$ such that $(g_2, g_3) = (\lambda ^4 g'_2, \lambda ^6 g'_3)$.
Hence, $\C P^1(6,4)\backslash \{ \text{one point}\}$ is a moduli space of complex tori.

According to Eq.(\ref{Bou1}), the $(X_3,Y_3)$-space is foliated by a family of elliptic curves
(including two singular curves $g_3 = \pm (8/27)^{1/2}$) defined by the Weierstrass normal form $X^2 = 4Y^3 + 2Y - g_3$.
By the $\Z_4$ action $(X_3, Y_3) \mapsto (iX_3, -Y_3)$, 
the normal form is mapped to $X^2 = 4Y^3 + 2Y + g_3$.
This means that by the $\Z_4$ action induced from the orbifold structure,
two elliptic curves having parameters $(-2, g_3)$ and $(-2, -g_3)$ are identified.
However, the equality $(-2, g_3 ) = (-2 \lambda ^4, g'_3 \lambda ^6)$ holds for some $\lambda $
if and only if $g_3 = g_3'$ or $g'_3 = -g_3$.
This proves that two elliptic curves identified by the $\Z_4$ action are isomorphic with one another,
and $\C^2 / \Z_4$ is foliated by isomorphism classes of elliptic curves
(including a singular one, while the case $g_2 = 0$ is excluded).
The set $\C P^1(6,4)\backslash \{ q\}$ is expressed as $\{ (u_2, v_2, w_2) \, | \, v_2 = w_2 = 0\}$
divided by the $\Z_2$ action $u_2 \mapsto -u_2$.
We can show that each isomorphism class of an elliptic curve $X^2 = 4Y^3 + 2Y \mp g_3$ intersects with the moduli space 
$\C P^1(6,4)\backslash \{ q\}$ at the point $(u_2, v_2, w_2) = (g_3/4, 0,0) \sim (-g_3/4,0,0)$.
This proves that the second column of (\ref{7-2}) gives an elliptic fibration whose base space is
the moduli space $\C P^1(6,4)\backslash \{ q\}$ and fibers are isomorphism classes of elliptic curves
including the singular curve, but excluding the curve of $g_2 = 0$.
\\[0.2cm]
\textbf{Theorem.\thedef.} The total space $\mathcal{M}_I$ is decomposed into the disjoint union of 
the fiber space for ($\text{P}_\text{I}$) divided by $\Z_5$, an elliptic fibration obtained 
from the Weierstrass normal form as above, and $\C P^1$.
\\

Similar results also hold for ($\text{P}_\text{II}$) and ($\text{P}_\text{IV}$), for which elliptic curves are not 
defined by the Weierstrass form but Hamiltonians represented in Table 4.

%%%%%%%%%%%%%%%%%%%%%%%%%%%%%%%%%%%%%%%%%%%%%%%%%%%%%%%%%%%%%%%%%%%%%%%%%%%%%%%%%%%%%%%%%%%%%%%%%%

\subsection{The extended Dynkin diagram}

It is known that an extended Dynkin diagram is associated with each Painlev\'{e} equation (Sakai \cite{Sak}).
For example, the diagram of type $\tilde{E}_8$ is associated with ($\text{P}_\text{I}$).
Okamoto obtained $\tilde{E}_8$ as follows:
In order to construct the space of initial conditions of ($\text{P}_\text{I}$),
he performed blow-ups eight times to a Hirzebruch surface.
After that, vertical leaves, which are the pole divisor of the symplectic form, are removed.
The configuration of irreducible components of the vertical leaves is 
described by the Dynkin diagram of type $\tilde{E}_8$, see Fig.\ref{fig3}(a).
Our purpose is to find the diagram $\tilde{E}_8$ hidden in the total space $\mathcal{M}_I$.

\begin{figure}
\begin{center}
\includegraphics[]{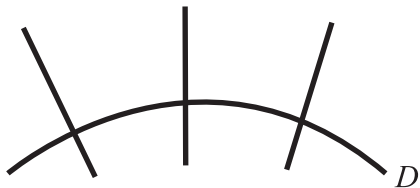}
\caption[]{The minimal resolution of singularities of the closure of $E(0)$.}
\label{fig2}
\end{center}
\end{figure}

\begin{figure}
\begin{center}
\includegraphics[]{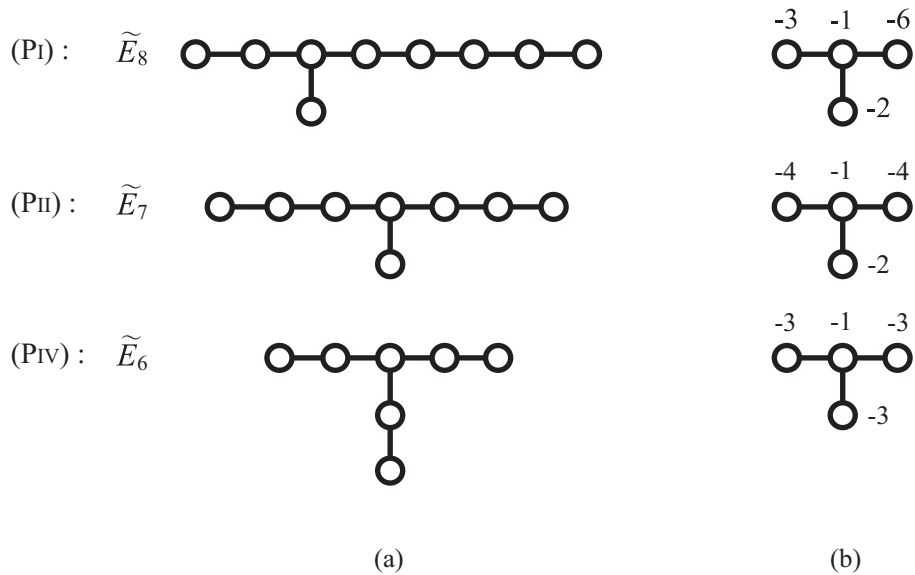}
\caption[]{(a) The usual extended Dynkin diagram.
Each vertex denotes $\C P^1$, and two $\C P^1$ are connected by an edge if they intersect with each other.
All self-intersection numbers are $-2$.
(b) The diagram obtained from our total space $\mathcal{M}_\text{J}$.
Each number denotes the self-intersection number.}
\label{fig3}
\end{center}
\end{figure}

Recall that $\mathcal{M}_I$ is covered by seven local coordinates; 
the inhomogeneous coordinates (\ref{3-1}) of $\C P^3(3,2,4,5)$ and $(u_i, v_i, w_i)$ defined by (\ref{4-3}).
These local coordinates should be divided by the suitable actions due to the orbifold structure.
The actions on $(Y_1, Z_1, \varepsilon _1)$, $(X_2, Z_2, \varepsilon _2)$ and $(X_3,Y_3,\varepsilon _3)$
are listed in Eqs.(\ref{3-2}) to (\ref{3-4}).
The action on $(u_1, v_1, w_1)$ is given by 
\begin{eqnarray}
(u_1, v_1, w_1) \mapsto (\zeta u_1, \zeta^2 v_1, \zeta w_1), \quad \zeta = e^{2 \pi i/6}.
\end{eqnarray}
Each fiber (the space of initial conditions) for ($\text{P}_\text{I}$) is not invariant under the $\Z_5$ action 
(\ref{1-4}) except for the fiber on $z = 0$.
Hence, we consider the closure of the fiber on $z = 0$ in $\mathcal{M}_I$.
The closure is a 2-dim orbifold expressed as
\begin{equation}
N:= \{ (x, y, 0)\} \cup \{(Y_1 ,0, \varepsilon _1)\} \cup \{ (X_2, 0, \varepsilon _2) \}
\cup \{ (u_3, 0, w_3)\} \cup \{ (u_1, 0, w_1)\}.
\end{equation}
$N$ is a compactification of the space of initial conditions $E(0) = \{ (x, y, 0)\} \cup \{ (u_3, 0, w_3)\}$
obtained by attaching a 1-dim space
\begin{equation}
D:= \{ (Y_1,0,0)\} \cup \{(X_2 ,0,0)\, | \, X_2 \neq \pm 2 \} \cup \{ (u_1, 0,0) \}, \quad
X_2 = Y_1^{-3/2} = u_1^6 \mp 2.
\end{equation}
$N$ has three orbifold singularities on $D$ given by
\begin{equation}
(Y_1, \varepsilon _1) = (0,0), \quad (X_2, \varepsilon _2) = (0,0), \quad (u_1, w_1) = (0,0).
\end{equation}
Let us calculate the minimal resolution of these singularities.
For example, the singularity $(X_2, \varepsilon _2) = (0,0)$ is defined by the $\Z_2$ action
$(X_2, \varepsilon _2) \mapsto (-X_2, -\varepsilon _2)$; i.e. this is a $A_1$ singularity,
and it is resolved by the standard one time blow-up.
The self-intersection number of the exceptional divisor is $-2$.
Similarly, singularities $(Y_1, \varepsilon _1) = (0,0)$ and $(u_1, w_1) = (0,0)$ are resolved by 
one time blow-ups, whose self-intersection numbers of the exceptional divisors are $-3$ and $-6$, respectively.
From the minimal resolution of singularities of $N$, we remove the space of initial conditions $E(0)$.
Then, we obtain the union of four  projective lines, whose configuration is described as Fig.\ref{fig3}(b),
see also Fig.\ref{fig2}.

Although the diagram Fig.\ref{fig3}(b) is different from $\tilde{E}_8$, it is remarkable that 
the self-intersection numbers $-2, -3, -6$ are the same as the lengths of arms from the center of $\tilde{E}_8$.

The same results hold for ($\text{P}_\text{II}$) and ($\text{P}_\text{IV}$). 
Let $\mathcal{M}_\text{II}$ be the total space for ($\text{P}_\text{II}$) 
obtained by two points blow-up with the weights $(4,2,3)$
of $\C P^3(2,1,2,3)$ constructed in Sec.4.2.
Consider a fiber $E(0)$ (the space of initial conditions) on $z = 0$ 
and take the closure $N$ of it in $\mathcal{M}_\text{II}$.
From the minimal resolution of $N$ at the orbifold singularities, we remove $E(0)$.
Then, we obtain the union of four  projective lines, whose configuration and self-intersection numbers
are described in Fig.\ref{fig3}(b).
Although the diagram Fig.\ref{fig3}(b) is different from $\tilde{E}_7$, 
the self-intersection numbers $-4, -4, -6$ are the same as the lengths of arms from the center of $\tilde{E}_7$.
A similar result is true for ($\text{P}_\text{IV}$). 

%%%%%%%%%%%%%%%%%%%%%%%%%%%%%%%%%%%%%%%%%%%%%%%%%%%%%%%%%%%%%%%%%%%%%%%%%%%%%%%%%%%%%%%%%%%%%%%%%%
%%%%%%%%%%%%%%%%%%%%%%%%%%%%%%%%%%%%%%%%%%%%%%%%%%%%%%%%%%%%%%%%%%%%%%%%%%%%%%%%%%%%%%%%%%%%%%%%%%

\end{document}